\documentclass[11pt]{ucsd}
\usepackage{amsmath}

\newtheorem{definition}{Definition}
\newcommand{\Bbb}{\bf}
\newcommand{\frak}{\bf}
\newtheorem{lemma}{Lemma}
\newtheorem{sub-lemma}{Sub-lemma}
\newtheorem{theorem}{Theorem}
\newtheorem{corollary}{Corollary}
\newtheorem{proposition}{Proposition}
\newtheorem{result}{Result}

\newtheorem{fact}{Fact}

\newtheorem{observation}{Observation}

\newtheorem{property}{Property}
\newtheorem{note}{Note}

\newtheorem{formula}{Formula}
\newcommand{\qed}{\hfill {\bf q.e.d} }
\input epsf

\title{Random Delaunay Triangulations, the Thurston-Andreev Theorem, and 
Metric Uniformization}
\author{Gregory Leibon}

\begin{document}
\degreeyear{1999}
\degree{Doctor of Philosophy}
\chair{Peter Doyle}
\othermembers{Professor Bruce Driver \\ Professor Kate Okikiolu \\
Professor Zheng-xu He 
\\ Professor Kenneth Intriligator \\ Professor Steven Shapin}
\numberofmembers{6} 
\prevdegrees{B.A. {University California, San Diego} 199? 
\\ M.A. {University of Washington, Seattle}199? } 
\field{Mathematics}
\campus{San Diego}

\begin{frontmatter}
\maketitle
\copyrightpage
\degreeyear{1999}
\approvalpage
\begin{dedication}
\setcounter{page}{4}
%\nul\vfil
{\large 
\begin{center}
To Squirl
\end{center}}
\end{dedication}

\tableofcontents
\listoffigures

\begin{acknowledgements}

I'd like to acknowledge that the story I present here 
would not exist without my thesis  
advisor Peter Doyle, and his inspired view of mathematics.  It was his sense 
of esthetics and understanding that led to the conjectured mathematical stories
for which this thesis is a confirmation.  
Trying to understand Peter's world 
was a huge part of 
my inspiration to become a part of mathematics and I can't thank him enough 
for sharing so many of his beautiful ideas with me over the years.
Thanks Peter!

UCSD was a wonderful community to do mathematics in, largely due to  
the huge number of supportive people, and I thank you all.  I would especially 
like to thank Bruce Driver, Jay Fillmore and Jeff Rabin
for their many years of support and mathematical inspiration.
I would also like to thank  Kate Okikiolu and Zheng-xu He
for the insights they have provided into the story presented here.
I'd also like to thank all my fellow   graduate students  who over the years 
have listened to my rantings and given sound advice 
or a good kick in the butt.     
I'd also like to thank Dave Collingwood and  
my master's thesis advisor John Lee 
at the University of Washington for their support; as well as my 
Washington math buddies 
Steve, Mike, Dave and Siva for showing me 
lots of neat stuff and being great friends.  I'd especially 
like to thank Albert Nijenhuis of University of Washington 
for his wonderful comments 
and his beautiful proof of the small circle intersection theorem  
which appears in chapter 3.  

I would also like to thank my outside committee members Steven Shapin 
and Ken Intriligator (as well as my inside department committee members) 
for putting up with my poor planning abilities and general silliness.  
On this note a 
{\bf HUGE} 
thanks goes out to Lois Stewart for saving 
my butt many many times;  I certainly would have been kicked 
out of UCSD long before completing my thesis  if not for her watchful eye.

The biggest thanks of all goes to my wife Nicole, for always 
being so supportive and  never letting me lose track of what's 
really important in life.

\end{acknowledgements}
\begin{vitapage}
\addcontentsline{toc}{chapter}{\protect\numberline{}Vita and Publications}
\begin{vita}
%\item [January 10, 1971] Born, Los Angeles, California
\item [1993] B.A., \emph{summa cum laude},  
 University of California San Diego 
\item [1993--1995] Teaching assistant, Department of Mathematics, 
University of Washington Seattle 
\item[1996] M.S.,
University of Washington Seattle 
\item [1995--1997] Mathematics Instructor, Department of Mathematics, 
University of Washington Seattle 
\item [1995--1999] Teaching assistant, Department of Mathematics, 
University of California San Diego
\item [1998--1999] Associate Instructor, Department of Mathematics, 
University of California San Diego
\item [1999] Ph.\ D., University of California San Diego 
\end{vita}
\begin{publications}
\item \textsl{The ideal Thurston-Andreev 
		  theorem and triangulation production}.   In Preparation. 
\item 
\textsl{Random Delaunay triangulations and metric uniformization}. In Preparation.
\item 
\textsl{Delaunay triangulation of a surface}. In Preparation.
\end{publications}
\end{vitapage}

\begin{abstract}
In this thesis a connection between the world of discrete 
and continuous conformal geometry is explored.  
The world of discrete conformal 
geometry is related to disk  pattern construction 
and triangulation production.  In particular we discover a generalization of the 
Thurston-Andreev theorem
to any angles in $(0,\pi]$  (when $\chi(M) <0$).. 
The proof of  this theorem
relies on an energy measuring how ``uniform'' the angle data of a 
triangulation is, where by uniform angle data I mean the data 
of a geodesic triangulation of a hyperbolic surface.  The connection to 
continuous geometry is via 
averaging 
this energy over geodesic triangulations, a process which 
forms an 
energy measuring how ``uniform'' (constant curvature) a metric is. 
In fact, the entire discrete energy proof carries over to produce a proof of 
the metric uniformization theorem of conformal geometry.  
The random triangulations used to produce this averaging 
are random Delaunay triangulations, and they are 
explored in some detail.  In particular the averaging techniques developed 
to construct the energy on metrics can be used to explore 
other geometric issues, including the production of a
probabilistic interpretation of the 
determinant of the Laplacian and a new
probabilistic proof 
of the Gauss-Bonnet theorem.

\end{abstract}
\end{frontmatter}
\begin{chapter}{Introduction}

It has been known for a while now that the world of disk patterns is 
intimately connected to the world of geometric  uniformization.  
To orient ourselves to what this 
statement means it is useful  to mention the most famous case.  
The disk pattern is the  Koebe theorem
which guarantees that a specified (reasonably nice)
graph  embedded in the closed unit disk can be realized as the nerve of
a disk pattern, which is in fact unique up to M\"obius equivalence.   
In \cite{Thu} Thurston presented the idea of (nicely) packing   
a bounded simply connected planar domain, recording the nerve of this packing,  and 
then using a suitably normalized  
solution to the Koebe disk pattern  with respect to this nerve in order
to approximate the conformal 
mapping of this region to the disk.  That such a mapping 
exits is our first example of a  geometric uniformization 
problem, the Riemann mapping theorem.  Thurston's procedure has 
been  verified  to work in several cases.  In fact this conjecture has 
been solved in several cases.  The case of using a finer and finer 
hexagonal disk
pattern for the approximation is very well understood; 
it was first shown to uniformly approximate the Riemann mapping in \cite{Su}
and is now seen in   \cite{He} to approximate 
uniformly all of the Riemann mapping's derivatives.

In this thesis another connection between these two stories is developed. 
For now I'll give a brief outline which can be viewed as an 
introduction to the introduction.   Precise statement of the below 
discussion can be found in the elaborated introduction,
sections 1.1 - 1.3.

In this thesis a generalization of a disk pattern theorem often called 
the Thurston-Andreev theorem is presented, which is the solution to
a disk pattern problem 
on a compact surface.  
To be precise, in chapter two of this thesis we prove...

\begin{result}\label{res1}
Given a geodesic triangular decomposition  of a
hyperbolic surface one can associate  combinatorial data: 
the 
topological class  of the triangular decomposition  and to each edge 
the  
angle  between the circumscribing 
circles of the two triangles sharing the edge. 
If the angle data is in $(0,\pi]$ then
this combinatorial  data uniquely determines both the 
hyperbolic structure and the geodesic triangulation.   
Furthermore, if  one starts with a topological 
triangular decomposition of  a compact  
surface with $\chi(M) < 0$, and  angles 
with values in $(0,\pi]$ assigned to each edge satisfying 
certain necessary linear conditions 
(conditions $(n_1)$ -- $(n_4)$ of section \ref{srla}),   then there is  a
hyperbolic surface and geodesic triangular decomposition 
realizing this data. 
\end{result}   

The angles in $(0,\frac{\pi}{2}]$  case is well known, 
see \cite{Th}.   This theorem has several interpretations in terms 
of not only a disk pattern problem, but also as a convex 
hyperbolic polyhedron production theorem, and perhaps 
most importantly (from the uniformization point of view)
as a geodesic triangulation production theorem. 
This triangulation production can be viewed as the second step in 
the proof of the above result.  The first step is a completely 
linear max-flow-min cut-type problem
guaranteeing that given such edge angle data one can produce triangle 
angle data  
which ``formally'' has
the same edge angle data associated to it.  
Given such triangle angle data 
one can measure how uniform  this data is by using an 
energy on the set of possible conformally 
equivalent angle data.  This energy  is larger when 
one is closer to the angle data of a
geodesic triangulation of a hyperbolic surface (so for aesthetic 
reasons I've chosen an energy which one tries to maximize, my apologies to the physicists).  
The proof becomes to guarantee the existence 
of a  unique set of angle data maximizing this energy 
and to show it  corresponds to the 
unique solution of the above disk pattern result.

The geometric  uniformization theorem which this result 
will be
shown intimately related to, is the uniformization for 
geometric surfaces, namely:

\begin{result}\label{res2}
Every Riemannian surface is conformally equivalent to a Riemannian 
surface with constant Gaussian curvature.
\end{result}

This theorem  is very well known.  The heart of what takes place here is 
that in chapter 4 there is a proof of the above uniformization 
theorem which is directly related to 
the above disk pattern  theorem and its proof.  
This relationship is to  first 
(as above) solve the completely linear problem of finding a varying negative 
curvature metric, then we can measure how distant this metric is from a 
constant curvature 
metric by averaging the energy associated to data provided by 
randomly selected 
triangulations.  In the end one finds the 
same arguments used to uniformize the discrete angle data 
can be applied to this averaged 
energy on the space of metrics, 
providing a nice   
proof of the uniformization theorem of surfaces. 
This is not the first energy method proof of this theorem, 
in fact there are many.  One particularly interesting energy 
is the $det(\Delta_g)$ as 
used in \cite{os}.  In fact these random energies are related to the    
$det(\Delta_g)$, and in section 1.3 this relationship is made explicit.

The main tool in articulating the averaging procedure is a detailed 
understanding of random Delaunay 
triangulations.  The details of this understanding form  chapter 
3. I think these triangulations are quite interesting 
independent of their connection to 
metric uniformization, and I hope to convince others 
of this by beginning the 
precise introduction with a 
 probabilistic
proof on the Gauss-Bonnet theorem.

It is worth noting that these ideas may be explored in three-dimensions. 
In fact much of what takes place here can be viewed as 
realization in two dimensions of an idea by 
Peter Doyle to prove that a three manifold accepting 
a varying negative curvature metric also  
accepts a constant sectional curvature metric. 
If proved this result would shed considerable light on Thurston's 
uniformization conjecture.
In section 1.3 I sketch the three dimensional strategy, as related to the 
two dimensional strategy presented in 1.1-1.2.

\begin{section}{A New Proof of the Gauss-Bonnet Theorem}\label{anew}

The  Gauss-Bonnet  formula in its most primitive form  is stated 
for a compact  boundaryless surface 
endowed with a geometry as: 
\[ \frac{1}{2\pi} \int_{M} k dA = \chi (M). \]
The quantity, $k$ ,  being integrated over the surface, $M$,  
is the Gaussian curvature, and on  the 
right-hand side is the Euler Characteristic. 
The Euler characteristic is a topological invariant of the surface
which can be computed from triangulation as
\[ \chi (M) = V - E + F, \] 
where $V$  is the number of vertices,  
$E $  is the number of  edges, and $F$ is the number of  faces
in the triangulation. 
In a somewhat unintelligible nut shell, 
the proof here is accomplished by
randomly triangulating the surface and then noting while the 
Euler characteristic is constant the other (now random variables)
$E$, $F$, and $V$ have expected values which can be computed and 
compared to this constant.  As the density of the randomly distributed 
vertices goes to infinity  one finds these expected values produce the 
Gauss-Bonnet formula, along with a probabilistic interpretation of curvature.

This proof's initial need is an articulation of what a random 
triangulation of a surface is. The first step is to ignore the fact that 
anything random  is going on here and to simply attempt to construct a 
a geodesic triangulation in a fixed metric $g$
from a given set of points,
${\bf p} = \{ p_1, ... ,p_n \}$.
To do this one examines all the triples  and pairs in $\{ p_1, ... ,p_n \}$ 
and decides whether or not to put in a face for a given triple or 
an edge for a given pair; this decision procedure will be relative to a 
certain positive number $\delta$ -- the decision radius.  The procedure 
(with its origin found in Delaunay's ``empty sphere'' method, see \cite{Da}) 
is to formally put in a face for a triple or an edge for a pair if the 
triple or pair lies on a disk of 
radius less $\delta$ which has its interior empty of points 
in $\{ p_1, ... ,p_n \}$.   Clearly every pair on a face forms an edge
(though not necessarily the converse), hence at this moment we have 
an abstract two complex. 
%(a set $K$ of triples, pairs, and 
%singletons of  $\{ p_1, ... ,p_n \}$ such that any triple in $K$ has 
%its pairs in $K$ and every pair has its singletons in $K$). 
Call this 
abstract complex $K_{\{ p_1, ... ,p_n \}}$ and denote a  polyhedral 
realization  as
$|K_{\{ p_1, ... ,p_n \}}|$ (for a realization of it place the points at 
$(0, \dots , 0)$ and $(0,\dots,0 ,1,0, \dots,0)$ in   $\Bbb{R}^{n-1}$ 
and take the 
convex hull of the subsets found in  $K_{\{ p_1, ... ,p_n \}}$).

Our goal now is to realize this polyhedron  inside the surface, 
which will involve restricting $\delta$.  To get going we must recall a 
well known geometric constant associated to the surface, the injectivity 
radius, which is the largest number $i$ such that if $d(p,q) < i$ then 
there is a unique geodesic between $p$ and $q$ of length less than $i$.   
If one assumes that $\delta < \frac{i}{6}$ one can 
(see the proof of  $\ref{dtri}$) produce a continuous map 
\[ R:|K_{\{ p_1, ... ,p_n \}}| \rightarrow M, \]
with the edges parameterizing geodesics in the $g$ metric. 
That our procedure forms a triangulation is  now equivalent to 
$R$ being a homeomorphism; and in such a case  we will call the 
resulting triangulation a Delaunay triangulation.  
Such triangulations typically will look locally like figure \ref{tril}.

\begin{figure*}
\vspace{.01in}
\hspace*{\fill}
\epsfysize = 2in 
\epsfbox{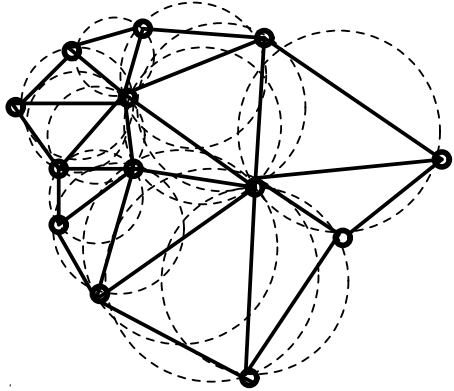}
\hspace*{\fill}
\vspace{.01in}
\caption{\label{tril} Part of a Typical  Delaunay Triangulation}
\end{figure*}

At this point it is useful to find a geometric criterion on a 
set of points ${\{ p_1, ... ,p_n \}}$ guaranteeing that it forms 
a triangulation.  To accomplish this it is useful to 
recall a second geometric 
constant associated to the surface: the strong convexity radius $\tau$. 
On a compact surface it is the largest number with the property that in 
a disk of radius less than $\tau$ the interior of the unique minimal length 
geodesic connecting any pair of points in the disk's closure is in the disk's 
interior.   
%In particular, if a triple lies on the boundary of such a 
%disk we may form the unique geodesic realizing its associated  edges, 
%and a simplex realizing its associated face (lemma $\ref{cones}$).  
We may now introduce our criterion for triangulation detection: we 
will call a set of points  ${\{ p_1, ... ,p_n \}}$ $\delta$-dense 
if each open ball of radius $\delta$ contains at least one $p_i$ and 
${\{ p_1, ... ,p_n \}}$ contains no four of its points is on a circle of 
radius less than $ min \{ \frac{i}{6}, \tau \}$.  
With this concept 
we will prove in section 2.2 the following theorem...
  
\begin{theorem}
\label{dtri}
If ${\{ p_1, ... ,p_n \}}$ is $\delta$-dense with $\delta \leq  min \{ \frac{i}{6}, \tau \}$, then ${\{ p_1, ... ,p_n \}}$ forms a Delaunay triangulation. 
\end{theorem}

 {\bf From here on out always assume the decision radius $\delta$ satisfies $\delta \leq  min \{ \frac{i}{6}, \tau \}$}.

Now comes the second step in this proof: to  put down random sets of points.
The points will be distributed with a density $\lambda$ via a Poisson 
distribution relative to $g$ (developed in the section 2.3), 
and it is convenient to also 
randomly select an ordering of the points after placing them down.  
In technical terms, any such   
point configuration will be denoted {\bf p} and 
will live in a measure space $\frak {P}$ 
with a probability measure $ \Bbb P_{\lambda}$ and specified measurable 
sets $\bf B$ - all of which are constructed in section 2.3.  For now 
all we need to understand are certain basic properties and terminology.  
Let  $U \subset M$ be measurable with area $A(U)$. The first property is 
the characteristic property of the Poisson distribution, namely the set of 
point configurations in  $\frak {P}$ with exactly $n$ points in $U$ is an 
element of  $\bf B$ with measure 
$\frac{{(A(U) \lambda)}^n e^{-\lambda A(U)}}{n!}$.

Recall a random variable on $\frak P$ is simply a 
$\bf B$ measurable function.  We will encounter several in this proof,
 namely $E({\bf p})$, $V({\bf p})$, $F({\bf p})$
, and the characteristic (indicator) 
functions of various measurable sets (denoted $1_{U}$ with $U$ measurable).
The expected value of a random variable $L({\bf p})$ 
is the integral of the function 
over $\frak P$ with respect to the measure $\Bbb P_{\lambda}$, and 
will be denoted $\Bbb E_{\lambda}(L)$.  
%When the choice of metric, $g$, used is to construct the Poisson distribution is vague this expected value will be denoted  $\Bbb E^{g}_{\lambda}(L)$. 
As an example of such an expected value, 
 using $A$ to denote the area on $M$, we have from the above  characteristic 
property
 \[  \Bbb E_{\lambda}(V) = \sum_{n=0}^{\infty} n \frac{ A^n \lambda^n}{n!} e^{- A \lambda}=  A\lambda \sum_{n=1}^{\infty} \frac{ {(A \lambda)}^{n-1}}{(n-1)!} e^{- A \lambda}= A \lambda, \]
which reveals the use of the term density for $\lambda$.   
The most important set for us in $\frak P$ is the set 
$\frak T_{\delta} \subset \frak P$ consisting of configurations which have Delaunay triangulations associated to them.   By theorem $\ref{dtri}$, when our points are put in with a high density we expect they will typically form Delaunay triangulations.  Allowing $ O(\lambda^{-\infty}) $ to mean a quantity decaying faster than any polynomial in $\lambda$, in section 2.3 we prove this with:
  
\begin{theorem}
\label{size} 
If $L$ is any one of  $1_{\frak T_{\delta}}$, $E$, $V$, or $F$, then $L$ is measurable and has expected value 
$ \Bbb E_{\lambda}  (L) = \Bbb E_{\lambda} (L 1_{\frak T_{\delta}} ) + O(\lambda^{-\infty}) $.   
\end{theorem}

Now one observes that the Euler characteristic of any configuration in $\frak T_\delta $ is the Euler characteristic of an actual triangulation, hence the  constant $\chi(M)$; so along with the above theorem $\ref{size}$ one has 
\[\chi (M) =  \Bbb E_{\lambda}(\chi(M) 1_{\frak T_{\delta}}) +   O(\lambda^{-\infty}) =  \Bbb E_\lambda ( (V - E + F) 1_{\frak T_{\delta}})  +  O(\lambda^{-\infty}). \]
Applying the mathematical triviality, yet philosophical miracle, that expected values add gives us

\[ \chi (M) = \Bbb E_\lambda (V 1_{\frak T_{\delta}}) -   \Bbb E_\lambda (E 1_{\frak T_{\delta}})  + \Bbb E_\lambda (F 1_{\frak T_{\delta}})  +  O(\lambda^{-\infty})  . \]

In an actual  triangulation we have $\frac{3}{2}E = F$ so
\[\chi (M) =  \Bbb E_\lambda ( V 1_{\frak T_{\delta}}) - \frac{3}{2} \Bbb E_\lambda (F 1_{\frak T_{\delta}}) +  \Bbb E_\lambda (F 1_{\frak T_{\delta}})  +  O(\lambda^{-\infty}). \]

As explored  above, we have the expected number of vertices is $\lambda A $; this along with theorem $\ref{size}$ gives us 

\[\chi (M)  = A \lambda  - \frac{1}{2}  \Bbb E_\lambda (F 1_{\frak T_{\delta}})  +  O(\lambda^{-\infty}) . \]

So finally using theorem $\ref{size}$ one last time we find: 
 
\begin{formula}[Euler-Delaunay-Poisson Formula]
\label{edpf}
\[  \chi (M) = \lim_{\lambda \rightarrow \infty} 
\left( A \lambda  - \frac{1}{2}  \Bbb E_\lambda (F)  \right). \]
\end{formula}

The goal now becomes to compute $\Bbb E_\lambda (F)$. 
Now a triple of points occurring in a Poisson configuration, 
$y \in M \times M \times M$, 
forms a face exactly when one of its associated disks
is empty of points.  
The first thing we'd like is that we are in fact 
talking about a unique such disk.  This is the most sophisticated 
fact needed, and occurs throughout 
the proof. I'll record it here for future reference; for
a couple of proofs see \cite{Le}:  

\begin{lemma}[The Small Circle Intersection Lemma]
\label{lit}
If a triple of points lies on a the boundary of a disk with 
radius less than $\delta$, then this disk is unique among disks of 
radius less than  $\delta$. (Recall $\delta \leq min\{\frac{i}{6},\tau\}$.)
\end{lemma}

Back to the computation, by the characteristic property of the 
Poisson distribution the probability any area $a$ disk is empty is 
$e^{-\lambda a}$.  So letting $V_{\delta} \subset \times^3 M$  
be the set of ordered triples living on circles of radius less than 
$\delta$, we now see the probability that $y$ has an associated face should be 
$\frac{1}{A^3} \int_{V_{\delta}} e^{-\lambda a(y)} dA^3 $, 
with $a(y)$ the area of $y$'s  uniquely associated disk.   
This together with the fact that the expected number of triples is  
$ \sum_{n=0}^{\infty} \binom{n}{3} 
\frac{{(A \lambda)}^n e^{-\lambda A }}{n!} = \frac{{(A \lambda)}^3}{6} $
leads one to hope, as is confirmed in the beginning of section 3, that 

\begin{eqnarray}
\label{ex1}
 \Bbb E_{\lambda}(F) =    \frac{{ \lambda}^3}{6} \int_{V_{\delta}} e^{-\lambda a(y)} dA^3. 
\end{eqnarray}

To compute this explicitly it is necessary to put coordinates 
on $V_{\delta}$.  First one chooses a way to discuss directions 
at all but a finite number of tangent planes of $M$ (via a orthonormal frame).
Then one can parameterize a full measure subset of 
$V_{\delta}$ with a  subset
of  
$(\theta_1, \theta_2,\theta_3, r,p) \in S^1 \times S^1 \times  S^1 
\times (0 , \delta) \times M$ 
by hitting the triple described by sitting at the point $p \in M$  
and moving a distance $r$ in each of the three directions $ \theta_1$, 
$\theta_2$, and $\theta_3$ .
 That this is a legal parameterization is easily deduced from 
lemma \ref{lit}, and discussed in section \ref{coords}.
Note that when fixing $p$ and varying 
$\theta_i$ in these coordinates one 
produces the Jacobi field $J_i(r)$; whose norm I will call  $j_i$.
% measures how quickly the the geodesics spitting out of $p$ are separating.
Using this notation,
letting   $d \vec{\theta} =  d \theta_1 \wedge d \theta_2 \wedge d \theta_3$,
and 
letting  $ \nu( \vec{\theta})$ be the area of the triangle 
in the Euclidean unit circle 
with vertices at the points corresponding to the $\{ \theta_i\}$,
% (i.e. \left|\frac{1}{4} \sin\left(\frac{\theta_2 - \theta_1}{2}\right)  \sin\%left(\frac{\theta_3 - \theta_2}{2}\right) \sin\left(\frac{\theta_3 - \theta_1}%{2}\right)\right| $,
we will see in section 3.1  that equation (1)  expressed in these 
new coordinates is

\begin{eqnarray}
\label{ex2}
 \frac{\lambda^3}{6}  \int_{M} \int_0^{\delta} \int_{\times^3 S^1} e^{-\lambda a(\vec{\theta},r,p)} j_{\theta_1}j_{\theta_2}j_{\theta_3}  \nu( \vec{\theta}) d \vec{\theta} dr dA.
\end{eqnarray}

Both $j_i$ and the area of a ball function, $a(\vec{\theta},r,p)$,
are directly related to the curvature; 
%(reviewed in the next subsection)
with this  relationship equation \ref{ex2} can be expressed as
\begin{eqnarray}
\label{ex3}
 \frac{\lambda^3}{6}  \int_{M} \int_0^{\delta} \int_{\times^3 S^1}  
e^{-\lambda \pi r^2}\nu(\vec{\theta})  \left(  r^3 -\frac{ k}{2}  r^5 + \frac{\pi\lambda k}{12} r^7 +  ``O"(r^6) \right) d \vec{\theta} dr dA,
\end{eqnarray}
with the  $``O"(r^6)$ a  well controlled
function of $\vec{\theta}$, $r$, $p$ and $\lambda$. 
In fact the $``O"(r^6)$ term is so controlled that 
there is a $\rho_M >0$ associated to the surface such that 
if $\delta < \rho_M$  after integrating one can reduce equation \ref{ex3} to

\begin{eqnarray}
\label{ex4}
  \Bbb E_{\lambda} (F) =  2 A \lambda -    
\frac{1}{ \pi } \int_{M} k dA  + O(\lambda^{-\frac{1}{2}}), 
\end{eqnarray}
with the $|O(\lambda^{-\frac{1}{2}})| \leq  C  \lambda^{-\frac{1}{2}}$
(see section \ref{dercomp}  for the details concerning \ref{ex3} and \ref{ex4}
and the 
related estimates). 

In particular we may now plug this into the Euler-Delaunay-Poisson Formula 
(formula \ref{edpf}) to give simultaneously  a probabilistic 
interpretation of curvature and new proof of the Gauss-Bonnet theorem.

\begin{formula}[Euler-Gauss-Bonnet-Delaunay Formula]
\label{egbdf}
Using a decision radius $\delta < \rho_M $ we have 
\[\chi (M) =  \lim_{\lambda \rightarrow \infty} ( A \lambda  - \frac{1}{2}  \Bbb E_{\lambda} (F)) = \frac{1}{2 \pi} \int_{M} k dA .\]
\end{formula}

Note this allows us to interpret the Gaussian curvature as the 
defect in the expected number of faces in a random Delaunay triangulation
in the surface's geometry 
from what would be expected in 
Euclidean space.

\end{section}

\begin{section}{Conformal Geometry and Uniformization}\label{conu}

The techniques used in the previous section
  to prove the Gauss-Bonnet  theorem can be used to examine certain 
aspects of conformal geometry.  
The ideas will involve the comparison of conformally equivalent metrics. 
 In terms of Riemannian metrics we say two 
 metrics on $M$, $g$ and $h$, are conformally 
 equivalent if  $h = e^{2 \phi} g$ for a 
 smooth function $\phi$.  We will always be thinking in terms of a fixed 
background metric
called $g$ and will label its associated geometric  objects like its gradient,
 Laplacian, curvature, norm, area element, or area  
as $\nabla$,  $\Delta$, $k$,$|\cdot|$ , $dA$ or $A$. 
 For the $h=e^{2 \phi} g$ metric we shall denote these object with an $h$
subscript.  
%Similarly for any of the  metric dependent objects  constructed in the previous section.  

The key to using the ideas of the previous section 
is the development of
a discrete analog of a conformal structure and a conformal 
transformation relative to a Delaunay triangulation.   
The most important observations in guessing this  discrete structure 
concerns two  triangles
 $t_1$ and $t_2$ in a Delaunay triangulation 
sharing an  edge $e$, and the intersection 
angle  $\theta^{e}$ between  the circles in
 which  $t_1$ and $t_2$  are inscribed (see figure \ref{zip}).
 The facts of use to us are  that 
 $\theta^e \in (0,\pi)$ and $\theta^e$ is preserved under a
conformal deformation  of the metric up to order $O(r)$
(see formula \ref{wind} in section \ref{triangulardef}).

\begin{figure*}
\vspace{.01in}
\hspace*{\fill}
\epsfysize = 1.5in 
\epsfbox{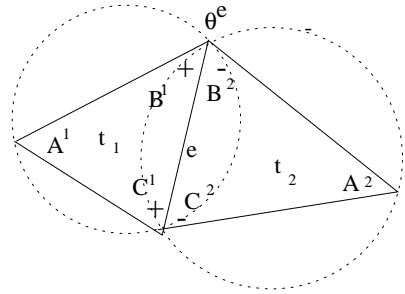}
\hspace*{\fill}
\vspace{.01in}
\caption{\label{zip} The Notation of Neighboring Triangles}
\end{figure*}

To use this information it is necessary to note 
$\theta^e$ can be written down 
accurately in terms of the angles within the triangles.    

\begin{formula} \label{theta} 
Letting $A^i$, $B^i$ and $C^i$ be the angles
in  $t_i$ as in  figure \ref{zip} we have
    \[\theta^e= \frac{\pi+A^{1}-B^{1}-C^{1}}{2} 
    + \frac{\pi+A^{2}-B^{2}-C^{2}}{2} + O(r^3).\]
    \end{formula}
    
It is worth noting this formula is exactly true on a surface with
constant curvature.  One proof of this formula in the negative curvature case  
%(see \cite{Le} for the proof in this and more general cases) 
utilizes an object, an ideal hyperbolic prism, that will turn out to be 
fundamental to all that takes place here; and it is useful to 
construct and discuss it now.   That it shows up at such a fundamental 
point in understanding the triangulations will hopefully inspire its otherwise
unexpected appearance in the next section. 

The prism  is constructed from the angle data of a 
hyperbolic triangle, namely a set of positive angles 
$\{A,B,C\}$ such that $A + B + C -\pi <0$.
To construct it first form a hyperbolic 
triangle with the $\{A, B, C\}$ data,  then  place the triangle on
 a hyperbolic plane and 
      then place the plane in hyperbolic three space.  Now 
union this triangle with   the geodesics 
       perpendicular to this 2-plane going through the 
vertices of the triangle.  The prism of interest is
 the convex hull of this 
arrangement,   
See figure \ref{pri1}.

\begin{figure*}
\vspace{.01in}
\hspace*{\fill}
\epsfysize = 1.5in 
\epsfbox{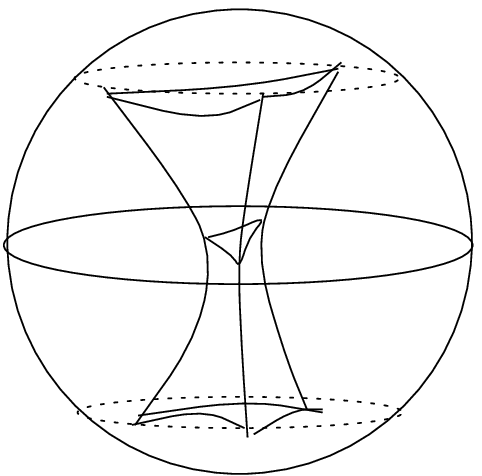}
\hspace*{\fill}
\vspace{.01in}
\caption{\label{pri1} The Ideal Prism}
\end{figure*}

With these  conformal geometry 
facts in mind we may set up our discrete conformal geometry.

\begin{subsection}{Discrete Uniformization}\label{dsc}

  The discrete object replacing a 
metric is the information naturally associated to 
a Delaunay geodesic triangulation living on the surface in this metric;
namely a topological triangulation, $\frak{T}$, and the angles 
    in all the triangles.
Given $\frak{T}$ there are 
$3F$ slots $\{ \alpha_i \}$ in which one can insert 
possible triangle angles, which we 
will place an order on and identify with a basis of a 
$3 F$ dimensional real vector space.  
With this basis choice we will 
denote this vector space as   
$\Bbb R^{3F}$, and  denote vectors in it as $x =\sum A^i \alpha_i$. 
Let 
$\alpha^i$ be a dual vector such that 
$\alpha^i(\alpha_j) = \delta_i^j$.  Sometimes we will refer to the 
angles in a fixed triangle, and will abuse 
notation in order to refer to figure \ref{zip} 
by letting  $\alpha^i(x) = A^i$, $\beta^i(x) =B^i$, and $\gamma^i(x) =C^i$. 

We will require that these angles satisfy certain requirements 
that the 
angles in an actual Delaunay triangulation would satisfy.
In particular we require that 
all the angles 
are  in $(0,\pi)$ and 
that the sum of all the angle at a vertex is $2 \pi$.  
%We call angles satisfying these condition an angle system.
 
Motivated from the previous section we define the intersection angle
at an edge to be the function  
\[  \theta^e(x) = \frac{\pi+A^1 -B^1 - C^1}{2} +
       \frac{\pi+A^{2}-B^{2}-C^{2}}{2},  \]
%where the angle values correspond the the angle slots 
with the notation inspired from  figure \ref{zip}.
The Delaunay quality is reflected in the fact that the 
``intersection angles'' must satisfy  $\theta^e(x) \in (0,\pi)$ 
for each edge $e$.  It is worth noting that as 
an immediate consequence of the fact that the angles 
sum to $2 \pi$ at a vertex we have that the 
angle discrepancies at a vertex  $\pi -\theta^e(x)$ 
sums to $2 \pi$ at a vertex as well,
see section \ref{srla}.
%(see \cite{Le1}). 

We will call angles satisfying 
the above requirements a Delaunay angle system.
With the previous section as  
motivation, we define a pair of  
Delaunay angle systems  $x$ and $y$ to be  
conformally equivalent if 
they arise form the same triangulation and if for each edge $e$ the 
$\theta^e(x) = \theta^e(y)$.

Motivated by the Gauss-Bonnet theorem 
we define the curvature of a triangle $t$ relative to the Delaunay 
angle system $x$ 
to be  
$k^t(x) = A + B + C - \pi$,  
where $\{A,B,C\}$ is the angle data corresponding to angle slots 
in $t$.  On 
a surface with $\chi(M) <0$  there 
is a very simple set of linear equations 
which will gaurentee that a 
Delaunay angle system is conformally equivalent to a Delaunay 
angle system where each triangle has negative curvature.  

Let $S$ be a set of triangles in $\frak{P}$, and denote the 
cardinality of $S$ as $|S|$.  Let condition ($\star$) on 
a set of angles $\{\theta^e(x)\}$ associated to all
the edges of $\frak{P}$ be the condition that  for any set of 
triangles $S$ we have 

\hspace{.5in} $(\star)$ \hspace{1in} 
$ \sum_{e \in S} \theta^e(x) > \pi |S|.$

As an immediate  consequence of theorem 
\ref{cir1}
% in  \cite{Le1} 
we have the following lemma.

\begin{lemma}[The Discrete Teleportation Lemma]\label{dtele}
A Delaunay angle system $x$ has  $\{\theta^e(x)\}$ satisfying $(\star)$ 
if and only if $x$ is conformally 
equivalent to a negative curvature  Delaunay angle system.
\end{lemma}   

We will only be concerned with Delaunay  angle systems satisfying 
the the conditions of this lemma, which note 
includes the angle data associated to a 
Delaunay triangulation of surface with varying negative curvature.
Denote the convex bounded set of negative curvature angle 
systems conformal to a given 
one, $x$, as 
$\frak{N}_{x}$.

%The negative curvature quality is simply the restriction that 
%$k^t(x) = A + B + C - \pi < 0$ 
%for all $t$,
%were $\{A,B,C\}$ is the angle data corresponding to angle slots in $t$. 
%A topological triangulation
%along with formal angle 
%data satisfying these conditions is what I call a negative curvature 
%Delaunay angle system, which I will refer to simply as a 
%{\bf angle system} for the remainder of this paper (see \cite{Le1} 
%for cases where the above conditions are relaxed). 
%Fixing a triangulation with $3F$ angles, 
%the set  of systems
% is clearly a convex subset of $(0,\pi)^{\bf A}$ - and will be viewed as 
%such throughout. 

The  pleasure derived from the Delaunay angle systems comes from  
 a beautiful Energy which lives on  $\frak{N}_{x}$.  Let $V^t(x)$ 
denote the volume of the ideal hyperbolic 
prism constructed from $t$'s angle data  relative to $x$,
as described in the previous section.
Now simply let the energy be
\[E(x) = \sum_{t \in \frak{P}} V_t(x). \]

The points where $E$ attains a  maximum will be of great interest to us.  
It is worth noting that  $E$  is a continuous function on the 
compact set $\bar{\frak{N}}_x$,  
so attains its  maximum.  
In fact it is an immediate consequence of 
lemma \ref{push} 
%in \cite{Le}
in section \ref{prot}
that the angle system with 
maximal energy is in fact in $\frak{N}$. 
An understanding of this as well as the energies  behavior at critical points 
comes from understanding  its
differential.  To interpret it let  $a$, $b$ and $c$ denote the
edge lengths opposite to the angle  $A$, $B$, and $C$ in the hyperbolic 
triangle determined by the angle data $\{A,B,C\}$ 
(which once again exists since $k^t(x) < 0$).  
In formula \ref{diff} of section \ref{prot} 
% of \cite{Le1}  
the following formula is produced:
\begin{formula} \label{dif}
$dE^x = \sum_{\alpha_i \in \frak{T}} E_i(x) \alpha^i$ with 
 \[ E_i(x) = \frac{1}{2} \left( 
   \ln\left(\frac{\cosh(a)-1}{2}\right) - \ln\left(\frac{\cosh(b)-1}{2}\right)
   -\ln\left(\frac{\cosh(c)-1}{2}\right) \right), \]
where  $\alpha^i(x)=A$, $\alpha_i \in t$, $t$'s angle data is
$\{A,B,C\}$,  and the $a$,$b$, and $c$ are determined as above.
 \end{formula}

 In order to best exploit this formula it is necessary to note 
(see the proof of corollary \ref{idpun}  
in  section \ref{prot})
%in \cite{Le1} 
that the 
tangent space at any point of
 $\frak{N}_{x}$ 
can be described 
explicitly as the translation to that point of 
the span over all edges of the vectors 
\[ w_e = \gamma_1 + \beta_1 - \left( \gamma_2 + \beta_2\right) \]
 as  in  figure $\ref{zip}$.
%where $\beta^i(x) = B_i$ and $\gamma^i(x) = C_i$.  
%\[ V_l = \frac{\partial}{\partial B_{1}} + 
% \frac{\partial }{\partial C_{1}} - \left(\frac{\partial}{\partial B_{2}} + 
% \frac{\partial }{\partial C_{2}} \right), \]
The  loveliness of this energy
can now be expressed in terms the following  observation about its 
critical points (i.e. where  $dE = 0$).

\begin{observation}\label{dob}
From the above formula and the above description of 
$T_p(\frak{N}_{x})$
we see  at a critical point $x$ of the energy in a conformal class satisfies 
\[0 =  dE(w_e) = \ln\left(\frac{\cosh(a_{2})-1}{2}\right) - 
\ln\left(\frac{\cosh(a_{1})-1}{2}\right),\]
for all edges $e$ with the notation coming form  as usual from 
figure \ref{zip}. 
%the hyperbolic triangle edge length opposite to $\alpha_i$.
Hence  at a critical point of $E$ we have that $a_{1} = a_2$
and  the set of $- 1$ curvature triangles  formed from the given
angle data
fit together to form an actual
constant curvature surface.  
      \end{observation}
   
I will call the Delaunay angle system  of such a critical point
   a uniform angle  system.  The question becomes: how many (if any) 
uniform structures
can be  associated to a given   
angle system? From above we 
 know there is at least one internal maximum  and in fact
$E$  is strictly concave down  
(see lemma \ref{lemer} in 
\ref{prot}). 
So any critical point is a maximum and 
unique, as needed.

\begin{theorem}[Discrete Uniformization Theorem]\label{pack}
If $\chi(M) < 0$ and $x$ is a Delaunay angle system 
where  $\{\theta^e(x)\}$ satisfies $(\star)$ then $x$ is 
 conformally equivalent to a unique uniform angle system.   
\end{theorem} 

%In fact one can flow an initial set of angles to the uniform one 
%using the gradient of the energy. To articulate we need to  
%associate the Gaussian curvature to a triangle (motivated by the Gauss-Bonnet th%eorem) as
%$k(t) = \frac{K(t)}{area(t)}$ In fact it is useful to use the discrete 
%$$H_2$ inner product (see section?). In this inner-product we have 
%\[ grad(E) = \sum_{t \in T} \ln(k)  .\]
%The natural log of this curvature naturally slits up into 

%Note if we start with a negative curvature angle system then we could use the gradient flow of $E$ with respect the metric induced from the standard metric$\Bbb R^{3f}$ to flow the given angle data to the uniform one

It is worth noting that 
there is a rephrasing of the above uniformization 
construction in terms of
a solution to a disk pattern problem, 
namely as a corollary of theorem \ref{cir1} in  
section \ref{srla}
%\cite{le1} we have:

\begin{theorem}\label{ddisk1}
 If one is given a topological triangulation of a surface with $\chi(M) < 0$
and a set $\{ \theta^e \}$ associated to the edges satisfying  
both condition $(\star)$ 
and that at each vertex 
$\sum_{\{e \in v\}}\pi - \theta^e = 2 \pi$,
then
there is a uniquely associated constant curvature surface on which 
the given triangulation is realized as a geodesic triangulation
and  the circles in which the triangles live meet with the specified 
$\theta^e$  angles. 
\end{theorem}

%See \cite{Le1} for generalizations and proofs of this theorem.  

 Solutions to problems similar to this are well known. The first solution to such a pattern problem goes back to Koebe and pattern problems closely related to the above corollary were implicit in the work of Andreev \cite{An} and then rediscovered and articulated in this language by Thurston (see \cite{Th}).The above theorem (actually a generalization of it found in \ref{srla})
%\cite{le1})
is a generalization to arbitrary angles of the convex ideal case of the Thurston-Andreev theorem.  The first solutions to these pattern problems using energy methods (as far as I'm aware) can be found in \cite{Co}, and 
energy methods using hyperbolic volumes (in the toroidal case)
have there earliest versions in \cite{Be} and later in \cite{Ri2}.

\end{subsection}

\begin{subsection}{Continuous Uniformization}

Now we'd like to mimic the uniformization procedure in the previous section 
for metrics. 
For metrics by uniform 
 structure I will mean 
a metric $e^{2 \phi}g$ with constant curvature.   It is useful to note that 
$k_h = e^{2 \phi}(-\Delta \phi + k)$.

For starters let us 
note in the metric world we still have 
\begin{lemma}[The Metric Teleportation Lemma]\label{mtele}
Every metric on a surface 
is conformally equivalent to a metric with  either 
negative, positive , or zero curvature.
\end{lemma}
{\bf Proof:}
  This is 
 a consequence of 
the Fredholm alternative, which says  if $f \in C^{\infty}(M)$ satisfies 
$\int_M f dA =0$ then  
$f \in \Delta(C^{\infty}(M))$. Since $
\int_M (k - 2 \pi \chi(M)) dA =0$ we have there is a smooth 
$\phi$ satisfying $\Delta \phi  = k - 2 \pi \chi(M) $
and hence 
$k_\phi = e^{2 \phi}(-\Delta \phi + k) = 2 \pi \chi(M)  e^{2 \phi}$ 
as required.   

\qed

With this observation 
in our $\chi(M) <0$ world we will
restrict our attention to metrics with 
strictly negative curvature.
It is worth noting that as in the discrete teleportation lemma 
the problem is
linear.

Now given a varying negative 
curvature metric $h$ if we chose a ${\bf p} = \{p_1, \dots ,p_n\}$ 
and topological triangulation  
$K_{\bf p}$ with ${\bf p}$ as vertices then  we 
could  measure how close to uniform $h$ is  
with the energy of the previous section. Namely  we could let 
\[ {\bf E}_{h}({\bf p}) = \sum_{t \in K_{\bf p}} V^h_t({\bf p}),\]
where $V^h_t({\bf p})$
is computed using 
the angle data associate to the 
``triangulation'' viewed in the $h$ metric.
This of course  means  connecting the  needed vertices of 
${\bf p}$ with $h$ geodesics 
and measuring the resulting $h$ angles.
Notice this makes sense for sufficient $h$ dense points 
since  
by the Gauss-Bonnet theorem for geodesic triangles 
the angles in a triangle  $\{A,B,C\}$ 
will indeed satisfy
\[  k^t({\bf p}) = A + B + C  - \pi <0 . \]

This gives us a measurement of uniformity relative to $\{\bf p\}$ and 
$K_{\bf p}$ and to rid this dependency it is natural to 
average
this measurement over all complexes as in section \ref{anew}.
To do so weight the point distribution and 
decide how to form the complexes 
using a fixed metric $g$  and density 
$\lambda$, and then compute
$\Bbb E_{\lambda}( {\bf E}_{h})$. 
It's worth noting that it may be necessary to shrink the decision 
radius down a bit to make sense out of this construction, since we'd like
points in a triple to be within the injectivity radius of each 
other in both metrics.
As with the random variables of section 
\ref{anew}  this 
random variable is naturally 
expressed as  the sum over all the faces $t \in  K_{\bf p} $
of a  function dependent only on the
data associated to an individual face, so 
we would expect the computations of \ref{anew} to  go through.

It is clear that the higher the vertex  density the  higher the  percentage 
of triangulations involved, hence  the more 
sensitive this energy should be to measuring the uniformity of the curvature. 
With this as inspiration we proceed as with
as with  our Euler characteristic computation and 
take the limit as the density of vertices goes to infinity. Let 
\[    E(h) = \lim_{\lambda \rightarrow \infty} \Bbb E_{\lambda}
({\bf  E}_{h} - {\bf E}_g), \]
and note the second term ${\bf E}_g$ is independent of $h$ and needed only 
to normalize the computation.  
We will show 
that among $h= e^{2 \phi} g$ with  with $k_h < 0$ this formula  
is equivalent to 
\begin{eqnarray}\label{cen}
 E (h)  = -  \int_M || \nabla \phi ||^2 + (\Delta \phi - k)
 \ln(\Delta \phi - k) + k \ln|k| dA .
\end{eqnarray}

The energy is clearly scale invariant and 
$k_{h} = e^{-2 \pi}(-\Delta \phi +k) $, so  we see that it 
is natural to consider the energy as a function 
on the ``energy norm'' closure of   
\[V =  \{ \phi \in C^{\infty} \mid 
 \int_M \phi dA =0 \mbox{ and } -\Delta \phi + k < 0\}. \]
Using the energy norm allows us to  use 
compactness arguments to gaurentee the existence of a maximum, just as 
in the discrete case (for the details of everything that takes 
place in this section see section \ref{bigproof}). 
%\cite{Le}).  
Also as in the discrete we can easily force the maximum into the
 the interior (where $esssup(-\Delta +k) < 0$).
This requires looking at the energy's 
 Frech\'et  differential, which is quite   revealing. 

\begin{formula}\label{cdif}
The Frech\'et derivative of $E$ at $\phi$ in the direction $\psi$ is 
 \[ D E (\psi) 
=   - \int_{M} \Delta \psi \ln|k_{e^{2 \phi}g }| dA  .\] 
 \end{formula}

%the curvature in the $e^{2 \phi} g$ metric we have,
%At a  metric $h = e^{2 \phi}g$  with $k_{h}<0$, we have for all $\psi \in C^{\infty}(M)$ 

This formula for  Frech\'et derivative guarantees that 
$k_{h}$ is constant and $\phi$ is smooth. 
%(see \cite{Le}).  
The easiest way to understand  this is to recall  
the Fredholm Alternative described 
in the poof of lemma \ref{mtele}
and observe...

\begin{observation}\label{cob}
The above formula insures us that if we knew the critical 
$\phi$ to be a $C^{\infty}$ function then
$\ln|k_{h}|$ would  $L^{2}$ orthogonal to the image of $C^{\infty}(M)$ under 
 the Laplacian, so by the Fredholm alternative $\ln|k_{h}|$ is constant.
Hence we indeed see that a critical point of
 $E$ should be a metric of constant negative curvature.
\end{observation}

Just as in the discrete case at this point concavity comes to the rescue 
to gaurentee  uniqueness. Since everything is $C^{\infty}$ 
we really can use the fact that 
the Frech\'et Hessian at 
$\phi \in V$ applied to 
$ (\psi,\psi)$
\begin{eqnarray}\label{hess} 
D^2{\bf  }E (\psi,\psi) =-  \int_{M} |\nabla \psi|^2  + 
\frac{( \Delta \psi)^2}{\Delta \phi - k} dA_0  
\end{eqnarray}
 is strictly negative 
at a non-zero $\psi$  to gaurentee the critical point 
is in fact unique.

So just as in the discrete case we arrive at the 
uniformization theorem (although only when $\chi(M) < 0$).

\begin{theorem}[The Metric Uniformization Theorem]\label{cuni}
 Every metric is conformally equivalent to a unique metric 
 of constant curvature.
  \end{theorem}

\end{subsection}

\begin{subsection}{The Determinant of the Laplacian and Entropy}

The use of energies to solve the uniformization 
problem is not at all new, and 
it is nice to relate this process to some other more familiar energies.
The method of producing other energies is to 
average over metrics using 
the varying metric to distribute the points. So in this section 
I will superscript the expected values with the metric 
used to construct the point distribution.  

%It's worth noting with these strange distribution choices that we ``get lucky'' to produce energies capable of uniformizing metrics, and lose 

The first variant  I will describe has the 
advantage of not needing to be restricted to a conformal class. 
It involves a 
function on the space of metrics I call the entropy
\[ H(h) = \int_M k_h \ln|h_h| dA_h .\]
The fact is that 
\[    E_1(h) = \lim_{\lambda \rightarrow \infty} (\Bbb E^h_{\lambda}
({\bf  E}_{h}) - \Bbb E^g_{\lambda}({\bf E}_g)) = H(h) - H(g). \] 
So viewing $g$ as fixed we have an energy on all metrics - 
which can be checked to optimize at a constant curvature metric 
and be strictly concave down 
in a conformal class.

More interestingly  is the following energy.

\begin{theorem} \label{lapy}
There is a constant $C$ such that on the set of all metrics $h$ conformally 
equivalent to $g$ of a fixed area
\[ E_2(h) = E_1(h) - E(h) 
 =   \lim_{\lambda \rightarrow \infty} (\Bbb E^h_{\lambda}
({\bf  E}_{h}) - \Bbb E^g_{\lambda}({\bf E}_h)) 
= \ln(\det(\Delta_h)) + C. \]
\end{theorem}

As we shall see in section \ref{encomp}, 
the key to being able to make this interpretation is the 
beautiful integral formulation of the $\ln(\det(\Delta_g))$ due to Polyakov
(see \cite{Po}).  The 
$ \ln((det(\Delta_g))$ 
was found to have the uniformization property by  
 Osgood, Phillips, and Sarnak 
(see  \cite{os}). 
In \cite{Co} Yves Colin de Verdi\'ere suggested that 
energies related to circle pattern problems might be related to 
the determinant of the Laplacian.  This procedure  
provides such a relationship.

\end{subsection}

\end{section}
\begin{section}{Three-dimensional Dreams}
The entire two dimensional story has a three dimensional wishful 
analog.  

\begin{subsubsection}{Discrete Uniformization}
The discrete case is essentially an idea do to Casson for deforming a 
triangulated three manifold into one of constant curvature (as I understand 
it the approach presented here  is the 
Lagrangian dual of his approach).

The idea is to suppose you have a topological 
triangulation of a three manifold, let $\frak{N}$ 
be the possible dihedral angle data such that the sum of the angles about 
an edge is $2 \pi$ and such that the data can be used to produce a 
hyperbolic simplex. Call the dihedral angles value at $e$ in $t$  
$\alpha^{t}_e$.  Now  we can deform our angle taking the 
span of the transformations $w^{t_1,t_2}_e$ 
which adds  one unit of angle to 
$\alpha^{t_1}_e$ as subtracts  one unit form $\alpha^{t_2}_e$, for  
neighboring simplexes $t_1$ and $t_2$.  Let $\frak{N}$ be the 
angles geometrically equivalent to $x$.

On $\frak{N}$ we can now place the energy 
\[ E(x) = \sum_{t \in \frak{P}} V_t(x) \]
 where the sum is over all the simplexes $t$ and $V_t(x)$ is 
the hyperbolic volume 
of $t$ relative to $x$'s angles data.

Now, as in the the two-dimensional discrete case, 
the critical points of this energy are a 
sets of simplexes which fit together. To see this we recall  
Schafli's formula which tells us 

\[ d V_t = \sum_{e \in t} l_e^t d \alpha^t_e , \]
where $l^t_e$ is edge length of $e$.

So at a critical point we have
\[ dE(w^{t_1,t_2}_E) = l^{t_1}_e - l^{t_2}_e .\]

Unlike in the two dimensional case we fail to have good boundary control (or 
Casson would have already uniformized the manifolds of interest). 
In particular simplexes may collapse and enough collapsing may take place that 
even the topological type of the complex could change.  
This is where the randomizing may help.

\end{subsubsection}

\begin{subsubsection}{From the Discrete to the Continuous}

Now just as before we  randomly Delaunay triangulate.
We use a Poisson point process to distribute the points and 
then assign a simplex to a quadruple  on a small enough 
sphere passing through the four  points which is  
empty of other points.  With high
probability this complex will form a  triangulation. 
However not as 
canonically as in the two dimensional case, since there 
are no natural geometric simplexes and 
faces only natural edges and vertices.   

For any set of vertices forming a triangulation 
we can construct a 
point in some $\frak{N}$, just as in the discrete case.  
It is essential to note that while we have no faces  
and hence no dihedral angles, we have dihedral 
angles at both end points (using the geodesic directions) 
so can average them and form the angle data of some $y$ in some
$\frak{N}$.  So we may use the energy to once again form a random variable 
which can be computed.

\end{subsubsection}

\begin{subsubsection}{Continuous Uniformization}

Just as in the two dimensional case we may  now use  
\[ E(h) = \lim_{\lambda \rightarrow \infty} \Bbb E_{\lambda}(E_h - E_g),\] 
to form an energy on the space of metrics.  
Currently not much is known about this energy, but the hope of course 
is that it 
will 
form an energy with which the techniques of section 
\ref{bigproof} can be carried out.

In the end hopefully one will be able to see that a three 
manifold admitting a metric of variable negative sectional 
curvature will accept one 
of constant negative sectional curvature.

\end{subsubsection}
\end{section}
\end{chapter}

\begin{chapter}{The Discrete Uniformization Theorem}\label{dunith}

This chapter is dedicated to proving the 
generalization of the Thurston-Andreev theorem
mentioned in the introduction.  It is presented in a potentially 
strange 
order from the point of view of presentation in section \ref{conu}.  
In section \ref{tdt}  the details of the energy argument 
presented in \ref{dsc}
are given ending in a  theorem giving  
conditions  under which triangle angles 
in a topological triangular decomposition can be conformally 
deformed to the angles of a  geodesic triangular decomposition of a hyperbolic 
surface.  Also in this  section a  connection between the disk patterns 
mentioned in section \ref{dsc} and
certain  hyperbolic polyhedra is explored and exploited. 

In section \ref{yipyeep} we explore 
the linear part of the discussion in section \ref{dsc}, and prove
a warm up case of the general ideal convex 
Thurston-Andreev theorem.  This warm up case 
is theorem   
 \ref{res1} form the introduction including the possibility of $M$ 
having boundary and has as an  immediate consequence 
the teleportation lemma (lemma \ref{dtele}).
In section \ref{itat} the general convex ideal case of the Thurston-Andreev 
theorem  is presented and 
dealt with.

\begin{section}{The Triangular Decomposition Theorem}\label{tdt}
\begin{subsection}{Statement and Notation}\label{san}
Throughout this paper $M$ will  denote a compact two-dimensional 
surface 
with $\chi(M) < 0$.  
By geometry I will mean a hyperbolic structure.  
Uniqueness of geometries, triangulations and disk patterns is 
of course up to isometry.  

The main theorem in this section  really 
should be stated for the following 
structure which generalizes the notion of triangulation.

\begin{definition}\label{decomp}
Let a  triangular decomposition, $\frak{T}$, 
be a cell decomposition of $M$  that  lifts to a triangulation
in $M$'s universal cover. 
\end{definition} 

We will 
keep track of the combinatorics of such a decomposition by denoting 
the vertices as $\{ v_1 , ... v_V \}$, 
the edges  as $\{ e_1  \dots e_E \}$  and the triangles 
as $\{ t_1 \dots t_F \} $.  For convenience I will let $E$, $V$, 
and $F$ denote 
both the set of edges, vertices, and faces and the cardinalities 
of these sets.  
Similarly for the subsets of $E$ and $V$ contained in $M's$ possible 
boundary, denoted $\partial E$ and $\partial V$. 
Let $\{e \in S\}$  denote   the set of  
edges on the surface in a collection  of 
triangles $S$,
and let $\{e \in v\}$  denote all the edges associated to a vertex $v$ 
as if counted in the universal cover.  The set of triangles containing a
vertex  $\{ t \in v\}$ has the special name of the flower at $v$.

\begin{figure*}
\vspace{.01in}
\hspace*{\fill}
\epsfysize = 1.5in  
\epsfbox{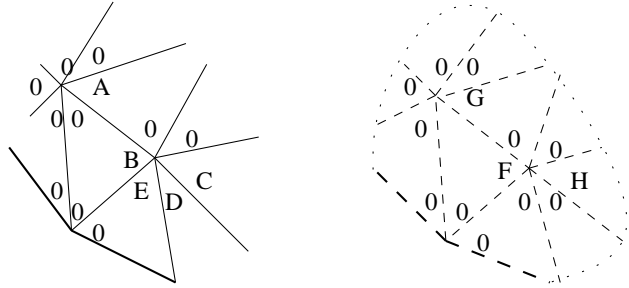}
\hspace*{\fill}
\vspace{.01in}
\caption{\label{the5} A Covector and a Vector}
\end{figure*}

Note that in a triangular decomposition there are 
$3F$ slots $\{ \alpha_i \}$ in which one can insert 
possible triangle angles, which we 
will place an order on and identify with a basis of a 
$3 F$ dimensional real vector space.  With this basis choice we will 
denote this vector space as   
$\Bbb R^{3F}$, and  denote vectors in it as $x =\sum A^i \alpha_i$.
Further more let 
$\alpha^i$ be a dual vector such that $\alpha^i(\alpha_j) = \delta_i^j$. 
With this we will view the 
angle at the slot $\alpha_i$ as  $\alpha^i(x) = A^i$.
It is rarely necessary to use this notation and instead 
to use the actual geometry. We will denote a vector 
by placing the $A^i$ coefficients in a copy of the triangulation 
with dashed lines and covector will contain its  
coefficients $A_i$ 
in a copy of the triangular decomposition with solid lines.  
Thick lines will
always denote a boundary edge, as  in lower left corner of 
vector and covector in figure  \ref{the5}.
If in the picture we mean the non-specified values 
to be arbitrary we will surround the picture with a loop (see the vector  
in figure \ref{the5})  and 
if we mean   the non-specified values  to be zero 
the picture will not be surrounded (see the covector in figure \ref{the5}).  
The pairing of a vector and a covector denoted 
$ \sum_{\{\alpha_i\}} A_i \alpha^i( \sum_{\{\alpha_i\}} A^j \alpha_j)$ 
can be viewed geometrically by placing the copy of the  
triangular decomposition corresponding to the vector on top of the 
triangular decomposition  corresponding to the covector and multiplying 
the numbers living in the same angle slots to arrive at 
$\sum_{\{\alpha_i\}} A^i A_i$ 
(see figure \ref{the6}).
For a  triangle $t$ containing the angle slots 
$\alpha_i$, $\alpha_j$, and $\alpha_k$ 
let $d^t (x)   = \{ A^i,A^j,A^k \}$
and call $d^t(x)$  the angle data associated to $t$. 

\begin{figure*}
\vspace{.01in}
\hspace*{\fill}
\epsfysize = 1.5in  
\epsfbox{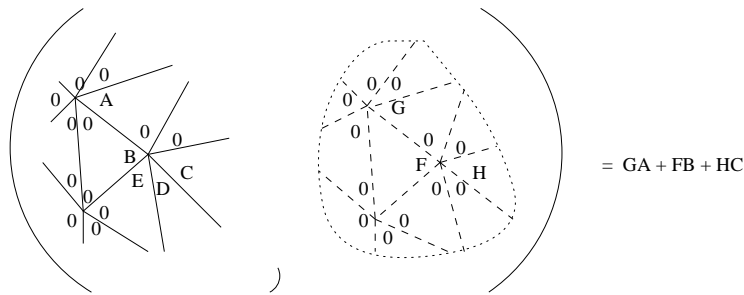}
\hspace*{\fill}
\vspace{.01in}
\caption{\label{the6} A Pairing}
\end{figure*}

\begin{figure*}
\vspace{.01in}
\hspace*{\fill}
\epsfysize = 1.5in  
\epsfbox{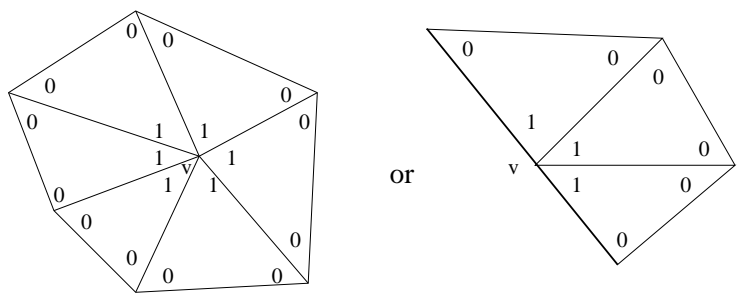}
\hspace*{\fill}
\vspace{.01in}
\caption{\label{co0} The $p^v$ Covector}
\end{figure*}

In order to live on an actual nonsingular geometric surface
all such angles should be required to live  in the subset of 
$\Bbb R^{3F}$ where the angles at an interior 
 vertex sum to $2 \pi$ and the angles at a boundary vertex sum to 
$\pi$. 
Let 
$p^v$ be the covector living 
in the flower at $v$  defined as in figure \ref{co0},
this encourages us to choose  our possible angles in the affine flat
\[ V = \{x \in \Bbb R^{3F} \mid p^v(x) = 2 \pi \mbox{ for all } v 
\in V -\partial V \mbox{ and } p^v(x) =  \pi 
\mbox{ for all } v \in \partial V  \}. \]

To further limit down the possible angle values we define 
the covector $l^t$ as in figure \ref{co2}
 \begin{figure*}
\vspace{.01in}
\hspace*{\fill}
\epsfysize = 1in  
\epsfbox{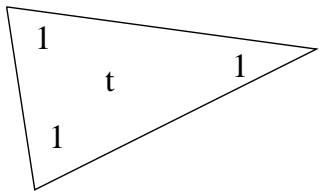}
\hspace*{\fill}
\vspace{.01in}
\caption{\label{co2} The $l^t$ Covector}
\end{figure*}
and note by the Gauss-Bonnet theorem that 
\[ k^t(x) = l^t(x) - \pi \] 
would be the curvature in  a geodesic triangle with 
angle data $d^t(x)$. 
We will now isolate the open convex subset of $V$ 
where the curvature is negative and all angles are realistic.
\begin{definition}\label{dcc}
Let an angle system be a point in  
  \[ \frak{N} = \{ x \in V \mid k^t(x) < 0 \mbox{ for all } t 
 \mbox{ and  } \alpha^i(x) \in (0, \pi) \mbox{ for all } \alpha^i\}. \]
\end{definition}
Note the actual angle data of a geodesic 
triangulation of a  
surface with negative curvature  
has its angle data living in this set.

Observe that from  any set of triangle data $d^t(x) = \{ A,B,C\}$ 
with the angles in $(0,\pi)$ and $A+B+C - \pi  < 0$ 
we may associate an
actual hyperbolic triangle, 
call this triangle $t(x)$.  
Suppose the triangles in $\{t(x) \mid t \in \frak{T}\}$ 
fit together in the sense that 
all the corresponding edges are the same lengths.  Then
$\frak{T}$ being a triangular decomposition implies every open flower 
is embedded in $M$'s universal cover and  when the edge lengths all agree
this flower can  be given  
a hyperbolic structure which  is consistent on flower over laps.  
So we have formed  a hyperbolic structure
on $M$. 

\begin{definition}\label{duas}
Call  an angle system $u$ uniform 
if all the hyperbolic realizations of the triangles in $u$
fit together to form a hyperbolic structure on $M$. Let 
$\frak{T}(u)$ denote the the geodesic triangular decomposition 
corresponding to $\frak{T}$ and $u$. Further more let
\[ \frak{U} = \{ u \in \frak{N} \mid u \mbox{ is uniform } \}.\]
\end{definition}

In section 2.4  we will  attempt to take a point in $\frak{N}$ and 
deform it into a point of $\frak{U}$.  Such deformations are located in  an
affine space and I will call them conformal deformations (see 
the introduction to section \ref{conu} 
%\cite{Le} 
to  motivate this terminology).  
To describe this  affine space for each edge $e \in E -\partial E$ 
construct a vector $w_e$ as in figure \ref{ve1}.
\begin{figure*}
\vspace{.01in}
\hspace*{\fill}
\epsfysize = 1in  
\epsfbox{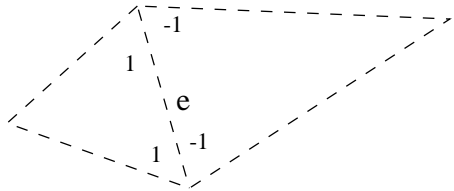}
\hspace*{\fill}
\vspace{.01in}
\caption{\label{ve1} The $w_e$ Vector}
\end{figure*} 
\begin{definition}\label{condef}  
A conformal deformation will be a vector in 
\[ C = span\{ w_e \mid \mbox{ for all } e \in E - \partial E \}, \]
and call $x$ and $y$ conformally equivalent if $x - y \in C$.
\end{definition}
The first thing worth noting is that if $x \in V$ and 
$y$ is conformally equivalent to $x$
then as
an immediate consequence of geometrically  
pairing the covector in figure \ref{co0} with the vector  in figure 
\ref{ve1} we have for $v \in V -\partial V$ that 
\[ p^v(y) = p^v\left(x + \sum_{e \in \frak{P}} B^{e} w_{e}\right) 
= p^v(x) + \sum_{e \in \frak{P}} B^{e} p^v(w_{e} ) 
= 2 \pi  + 0, \]
so $y$ is also in $V$.    Similarly for $v \in \partial V$.

To combinatorially understand the points in $\frak{N}$ 
which we may  conformally deform into uniform 
structures it is useful to express a particularly 
nasty set in the boundary of $\frak{N}$.  
\begin{definition}\label{legal}
Let $t$  be called a legal with respect to $x \in \partial \frak{N}$ if 
$d^t(x)  = \{A^1,A^2,A^3\}  \neq \{ 0,0,\pi \}$ yet either $k^t(x) =0$
or for some $i$ we 
have $A^i =0$. Let  
\[ B = \{x \in \partial \frak{N} \mid  \mbox{ x contains no legal  } t \}. \] 
\end{definition}

We will prove the following theorem.
\begin{theorem}\label{tri1}
If there a uniform angle system conformally 
equivalent to $x$ then it is unique, and  for any angle 
system $x$ with   $(x + C) \bigcap B$ empty 
there exists a conformally equivalent uniform angle system. 
\end{theorem}

%\begin{note}\label{bound1}
%It is easy to modify this set up so that this 
%theorem  applies to surfaces with 
%boundary if $\chi(M) <0$.  First to enforce a straight boundary 
%our negative curvature angle systems  must be modified 
%so that at a boundary vertex $p^v(x) =\pi$.  The only real change will
%be that a  conformal transformation 
%will be an element of the the span of the $w_e$ over only the interior edges.
%\end{note}

Much of what takes place here relies on certain basic 
invariants of conformal 
deformations. 
To describe them  for each triangle $t$ and $e \in t$
we form the covector 
$\psi^e_t$ as in figure \ref{co3}.  For each edge $e \in \partial V$ 
we will denote  $\psi^e_t$ as $\psi^e$ while 
for each edge $e \in V - \partial V$ associated with triangles
$t_1$ and $t_2$ we will let
\[ \psi^{e} =  \psi^{e_1}_{t_0} + \psi^{e_1}_{t_1}.\]    
We will call the  $\psi^e$ covector  the formal angle defect at $e$.
Let the formal intersection angle be determined by
\[ \theta^e(x) = \pi - \psi^e(x)\]
when $e \in E -\partial E$.
 and 
\[ \theta^e(x) = \frac{\pi}{2} - \psi^e(x)\]
when $e \in \partial E$.

\begin{figure*}
\vspace{.01in}
\hspace*{\fill}
\epsfysize = 1in  
\epsfbox{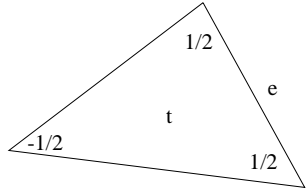}
\hspace*{\fill}
\vspace{.01in}
\caption{\label{co3} The $\psi^e_t$ Covector}
\end{figure*} 
Looking at the pairing 
between a covector  $\psi^e$ and  
the vectors 
spanning $C$ in  figure \ref{ve1},
we see if $y$ is conformally equivalent to $x$ then 
\[ \psi^e(y) = \psi^e\left(x + \sum_{f \in \frak{P}} B^f w_{f}\right) 
= \psi^e(x) + \sum_{f \in \frak{P}} B^f \psi^e(w_{f} ) 
= \psi^e(x),\]
and indeed for each relevant edge $e$ we see 
$\psi^e$ and $\theta^e$ are conformal invariants.

It is worth noting the trivial but extremely useful fact that the 
curvature assumption gives us some control of the angle discrepancy.

\begin{fact}\label{ob1}
When $x \in \frak{N}$ we have
$\psi_{t_i}^e(x) \in \left(\frac{-\pi}{2},\frac{\pi}{2} \right)$, and 
when $x \in \partial \frak{N}$ we have
$\psi_{t_i}^e(x) \in \left[\frac{-\pi}{2},\frac{\pi}{2} \right] $.
\end{fact}
{ \bf Proof:}

Let  $d^{t_i}(x) = \{A,B,C\} $ and note since 
$B +C \leq A+B+C = l^t(x) < \pi$ and $A<\pi$ we have  
\[ -\frac{\pi}{2} < -\frac{A}{2} \leq  \psi^e_{t_i} = \frac{B+C-A}{2}
\leq \frac{B +C}{2} < \frac{\pi}{2} . \]
The boundary statement follows from the possibility 
of these inequalities becoming equalities.
\qed
%\begin{note}\label{bound2}
%In the surface with  boundary case at a boundary edge $e$ one uses the 
%covector 
%$\psi^e(x)$ which is an interior $\theta^e$ restricted to the only 
%triangle present at the boundary edge.
% is the covector in figure \ref{bound}.     
%\end{note}

%\begin{figure*}
%\vspace{.01in}
%\hspace*{\fill}
%\epsfysize = 1.5in  
%\epsfbox{bound.eps}
%\hspace*{\fill}
%\vspace{.01in}
%\caption{\label{bound} A Boundary Covector}
%\end{figure*}

As we shall see, the  fundamental 
reason why the above triangulation theorem is  
related to circle patterns and polyhedra construction is that 
at a uniform angle system $\theta^e(u)$ can be 
geometrically interpreted as the intersection angle
of the circumscribing circles of two 
hyperbolic triangles meeting along $e \in E -\partial E $.

Equivalently $\theta^e$ will be realized as  
the dihedral angle at an edge  of a polyhedra  in the class $I$ 
corresponding to 
$e$.  
It's high time to describe this class $I$.

\end{subsection}

%\begin{section}{Proof of the Triangular Decomposition Theorem

\begin{subsection}{Polyhedra in the Class $I$}\label{poly}
In this section we will construct and examine a bit of the geometry
of the  infinite sided  ideal polyhedra 
that will arise in the proof of theorem \ref{tri1}.
These polyhedra can be constructed out of building blocks each in the form 
of an ideal prism.
%To describe these polyhedra it is useful to let $H^2$ denote  a specified copy of $H^2$ in $H^3$ which will always be drawn as the $z =0$ plane in the ball model of $H^3$, and extend any isometries of $H^2$ to $H^3$ via this embedding.   

\begin{definition} \label{pridef}
Place a given $t(x)$ on a copy of $H^2 \subset H^3$. 
Let $P_t(x)$ be the convex hull of 
the  set  consisting of $t(x)$ 
unioned with   the geodesics 
perpendicular to this $H^2 \subset H^3$ 
going through $t(x)$'s 
vertices. See figure \ref{pri}.
\end{definition}

Given a triangular decomposition and uniform angle system $u$ let 
$\tilde{\frak{T}}(u)$ 
denote the lift of $\frak{T}(u)$ to $H^2$.  With this 
notion we may now construct the polyhedra in $I$. 

\begin{definition}\label{polyhedra}
Let the class of polygons $I$ be those 
constructed by placing
$\tilde{\frak{T}}(u)$  on an $H^2 \subset H^3$ and forming  
$\bigcup_{t \in \tilde{\frak{T}}} P_t(u)$.  
\end{definition}

\begin{figure*}
\vspace{.01in}
\hspace*{\fill}
\epsfysize = 1.5in 
\epsfbox{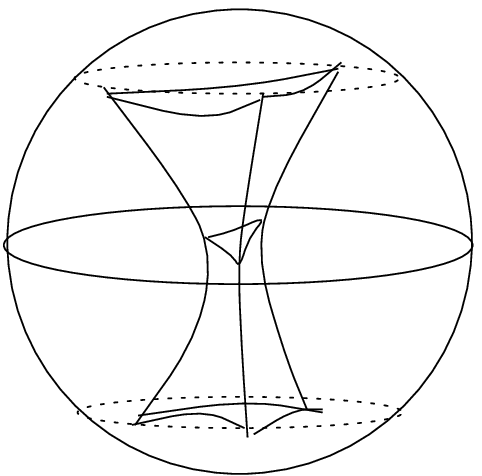}
\hspace*{\fill}
\vspace{.01in}
\caption{\label{pri} The Ideal Prism $P_t(x)$}
\end{figure*}

Notice abstractly a polyhedron $P$ in $I$  
is an ideal polyhedra symmetric 
upon  the reflection through some plane  
upon which lives a group of isometries
forming a compact surface  which 
when extended to $H^3$  are isometries of $P$.  

Our first observation will be that if 
we happen to know a  $\frak{T}$ and $u$ forming $P$,
then from this data we can  easily construct the dihedral angles.  
We will also refer to {\bf the} dihedral angle associate to an 
edge of $\frak{T}$, despite the fact there is 
always in fact a symmetric pair of such angles.

\begin{formula} \label{ptheta}
If $P$ is constructed from $\frak{T}$ and $u$ as in definition 
\ref{polyhedra}
then the dihedral 
angle at an edge of $P$  associate to the  edge $e$ of $\frak{T}$ 
is given by  $\theta^e(u)$. 
\end{formula}
{\bf Proof:}
Note that the needed dihedral  angle is the sum of the angles in  
$P_{t_1}(u)$ and $P_{t_2}(u)$ corresponding to $e$. 
It is these angles that will be computed.

To do this simply note that $P_{t_i}(u)$ can be decomposed as  
into three ideal tetrahedra as in figure \ref{tetdecomp}, where  
$d^{t_i}(u)=\{A,B,C\}$ with the angle slot containing 
the $A$ coefficient across form $e$.  
The labeled angles in the figure \ref{tetdecomp}
are the internal angles closest to the label.

\begin{figure*}
\vspace{.01in}
\hspace*{\fill}
\epsfysize = 3in 
\epsfbox{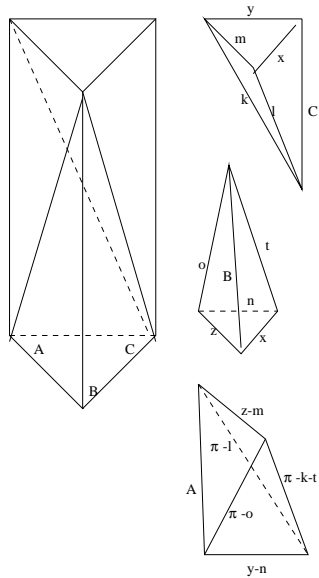}
\hspace*{\fill}
\vspace{.01in}
\caption{\label{tetdecomp} A Decomposition of $P_t(x)$}
\end{figure*}

Now use the fact that at the vertex of an ideal tetrahedron 
the angles sum to $\pi$, to 
form ten linear equations in the labeled unknowns.
Solving in terms of $A = A^i$,$B = B^i$, and $C=C^i$ one finds
\[ x = l =\frac{\pi + A -B -C}{2 } \]
as needed. 

\qed

\begin{note} \label{bound3}
The polygons in the surface with boundary case now are extremely 
non-convex and have faces between the two hemispheres.  Note 
that the $\{ \theta^e(x) \}$
still represent the dihedral angles in these faces.
\end{note}

Let $G$ be a compact surface group 
extend form the $H^2$ through which $P$ is symmetric 
and under which $P$ is invariant.  To be in the class $I$ such a group exists.
When such a $G$ is chosen we  will be interested in the volume of 
$P/G$.  Notice given such a $G$ we may choose a geodesic  
triangulation 
which descends to a triangular decomposition $\frak{T}$ of $H^2/G$ with an 
associated uniform structure $u$ from which 
$P$ is constructed as in definition 
\ref{polyhedra}.  Using this triangulation we have. 
\begin{formula}\label{volume}
The volume of $P/G$ is $\sum_{t \in \frak{T}} V_t(u)$ where  
\[V_t(u) = \Lambda(A) +\Lambda(B)+\Lambda(C)+ 
\Lambda\left(\frac{\pi-A-B-C}{2}\right) \]
\[
+\Lambda \left(\frac{\pi+A-B-C}{2}\right)
+\Lambda\left(\frac{\pi+B-A-C}{2}\right)+
\Lambda\left(\frac{\pi+C-A-B}{2}\right) ,\]
with  $\Lambda $ is the Lobacevskii function
 \[ \Lambda(\alpha) = -\int_{0}^{\alpha} \ln(2 |\sin(t)|) dt. \]
\end{formula}
{\bf Proof:} 
First note the volume of 
$P/G$  can be reduced to the volumes of individual $P_t$ 
by noting  the volume is 
$\sum_{t \in \frak{T}} V_t(u)$, 
where $V_t(u)$ is the volume of $P_t(u)$.

Its useful to get a formula for $V_t(u)$
Recall if a tetrahedra has angles $\alpha$, $\beta$, and $\gamma$
 meeting at 
 an ideal vertex then
 its volume is \[ \Lambda(\alpha) + \Lambda(\beta) + \Lambda(\gamma)\]
where $\Lambda $ is the Lobacevskii function.

In the proof of formula  \ref{ptheta}  
we decomposed $P_t(u)$ with $d^t(u) = \{A,B,C\}$  
into three ideal tetrahedra,  and found linear equations 
determining all the angles in these ideal tetrahedra 
in terms of $A$, $B$, and $C$.   We already wrote down 
the angle corresponding to $x$ and $l$ and can further note 
\[ y=\frac{\pi+ A -B -C}{2 } \]
  \[ z=\frac{\pi + C -A -B}{2 } \]
  \[ y-n=\frac{\pi - A -B -C}{2 } \]
  \[m = C.\]
Plugging  these angles in the tetrahedra's volume formula 
gives us the needed formula.
%\[V_t(u) = \Lambda(A) +\Lambda(B)+\Lambda(C)+ 
%\Lambda\left(\frac{\pi-A-B-C}{2}\right) \]
%\[
%+\Lambda \left(\frac{\pi+A-B-C}{2}\right)
%+\Lambda\left(\frac{\pi+B-A-C}{2}\right)+
%\Lambda\left(\frac{\pi+C-A-B}{2}\right) ,\]

%as needed.
\qed

Notice that when $P$ is constructed form $\frak{T}$ and an edge  $e$ 
of $\frak{T}$ is associate a dihedral angle $\pi$ that  
$P$'s edge corresponding to $e$ is fake in the sense that the two
triangular faces meeting at this edge share the same polyhedron face.
If a polyhedra's face is triangulated then a change in this 
triangulation will not 
effect the polyhedra so   
it will turn out convenient to have articulated the 
topological cell divisions which naturally arise form polyhedra in $I$. 

\begin{definition}\label{polytion}
Let a polygonation of a surface be a locally finite cell division 
such that the closure of a cell is 
an embedded polygon and the intersection of two polygons 
is empty, 
contains a single point, or contains 
one edge and the two vertices associated to the edge.
Let  a  polygonal decomposition $\frak{P}$ be 
a cell decomposition of a compact surface 
which lifts to a polygonation in its universal cover. 
\end{definition}

To keep track of all the the combinatorics of such a cell division we
will use the same notation as we did for triangular decompositions.

At this point it is useful to name what turns out to be the appropriate 
 home of the possible dihedral angle assignments.
Let
$\{ e_i \}$ be the set of edges 
in a polygonal decomposition and just as we did with the angle 
slots 
let them correspond to the the basis vectors of 
an $E$  dimensional vector space, which we will denote $\Bbb R^E$ with this 
basis choice.  We will be viewing this as the space of 
possible angle discrepancies. Denote these vectors as 
$p = \sum \psi^{e_i} e_i$.

Notice the $\psi^e$ were covectors in the previous section. 
We can motivate this abuse 
of notation in the case that $\frak{P}$ is triangular decomposition 
by letting
\[ \Psi: \Bbb R^{3 F} \rightarrow \Bbb R^{E} \]
be the linear mapping 
given by 
\begin{eqnarray}\label{mapping}
 \Psi(x) = \sum \psi^{e_i}(x) e_i,
\end{eqnarray}
and noting by the above lemma that we do indeed hit the dihedral angle
discrepancies when using a uniform angle system.
As a further justifiable abuse notation we let 
$\theta^{e_i}(p) = \pi - e^i(p) = \pi - \psi^{e_i}$ when 
$e^i \in E -\partial E$ and 
$\theta^{e_i}(p) = \frac{\pi}{2} - e^i(p) =\frac{\pi}{2} - \psi^{e_i}$
 when 
$e^i \in \partial E$.

At this point it is convenient associate explicitly  
the data contained in polyhedra in the  class $I$ with certain disk patterns.

\end{subsection}
\begin{subsection}{Ideal Disk Patterns}\label{idp}

In this section we will discuss the relation of these polyhedra to 
disk patterns specified by combinatorial and topological data.   
The topological data comes in the form of 
a polygonal decomposition $\frak{P}$
and the extra data 
 associated to $\frak{P}$ 
will be a point $p \in \Bbb{R}^E$.

To begin to articulate the pattern here It's necessary to 
define a particularly nice polygonal decomposition.

\begin{figure*}
\vspace{.01in}
\hspace*{\fill}
\epsfysize = 2in 
\epsfbox{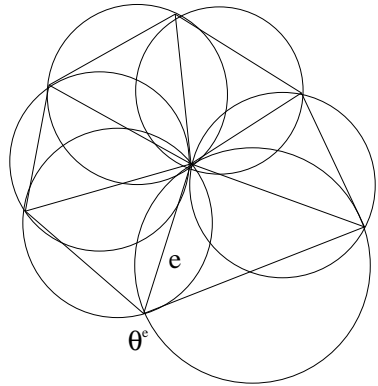}
\hspace*{\fill}
\vspace{.01in}
\caption{\label{the0} An Ideal Disk Pattern}
\end{figure*}

\begin{definition} \label{circum}
A polygonal decomposition on a geometric surface 
is called circumscribable 
if it is a geodesic polygonal decomposition where 
each polygon is circumscribed by circle in $M$'s universal cover.
\end{definition}

Given a circumscribable $\frak{P}$  we may use 
a point
$p \in \Bbb R^E$  to keep track of all angles of intersection between 
the circumscribing circles of the polygons meeting at $e$ (namely at $e$ this
angle is $\theta^e(p)$). 
Simply to articulate the sense of angle to be used here 
it is convenient to introduce $\frak{P}^t$ which is a triangular decomposition 
associated to $\frak{P}$ by triangulating  each polygon. 
It is occasionally useful to have an explicit  grip on 
this triangulation so we may assume that  we 
triangulate each polygon with a fan as in figure \ref{fan}.
Note this triangulation can be chosen to be a geodesic triangulation 
if $\frak{P}$ is circumscribable,  since all the polygons are then convex.
Also if we are given a $p \in \Bbb R^E$ relative to $\frak{P}$ 
we will let  $\hat{p}$ be the 
angle discrepancy assignment on $\frak{P}^t$ 
which is $\psi^e(p)$ on the the edges of $\frak{P}$ 
and $0$ on the new edges, and hence corresponds 
to geodesically triangulating a circumscribable decomposition.

\begin{figure*}
\vspace{.01in}
\hspace*{\fill}
\epsfysize = 1in 
\epsfbox{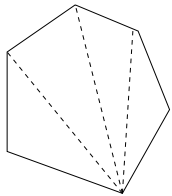}
\hspace*{\fill}
\vspace{.01in}
\caption{\label{fan} A Fan Associated to a Polygon}
\end{figure*}

Note given a uniform structure $u$ associated to $\frak{P}^t$ 
that every edge sits between two 
triangles $t_1$ and $t_2$, 
and the sense of the angle used here can be chosen   
relative to  these geodesic triangles.  
Let the intersection angle $\theta^e(\hat{p})$  be  a 
number in $(0,\pi)$ if the vertex of $t_1$ not on $e$ 
is out side $t_2$'s circumscribing disk, in 
$(\pi,2\pi)$ if this vertex is in $t_2$'s circumscribing disk, 
and exactly $\pi$ when this vertex is on $t_2$'s circumscribing circle.

\begin{definition} \label{idealpat}
Let an ideal disk pattern be 
a collection distinct 
disks on  a geometric surface 
whose boundary circles are 
in one to one correspondence with  the circumscribing circles 
of a circumscribable polygonal 
decomposition.
Notice in such a situation using a $\frak{P}^t$ 
we naturally have the triangles needed to 
associate  $\hat{p} \in (0,2 \pi)^E$ to keep track 
intersection angles at
each of the polygonal decomposition's edges.
\end{definition}

The pattern problem is given a topological polygonal  
decomposition  
and a 
suitable  $p \in \Bbb R^E$ to assert the existence 
and uniqueness of an ideal  disk  pattern.    
Such assertion are equivalent to such 
assertions about polyhedra in the class $I$.
 
\begin{observation}\label{humding}
Every ideal disk pattern can be 
associated a unique polyhedra  in $I$.  The polyhedra associated 
to a pattern is convex 
if and only if $p \in (0,\pi)^E$, in which case we will
also call the pattern convex.
To each polyhedra 
in $I$ and choice of $G$ as described  preceding formula \ref{volume} 
there is a uniquely associated ideal disk 
pattern.  Furthermore  under these correspondences 
the dihedral angle of the 
polyhedra  associated to an edge $e$ is precisely the 
intersection angle between the 
circumscribing circles 
of the polygons sharing 
$e$ in the ideal disk pattern.  
\end{observation}
{\bf Proof:}
Note that given a  circumscribable polygonal decomposition 
that we may form a geodesic triangular decomposition  
$\frak{P}^t(u)$ which can be associated an element of $I$ via the construction 
in definition \ref{polyhedra}.  Similarly given an element $P$ in $I$
and a $G$ as in the discussion preceding formula \ref{volume} 
we know there is a geodesic triangulation 
invariant under  $G$ which descends to a geodesic 
triangular decomposition 
of  $H^2/G$.  Let $\frak{P}$ be the circumscribable decomposition  
formed by ignoring the edges where $\psi^e(u) = 0$.

To see the angle correspondence   we will recall 
a map from a  specified $H^2 \subset H^3 $  
to the upper half of 
the sphere at infinity $S^{\infty}_u$  (using
the usual  conformal structure of the sphere at infinity) which sends 
circles in $H^2$ to circles in $S^{\infty}_u$.  
To form this map  send a point in this fixed $H^2$ to where the geodesic 
perpendicular to this $H^2$  hits $S^{\infty}_u$ (as in figure \ref{hit}). 
Note any circle can be sent to what 
we view as the center of $H^2$ via a hyperbolic isometry preserving $H^2$, 
were by symmetry it is sent under this mapping to a circle at infinity. 
This  isometry bringing
the circle to what we view as the center 
induces a Mobius transformation on the 
sphere at infinity and preserves the set of geodesics used to form this
map's image - so the fact that Mobius transformations send circle 
to circles on $S^{\infty}$ now gives us that indeed the image of any   
circle is a circle.
\begin{figure*}
\vspace{.01in}
\hspace*{\fill}
\epsfysize = 1.5in 
\epsfbox{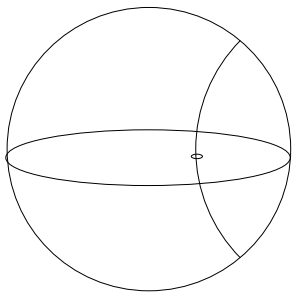}
\hspace*{\fill}
\vspace{.01in}
\caption{\label{hit} A Conformal Mapping}
\end{figure*}

%Using this map (and its lower hemisphere counter part)
%the  vertices of the triangle are sent to six points whose convex hull
%forms $P_t$.  
%Note that the the spheres perpendicular to $H^2$ whose intersections with 
So a neighboring pair triangles $t_1$ and $t_2$ have circumscribing 
circles in 
$H^2$ sent  under this map  sent to   
circles at infinity intersecting at the same angle 
and going through the ideal
points of the neighboring $P_{t_1}(u)$ and $P_{t_2}(u)$.  
But these  circles at infinity  are also the intersection 
of $S^{\infty}$ with the spheres representing 
the hyperbolic planes forming the top faces of 
$P_{t_1}(u)$ and $P_{t_2}(u)$.   
So the intersection angle of these spheres is precisely the dihedral angle, 
which is now seen to be the intersection angle of the 
circles on the sphere at infinity, or finally 
the original intersection angle of the circumscribing 
circles.  

The convexity assertion can be seen immediately by 
looking at the two neighboring $P_t(u)$ prism's.

\qed

At this this point there are two natural questions about 
such patterns and polyhedra, namely when they exist are they unique
and are there nice way to insure existence? 
The first question will be answered in the next section 
with the following fact:
a pattern/polyhedra  is determined uniquely by 
its/an associated topological polygonal decomposition 
$\frak{P}$ and $p \in (0,2\pi)^E$.
The existence issue will be handled in sections 
\ref{yipyeep} and \ref{pitat} 
where necessary and sufficient conditions for 
$p \in (0,2\pi)^E$  relative to $\frak{P}$ to 
be associated to a convex pattern/polyhedron 
with this given data will be presented.

With respect to these existence and uniqueness results 
I will only deal with the terminology of ideal 
disk patterns from here on out.

\end{subsection}
\begin{subsection}{Proof of Theorem \ref{tri1}}\label{prot}
Now we will prove theorem \ref{tri1}.
The  proof relies on   
the energy introduced in section \ref{dsc} 
which lives on  $\frak{N}$, and which at a uniform $u$ agrees 
with  the  volume of $P/G$ with  $P$ constructed from
$\frak{P}(u)$ as in definition \ref{polyhedra} and $G$ the surface group 
associated to $\frak{P}(u)$. We may assume we that 
either $\frak{P}$ is a triangulation or form it we have constructed the
 triangulation $\frak{P}^t$, so we will denote in $\frak{T}$.

Let the energy be the of the volume of the 
abstract disjoint union $\bigcup_{t \in \frak{T}} P_t(x)$
given by 
\[E(x) = \sum_{t \in \frak{T}} V_t(x). \]
Since hyperbolic objects optimize at fat objects we will be maximizing 
this energy, and  I suppose to call this an energy (in the physical sense)
I really should negate it.  However I like both the fact that hyperbolic 
objects like to be fat and the term energy, 
so I will simply warn the reader about this odd terminology.  
Note we get 
an explicit description of the energy from formula \ref{volume}.
Using this  formula
   we can differentiate to find $E$'s differential, $dE^x$.  As usual for 
a function in a linear space like $\Bbb R^{3F}$
we  use translation to 
identify the tangent and cotangent spaces at
every point with  $\Bbb R^{3F}$ 
and  $(\Bbb R^{3F})^*$ and express our differentials in the chosen basis.   
From the formula for the Lobacheski function we have    
$dE^x = \sum_{\alpha_i \in \frak{T}} E_i(x) \alpha^i$
with  
   \[E_i(x) = - \frac{1}{2} ln\left(\frac{
 \sin^2(A) \sin\left( \frac{\pi+A-B-C}{2}\right)}
   { \sin\left( \frac{\pi+B-A-C}{2}\right)
   \sin\left( \frac{\pi+C-A-B}{2}\right)
   \sin\left( \frac{\pi-A-B-C}{2}\right) } \right),\]
where $\alpha^i(x)=A$, $\alpha_i \in t$, and $d^t(x) = \{A,B,C\}$.  
After  a little bit of trigonometry this  can be simplified to
\[ E_i(x) = - \frac{1}{2} ln \left( \frac{\sin(A)^2 (\cos(B+C) + \cos(A))}
 {(\cos(A+C) + \cos(B))(\cos(A+B) + \cos(C))}\right). \]
   
  To compute further  let  $a$, $b$ and $c$ denote the
 edge lengths opposite to the angles  $A$, $B$, and $C$ respectively 
in the hyperbolic 
triangle determined by $\{A,B,C\}$. 
From elementary hyperbolic geometry we know
    that
    \[ \cosh(a) = \frac{\cos(C) \cos(B) + \cos(A)}{\sin(B) \sin(C)},\] 
    allowing us to simplify  $\cosh(a) -1 $
    to \[ \frac{\cos(B+C) + \cos(A)}{\sin(B)\sin(C)}.\]  
    So plugging this in the above  formula
we arrive at...
\begin{formula} \label{diff}
$dE^x = \sum_{\alpha_i \in \frak{T}} E_i(x) \alpha^i$ with 
 \[E_i(x)   = - \frac{1}{2} \left( 
   \ln\left(\frac{\cosh(a)-1}{2}\right) - \ln\left(\frac{\cosh(b)-1}{2}\right)
   -\ln\left(\frac{\cosh(c)-1}{2}\right) \right) . \]
where  $\alpha^i(x)=A$, $\alpha_i \in t$, $d^t(x) = \{A,B,C\}$,  and the $a$,$b$, and $c$ determined as above.
 \end{formula}

Now let's look at  the conformal class of a point 
$x \in \frak{N}$, and define

\begin{definition} 
Let
\[\frak{N}_x = (x + C) \bigcap \frak{N}\]
and call it the conformal class of $x$.
\end{definition}\label{class}
Recall that  $T_p(\frak{N}_x) = span \{w_e \mid e \in E -\partial E\}$ 
from definition \ref{condef}.
Let $e$ be the  edge between $t_1$ and $t_2$  
with $d^{t_i} = \{A_i,B_i,C_i\}$,
$A_i$ corresponding to the angle slot across from $e$ 
in $t_i$, and the $a_i$, $b_i$ , and $c_i$ determined as above.  
Then we may observe... 

\begin{observation}\label{disobse}
 From the above formula and the description of 
$T_p(\frak{N}_x)$
a critical point $y$ of the energy when restricted to 
$x$'s conformal class
satisfies 
\[0 =  dE^y(v_e) = \ln \left(\frac{\cosh(a_{2})-1}{2}\right) - 
\ln \left(\frac{\cosh(a_{1})-1}{2}\right),\]
for all edges $e$.
Hence  at such a critical point we have that $a_{1} = a_2$
and  $y$ is a uniform angle system. 
\end{observation}
   
With this observation the the existence and 
uniqueness of such critical points is equivalent 
to the existence and uniqueness of uniform 
structures conformal to a fixed one.
In particular we may prove the  uniqueness assertion 
in  theorem \ref{tri1}. 
\begin{lemma}\label{lemer}
If $\frak{N}_x$ contains a uniform angle system this angle system is unique in
$\frak{N}_x$. 
\end{lemma}
{\bf Proof:}
This will follow form the 
fact $E$ is concave down.  What we are really asking for is 
that the energy's  Hessian is negative
in the $C$  directions at all points in $\frak{N}_{x}$.
In fact we will show something considerablely stronger, namely that  
the Hessian is negative throughout all of the open subset of
$\Bbb R^{3F}$ satisfying 
$\frak{N}$'s open condition and in all of $\Bbb R^{3F}$'s directions.  
(A fact 
which has several interesting application to the production of triangulations with special symmetries.)
%(see \cite{Le}).
Note that in this setting $E$ and hence its 
Hessian splits up into a sum of independent functions associated
to each triangle - so we only need to show the $3 \times 3$ 
matrix corresponding to a fixed triangle is negative definite.  
Assume we are in the triangle 
with angle slots $\alpha$, $\beta$, and $\gamma$ and $d^t(x) = \{A,B,C\}$. 
Recalling that a linear change of coordinate 
will not effect whether the Hessian is negative definite or not, I
 found it useful to use the coordinates
$(A+B,A+C,B+C)$ for this computation.
To write down the Hessian in these coordinates  it 
is use full to introduce the functions

        \[ F(x,y,z) = \frac{\cos \left(\frac{-x+y+z}{2}\right)}
        {\sin\left(\frac{x}{2}\right)\cos\left(\frac{x}{2}\right)} \]
   and     
        \[ G(x,y,z) = \frac{\cos\left(\frac{x+y}{2}\right)
        \sin\left(\frac{x+y}{2}\right)
        \cos\left(\frac{-x+y+z}{2}\right)
        \cos\left(\frac{x-y+z}{2}\right)}
        {\cos\left(\frac{x+y-z}{2}\right)
        \cos\left(\frac{x}{2}\right)\cos\left(\frac{y}{2}\right)
        \sin\left(\frac{x}{2}\right)\sin\left(\frac{y}{2}\right)} .\]
 With these named a direct computation shows       
        \[Hess = 
        \frac{-1}{2 \cos\left(\frac{A+B+C}{2}\right)}
        \left( \begin{array}{lll}  G(B,C,A) &  F(C,B,A) & F(B,C,A) \\ 
F(C,B,A)  &  G(A,C,B) & F(A,B,C) \\ F(B,C,A) & F(A,B,C) & G(A,B,C)
 \end{array}\right). \]

  To see this matrix is negative definite we can use
  the following easily derived 
  condition: a symmetric  3 by 3 matrix $a^{i}_{j}$ is  
  negative definite if $a^{3}_{3}< 0$,
$ a^{3}_{3}a^{2}_{2} - (a^{2}_{3} )^2 > 0$, 
   and 
  \[ a^{3}_{3}a^{2}_{2}a^{1}_{1}+2a^{1}_{2}a^{2}_{3}a^{1}_{3}-a^{3}_{3}
  (a^{1}_{3})^2 
  -a^{1}_{1}(a^{2}_{3})^2 -a^{2}_{2}(a^{1}_{3})^2  >0. \]
   Note through out this computation that the assumption that $k^t(x) <0$ or rather 
$A +B+C < \pi$
implies all the $\sin$ functions 
 are evaluated at positive angle
   less  than $\frac{\pi}{2}$ and all the  $\cos$ 
   functions are evaluated at sums of angle of absolute value 
   less than $\frac{\pi}{2}$; so all such evaluations are  positive.
    In particular the  negative sign outside guarantees that 
   all terms in the matrix, including the needed $a^{3}_{3}$ term, 
are negative.

  Now we need to satisfy the remaining two conditions; namely we'd like

    \[ G(A,B,C) G(A,C,B)-(F(A,B,C))^2 =\]
  \[\frac{\cos\left(\frac{A+B+C}{2}\right)
  \sin\left(\frac{A+B+C}{2}\right) 
  \left(\cos\left(\frac{-A+B+C}{2}\right)\right)^2}
  {\sin\left(\frac{A}{2}\right)\sin\left(\frac{B}{2}\right)
  \cos\left(\frac{A}{2}\right)\cos\left(\frac{B}{2}\right)}\]
  to be positive, which by the observations made above it is.

  As well as  needing 
  \[ G(A,B,C) G(A,C,B) G(B,C,A) + 2 F(C,B,A) F(B,C,A) F(A,B,C) 
  \]
  \[-G(A,B,C) (F(C,B,A))^2 -  G(B,C,A) (F(A,B,C))^2 - G(A,C,B) (F(B,C,A))^2 =\]
  \[ \frac{32  \left(\cos\left(\frac{A+B+C}{2}\right)\right)^2
   \cos\left(\frac{-A+B+C}{2}\right)
  \cos\left(\frac{A-B+C}{2}\right)
  \cos\left(\frac{A+B-C}{2}\right) }
  { \sin\left(\frac{A}{2}\right)
   \sin\left(\frac{B}{2}\right) \sin\left(\frac{C}{2}\right)}\]
  to be positive which once again as observed above is true.
  
  So indeed the volumes Hessian is  negative definite 
throughout the open subset of $\Bbb R^{3F}$ satisfying 
$\frak{N}$'s open conditions.

   \qed

From this we will be able to show the uniqueness 
statement claimed in the previous section.

\begin{corollary}\label{idpun}
An ideal disk pattern is uniquely determined by its topological 
ploygonal decomposition and associated angle discrepancy system $p$.  
An ideal disk pattern can be constructed form a 
polygonal   
decomposition $\frak{P}$ and associated  
angle discrepancy system $p$ 
precisely when 
there  is a critical 
point of $E$ in  $ \Psi^{-1}(\hat{p}) \bigcap \frak{N}$. 
\end{corollary}
{\bf Proof:}
Notice the choice of triangulation in forming $\frak{P}^t$ from $\frak{P}$ 
at this step is irrelevant 
for both uniqueness and existence, since in an actual pattern 
any edges with a $\theta^e(\hat{p}) = \pi$ will not be 
included in the polygonal decomposition's description.

\begin{figure*}
\vspace{.01in}
\hspace*{\fill}
\epsfysize = 1in  
\epsfbox{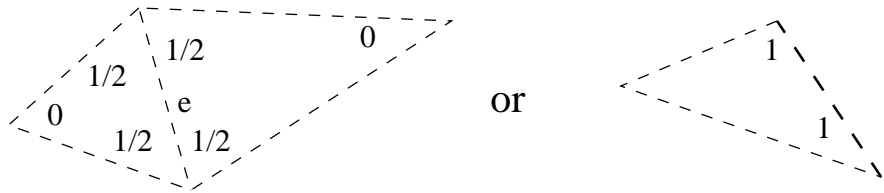}
\hspace*{\fill}
\vspace{.01in}
\caption{\label{ve2} The $m_e$ Vector}
\end{figure*}

Using such a triangular decomposition note 
$\Psi$ has rank $E$ since the pairing of $\psi^{e_i}$ with 
the vector $m_{e_j}$  in figure \ref{ve2} 
satisfies   $\Psi(m_{e_j}) = e_j$ for each $j$.  Further note from  
section 2.1 
that the null space contains the 
$E  - \partial E$ dimension space $C$ and is 
\[ 3F - E = 2E - \partial E - E= E -\partial E\]
dimensional - so $C$ is precisely the  null space.
In particular all angle systems 
which could conceivably hit a specified set of  discrepancy 
angles $\{\psi^e(\hat{p})\}$ is in
$\Psi^{-1}(\hat{p})$, which is $x +C $ for some $x$.
So the above lemma gaurentees the uniqueness and existence of the geodesic 
triangulation necessary to construct the ideal disk pattern 
under the above  conditions.

%(It is worth noting that the surface with  boundary case requires the 
%fact that the $3F = 2E_{int} + E_{\partial}$
%and the observation that the conformal class 
%$C$ in this case is $E_{int}$ dimensional.)
\qed

Now its time to explore the existence  of critical points.
Given a pre-compact open set $O$ 
and a 
continuous function $E$ on $\bar{O}$
we automatically achieve a maximum.  For this maximum to be a 
critical point it is enough to know that 
$E$ is differentiable in $O$ and 
that the point of maximal $E$ is in the open set
$O$.   

One way to achieve this is to show that for any boundary point $y_0$ 
that there is a direction $v$,
 an $\epsilon >0 $ and a $c >0$  such that $l(s) = y_0 + s v$ satisfies 
\[ E(l(0,\epsilon)) \subset \frak{N}\] 
and  
\[ \lim_{s \rightarrow 0} \frac{d}{ds}E(l(s)) \mbox{  } > \mbox{  } c, \]
for all $s \in (0,\epsilon)$.  
This works because under these hypothesize   $E(l(s))$ is continuous and 
increasing on $[0,\epsilon)$ and $y_0$ certainly could not have 
been a point where $E$ achieved its maximum.

%Then $E$ has a critical point if and only if for every such line satisfying  $l(s) \bigcap  O \neq \phi$\begin{eqnarray*}\lim_{s \rightarrow 0} \frac{d}{ds}E(l(s)) & < & 0 \\\lim_{s \rightarrow 1} \frac{d}{ds}E(l(s))&  > & 0. \end{eqnarray*}or equivalently simply that \[ \lim_{s \rightarrow 0} \frac{d}{ds}E(l(s)) \mbox{  } < \mbox{  } 0, \]since the other condition follows from reversing the roles of $y_0$ and $y_1$. 

It is useful to note that the  compactness of $\bar{O}$ gaurentees us that 
$l(s)$ eventually hits the boundary again at $y_1$ for some $s > 0$.  
So we may change the speed of our line and assume  
\[l(s) =
    (1-s)y_0 +sy_1\]
is the line connecting 
the two  boundary points. 
So the remainder of theorem \ref{tri1} follows by applying this above criteria 
to $\frak{N}_x$ and $E$ and noting...

\begin{lemma}\label{push}
For every pair of points $y_0$ and $y_1$ in $\partial \frak{N}$  
but not in $B$
with $l(s) \bigcap \frak{N} \neq \phi$ we have  
\[ \lim_{s \rightarrow 0} \frac{d}{ds}E(l(s))   
\mbox{  } =   \mbox{  }   \infty. \]
\end{lemma} 
{\bf Proof:}
By using $k^t(x)$ and the previous lemma's 
notation for the angles in a triangle $t$, we can write the equation for   
the dual of  $\alpha$'s  
coefficient $E_{\alpha}(x)$ in formula \ref{dif}  as  
      \[-\ln \left|  \frac{\cos(A - k^t(x)) - \cos(A)}
      {(\cos(B-k^t(x)) - \cos(B))(\cos(C- k^t(x)) - \cos(C))} \right|
 + 2\ln(\sin(A)). \]
Fixing a triangle $t$ let $d^t(y_i) = \{A_i,B_i, C_i \}$; 
and note with this notation that the contribution to
$\frac{d E(l(s))}{ds}$ coming from the triangle $t$ is given as 
\[ D^t(s) = \frac{d V^t(l(s))}{ds} \]
\[ = (A_1 - A_0)E_{\alpha}(l(s)) +
  (B_1 - B_0)E_{\beta}(l(s)) +  (C_1 - C_0)E_{\gamma}(l(s)) . \]

Since $y_0$ is on the boundary of 
$\frak{N}$ and not in $B$ there is some triangle $t$ such that   
$\{A_0,B_0,C_0\} \neq \{0,0,\pi\}$ however either $k^t(y_0) =0$ or
some angle is zero.

We will show that for any triangle $t$ in this case 
$\lim_{s \rightarrow 0} D^t(s) =  \infty $.
Its useful to divide the possibilities into the following 
three cases.

\begin{enumerate}
\item
Where  $d^t(y_0)$ contains zeros but  $k^t(y_0) \neq 0$. 
\item
Where $k^t(y_0)=0$ and  no angle is zero.
\item
Where $k^t(y_0) = 0$  and one angle in $d^t(y_0)$ is zero.
\end{enumerate}

In the first case note that the only pieces of $D_t$ 
which become infinite are 
of the form $(A_0 - A_1)2 \ln|\sin(A)|$ and further note 
that if  $A_0 =0$ and $l(s) \bigcap \frak{N} \neq \phi$ then $A_1-A_0 > 0$.  
So we indeed  have $\lim_{s \rightarrow 0}D_t(s) =  \infty$ as required.

To understand the second case note that  if $k^t(x)$ tends to zero then 
then we may rewrite $L_A$ as 
\[ L_A =  -3 \ln(\sin(A))  +  \ln(\sin(B)) +  \ln(\sin(C))  + \ln|k^t(l(s))| .   \]
From this we find  the part of $D_t(s)$ that is not bounded is in the form 
\[ - (A_0+B_0+C_0 -(A_1+B_1+C_1)) \ln|k^t(l(s)|. \]
Further note  when  $k^t(y_0) = 0$ and $l(s) \bigcap \frak{N} \neq \phi$ 
that $ k^t(y_1) =A_1+B_1+C_1 -( A_0+B_0+C_0 ) < 0$.
So once again we have 
$\lim_{s \rightarrow 0}D^t(s) =  \infty$ as required.

The final case is a combination of the above two were we
find  the part of $D_t$ that is not bounded is in the form 
\[ D^t(s) = - (A_0+B_0+C_0 -(A_1+B_1+C_1)) (\ln|k^t(l(s))| + \ln(\sin(A))
 + 4(A_0 - A_1) \ln(\sin(A)). \]
The arguments above immediately imply the correct derivative behavior.

Now note that $\lim_{s \rightarrow 0}D^t(s)$ is clearly 
bounded on  triangles with angles 
not not satisfying any boundary conditions, hence if it were  
bounded or $\infty$ on triangles where $d^t(x) = \{\pi,0,0\}$ we would 
be done.  
The whole reason the $B$ is 
bad set is that it is in fact finite.  
Using the same argument as above and Taylor expanding you 
find
\[ D^t(s) =  (4k^t(y_1) - 4k^t(y_0))\ln|s| + g(s) = g(s), \]
with $g(s)$ bounded.

So the proof is complete.
\qed

%So $D_t(s)$A_0+B_0+C_0 > A_1+B_1+C_1$.  Note
%that $\pi = A_0+B_0+C_0 > A_1+B_1+C_1$ if $k_t$ tends to 
%zero as $s$ approaches zero and similarly 
% $\pi = A_1+B_1+C_1 > A_0+B_0+C_0$ if $A$ tends to zero 
%as $s$ approaches one - once again as needed.

%Namely we can eliminate the possibility that the angles in a triangle
% deforms to $\{ \pi, 0,0\}$ by noting that  if they did then  we would have
%$\pi +  \frac{\pi + A -B -C}{2} = \theta_l < \pi$ forcing $\pi+A -B -C < 0$
%and hence $A+B+C - \pi > 0$ - forcing us to 
%have already left the set negative curvature angle arrangements.  

%This is the only use of conformality so as above we now simply assume that
%we are in the 
% region where $k_t < 0$ and all angles are in $(0,\pi)$; with the 
%added criteria  that 
%$n_0$ and $n_1$ are on the boundary of this region with no triangle's angles 
%are deforming to the eliminated condition.  
%As in the above lemma  we can work triangle by triangle, and use the above notation 

%if $A_0 > A_1$ and $- \infty$
% if $A_1 >A_0$.  Note that $A_1 > A_0$ if $A$ tends to 
%zero as $s$ approaches zero and $A_0 > A_1$ if $A$ tends to zero 
%as $s$ approaches one - exactly as claimed.

\end{subsection}

\end{section}
\begin{section}{A Warm Up Thurston-Andreev Theorem}\label{yipyeep}
In this section we will prove a warm up Thurston Andreev theorem.
This section and the next  have been made independent of 
each other so there is a bit of repetition.

Corollary   \ref{lemer} tells us that an ideal disk pattern
is always unique when it exists, and we are now left to deal with the 
dilemma of finding good existence 
criteria. 
Here I will describe in detail  the strict convex case where 
$p \in (0,\pi)^{(E - \partial E)}  \times (0,\frac{\pi}{2})^{\partial E}$ relative to  $\frak{T}$
a triangular decomposition.  
Linear conditions on the possible $p$ will be produced   
which are  necessary and 
sufficient for the $p$ relative to $\frak{T}$ to be the data 
of an ideal disk pattern.

Let
\[\frak{D} = \left\{ x \in \frak{N} \mid  
\Psi(x) \in  (0,\pi)^{(E - \partial E)}  
\times \left(0,\frac{\pi}{2} \right)^{\partial E} \right\}. \]
and call this the set of negative curvature 
Delaunay angle systems. 
These angle systems   are remarkably 
easy to work with such angle systems and in fact...

\begin{observation}\label{delnice}
Every point of $\frak{D}$ has unique ideal 
disk pattern associated to it.
\end{observation}
{\bf Proof:}  
This observation relies on the following fact which 
will be of interest in its own right. 

\begin{fact}\label{thing}
If $x \in \frak{D}$ is conformally equivalent to 
a point in $\partial \frak{N}$
where for some triangle
$d^t(x) = \{0,0,\pi\}$.
\end{fact}

To see this  fact assume to 
the contrary that for some $t$  and $c$ we have 
$d^t(x + c) = \{0,0,\pi\}$. 
Let $e$ be the  edge of $t$ across form 
$t$'s  $\pi$
and let $t_1$ be $t$'s neighbor next to $e$ if it exists.
Note by fact \ref{ob} that
that the conformally invariant $\psi^e(x) \in (0,\pi)$ would (even 
in the best possible case when $e$ is not on the boundary) 
have 
to satisfy
the contradictory inequality 
$\psi^e(x+C) = - \frac{\pi}{2} + \psi_{t_1}^e \leq 0 $.

From this fact we have 
that if $x \in \frak{D}$ then no element in 
$x +C$ could possibly be in $B$ and the observation follows form 
 theorem   
\ref{tri1}.

\qed

Now lets explore certain two
necessary conditions on a $p = \Psi(x)$ with $x \in \frak{D}$.   
The first condition 
is the condition 
related to the fact that the angles at the internal 
vertex in a geometric 
triangulation sum to $2 \pi$ and at a  boundary vertex sum to $\pi$.

\[ \mbox{(}n_1\mbox{)        } \left\{ \begin{array}{ll}
\sum_{e \in v} \psi^e = 
2 \pi   &   \mbox{  if   } v\in V - \partial V \\
 \sum_{e \in v} \psi^e = 
\pi   &   \mbox{ if  }  \partial V
\end{array} \right. \]

%\begin{eqnarray*}
%(n_1)  & & \sum_{e \in v} (\psi^e(p)) = 
%2 \pi   &  & \mbox{ such that } v\in V - \partial V. \\
%\end{eqnarray*}
%\begin{eqnarray*}
%\mbox{  }   & & \sum_{e \in v} (\psi^e(p)) = 
%\pi   &  & \mbox{ such that }  \partial V. \\
%\end{eqnarray*}

This condition  is 
equivalent to the following simple lemma.

\begin{lemma}\label{n11}
\[ \Psi(V) = \{p \in \Bbb R^{E} \mid p \mbox{ satisfies } (n_1)\}.\]
\end{lemma}
{\bf Proof:}
First note that if $p = \Psi(x)$  then 
\[ \sum_{e_i \in v} e^j(p) =  \sum_{e_i \in v} \psi^{e_i}(x) =
p^{v} (x).  \]
So by choosing $x \in V$ we see $\Psi(V)$ is included in 
\[ W = \{p \in \Bbb R^{E} \mid p \mbox{ satisfies } (n_1)\} .\]
Recall form the proof of corollary \ref{idpun} that
$\Psi(\Bbb R^{3F}) = \Bbb R^E$. So 
we may express any $p \in W$ as $p = \Psi(x)$ 
and  the above computation  gaurentees $x \in V$ as needed. 

\qed

The second necessary  condition is a global one; namely 
an  insistence that  for every set $S$ of
$|S|$  triangles  in $\frak{T}$ that 
\begin{eqnarray*}
(n_2)  & &   \sum_{e \in S}  \theta^e > \pi |S|.  \\
\end{eqnarray*}

Verifying  $(n_2)$ relies on the following formula.
\begin{formula}\label{sett}
Given a set of triangles $S$ 
\[ \sum_{\{e \in S \}} \theta^e(x)   = 
\sum_{t \in S} \left( \pi  - \frac{k^t(x)}{2}\right) 
+ \sum_{e \in \partial S - \partial E}  
\left(\frac{\pi}{2} - \psi^e_{t}(x)\right) , \]
with the   $t$ in $\psi^e_{t}(x)$ term being the triangle on 
the  non-$S$ side of $e$.   
\end{formula}
{\bf Proof:}
\[ \sum_{\{e \in S \}} \theta^e(x) = 
\sum_{e \in S - \partial E} (\pi - \psi^e(x)) 
+ \sum_{e \in S \bigcap \partial E} 
\left(\frac{\pi}{2} - \psi^e(x) \right) \]
\[ =
\sum_{e \in S - \partial E} 
\left( \left(\frac{\pi}{2} - \psi_{t_1}^e(x)\right) 
+ \left(\frac{\pi}{2} - \psi_{t_2}^e(x)\right)\right) 
+ \sum_{e \in S \bigcap \partial E} 
\left(\frac{\pi}{2} - \psi^e(x) \right)\] 
\[ = \sum_{t \in S} \left( \pi  + \frac{ \pi - l^t(x)}{2}\right) 
+ \sum_{e \in \partial S - \partial E}  
\left(\frac{\pi}{2} - \psi^e_{t}(x)\right)  \]
with the   $t$ in $\psi^e_{t}(x)$ term being the triangle on 
the  non-$S$ side of $e$.   Substituting the definition of 
$k^t(x)$ gives the needed formula.

 \qed

Note for any point $x \in \frak{N} $ that 
 $- k^t(x) > 0$ and form  fact \ref{ob} 
 that $\frac{\pi}{2} - \psi^e_{t}(x) >0$.  
So removing these terms form the above formula strictly 
reduces its size and when summed up we arrive at $(n_2)$. 

With these two necessary condition 
 we have our first pattern existence theorem:

\begin{theorem}\label{cir11}
If 
\[ p \in D  = \left\{q \in (0,\pi)^{(E - \partial E)}  
\times \left(0,\frac{\pi}{2}\right)^{\partial E} \mid  q \mbox{ satisfies } 
(n_1) \mbox{ and } (n_2) \right\} \] 
then $p$  is realized by  a unique 
ideal disk pattern.
\end{theorem}

By observation \ref{delnice} above 
%and the fact that all the condition in $D$ are necessary  
this would follow if 
we knew the following proposition. 

\begin{proposition}
\[\Psi(\frak{D}) = D.\]
\end{proposition}

It is this bit of linear algebra we now will tackle.  
Notice the fact $(n_1)$ and $(n_2)$ are 
necessary gaurentees that $\Psi(\frak{D}) \subset D$, 
and we left 
to explore  $\Psi$'s surjectivity.

\begin{subsection}{The Surjectivity of $\Psi$: the Delaunay Case}
To see the surjectivity of $\Psi$ 
let's assume the contrary that that $\Psi(\frak{D})$ 
is strictly contained in $D$ and produce a contradiction.
With this assumption we have
 a point $p$ on the boundary of $\Psi(\frak{D})$  
inside $D$.
Note $p = \Psi(y)$ for some 
$y \in \partial \frak{D}$.  Furthermore note
$(C+y) \bigcap \frak{D}$ is empty, since other wise for some 
$w \in C$ we would have $(y + w) \in \frak{D}$ which along with the fact
that $\Psi$ is an open mapping when restricted to $V$  
would force  $p = \Psi(y) =\Psi(y +w)$ to be 
in the interior of $D$.

At this point we need to choose a particularly nice conformal version of 
$y$, which requires the notion of a 
stable boundary point of $\frak{D}$.  
Before defining stability note since  $\frak{D}$ is a convex 
set with hyperplane boundary if 
$x \in \partial \frak{D}$ such that  
$(x + C) \bigcap \frak{D} = \phi$, then 
$(x + C) \bigcap \partial \frak{D} $
is its self a convex $k$ dimensional set.    

\begin{definition} \label{stable1}
A point in $x \in \partial \frak{D}$ is stable if 
$(x + C) \bigcap \frak{D} = \phi$  and $x$ is in the interior of 
$(x +C) \bigcap \partial \frak{D} $ as a $k$ dimensional set.
Any inequality forming  $\frak{D}$ violated in order to
make $x$ a boundary point  will be 
called a violation.
\end{definition} 

The key property of a stable point is that 
a conformal change $w \in C$ 
has $x+\epsilon w \in \bar{\frak{D}}^c$  for all $\epsilon > 0$ 
or for some  sufficiently small  $\epsilon > 0$ we have 
$x + \epsilon w$ must still be on $\partial{\frak{D}}$ and 
experience exactly the same violations as $x$.
  The impossibility of any other
phenomena when conformally changing a stable point is at the heart of the 
arguments in lemma \ref{sub1} and lemma \ref{sub21} below.
At this point subjectivity would follow if
for a stable  $x \in \partial {\frak{D}}$ 
we knew that  $\Psi(x)$ could not be in  $D$, contradicting the choice of 
$p = \Psi(x)$ as needed. 
 
We will prove this by splitting  
up the possibilities into the two cases in
lemma \ref{sub1} and lemma  \ref{sub21}.

\begin{lemma}\label{sub1}
If $x \in \partial {\frak{D}}$ is stable
and $\alpha^i(x) =0$ for $\alpha^i$ in a triangle where 
$k^t(x) <0$, then $\Psi(x)$ is not in $D$.
\end{lemma}
{\bf Proof:}
Look at an  angle slot which is zero in  triangle $t_0$ satisfying  
$k^{t_0}(x) < 0$.
View this angle as living  
between the edges $e_0$ and $e_1$.  
Note that in order for $x$ to be stable that  
either $e_1$ is a boundary edge  or 
the $\epsilon w_{e_1}$ transformation 
(with its positive side in $t_0$) 
 must be protected 
by  a zero on the  $-\epsilon$   side forcing the condition that 
$\epsilon w_{e_1} \in \bar{\frak{D}}^c$ , or else for small enough 
$\epsilon$ we would have 
$x + \epsilon w_{e_1}$  would be a conformally 
equivalent point on $\partial{\frak{D}}$  with 
fewer violations. When $e_1$ is not a boundary edge
call  this neighboring triangle $t_1$ and when it is a 
boundary edge stop this process.
If we have not stopped let  $e_2$ be another 
edge bounding a zero angle slot in $t_1$ and stop if it is a boundary edge.
If it is not a boundary edge  
then there are two possibilities.
If
$k^{t_1} < 0$ repeat the above procedure letting $e_1$ play the role of
 $e_0$ and $e_2$ the role of $e_1$ and constructing an $e_3$ in a 
triangle $t_2$.
If  $k^{t_1}(x) = 0$ conformally change $x$ to 
\[x +  \epsilon w_{e_1} + \epsilon w_{e_2}.\]
Notice  no triangle with 
$k^t(x)=0$ can have two zeros by fact \ref{thing}, so
for the initial zero violation to exist 
there most be a zero on the $-\epsilon$ side of 
$\epsilon w_{e_2}$.  Once again we have determined an $e_3$ and $t_2$.

Using this procedure  
to make our decisions  we may continue this process  
forming a set of edges 
$\{e_i\}$ with the angle between $e_i$ and $e_{i+1}$, $A^{i,i+1}(x)$, always 
equal to zero. 
Since there are a finite number of edges either we stop at a boundary edge or 
eventually in this 
sequence will have some $k<l$ such that $e_k = e_l$ and $e_{k+1} =e_{l+1}$.
(This by the pigeon hole principle since 
some edge $e$ will appear an infinite number of times in this 
list and among its infinite neighbors there must be a repeat).

In the case  the sequence never stops we can produce a contradiction.
To do it first note if  $e_{i}$ and $e_{i+1}$ are in $t_i$ then
$A^{i,i+1}(x) = \psi_{t_i}^{e_i} + \psi_{t_i}^{e_{i+1}}$.  So for 
the set of edges $\{e_i\}_{i=k}^{l-1}$ we have 
\[ 0 = \sum_{i = k}^{l-1} A^{i,i+1}  = \sum_{i=k}^{l-1} \psi^{e_i}(x) > 0 \]
our needed contradiction. 

In the case the sequence did hit the boundary 
perform the construction  in the opposite direction.  If we don't 
stop in this direction we arrive at the same contradiction.  If we did then
this computation still produces a contradiction on the path with the 
two boundary edges, since for a boundary edge in the triangle $t$ we 
have $\psi^e_t(x) = \psi^e(x) \in \left(0,\frac{\pi}{2}\right)$. 
%which  is assumed in $\$ for all edges.

\qed

\begin{lemma}\label{sub21}
If a stable $x$ satisfies the condition that 
if $\alpha^i(x) =0$ then $\alpha_i$ is in a 
triangle $t$ 
with
$k^t(x) = 0$, 
then $\Psi(x)$ is not in $D$. 
\end{lemma}
{\bf Proof:}
In this case, in order for $x$ to be a boundary point of $\frak{D}$ 
for some $t$ we have that 
$k^t=0$. We will be looking at the nonempty set of 
all triangles with $k^t = 0$, $Z$.  
The first observation needed about $Z$ is that it is
not all of $M$ and has a non-empty internal boundary 
(meaning $\partial Z - \partial M$).
To see this note 
\[\sum_{t \in \frak{P}} k^t(x) = \sum_{e \in v} A^i -\pi F 
 = \pi {\partial} V + 2 \pi (V - \partial V) -  \pi F \]
\[  =
2 \pi V - (\pi {\partial}  V+ 3 \pi F) +  2 \pi F  
=   2 \pi V - 2 \pi E + 2 \pi F  = 2 \pi \chi(M) < 0,\]
so there is negative curvature somewhere.

By the stability  
of $x$ once again 
there can be no conformal transformation capable of moving negative 
curvature into this set. Suppose we are at an internal boundary $e_0$
edge of $Z$, call the triangle  on the $Z$ side of 
the boundary edge $t_0$ and the triangle on the
non-boundary edge  $t_{-1}$.      Since  $t_{-1}$ has negative 
curvature the obstruction to 
the $\epsilon w_{e_0}$ transformation
being able to move 
curvature out of $Z$ must be due to  $t_0$.  
In order for $t_0$ to protect against 
this there must be zero along $e_0$ on the $t_0$ side.

Now we will continue the attempt to suck
curvature out with a curvature vacuum.
Such a vacuum is an element of $C$ indexed by a set of $Z$ edges.  
The key observation in forming this vacuum is once again fact 
\ref{thing} telling us
if an angle in $t$ is zero and $k^t(x) = 0$ 
then there is only one zero angle in $t$.  
Let $e_1$ be the other edge sharing the unique zero angle along 
$e_0$ in $t_0$ and if 
$e_1$ is another boundary edge we stop. If $e_1$ is not a boundary edge
 use $w_{e_1} + w_{e_0}$ 
to continue the effort to remove curvature.  Continuing this 
process forms a completely determined
set of edges  and triangles, $\{e_i\}$ and $\{t_i\}$,   and 
a sequence of conformal   
transformations $\epsilon \sum_{i =0}^n w_{e_i} \in C$.

We will now get some control over this vacuum.  Note a vacuum never 
hits itself since if there is a first pair $k <l$ such that
$t_k =t_l$ then $t_k$ would have to have to have two zero 
and zero curvature, which fact \ref{thing} assures us is  impossible. 
So any vacuum hits a boundary edge or  
pokes  through $Z$ into $Z^c$.

In fact with this argument we can arrive at 
the considerablely sponger fact 
that two vacuums can never 
even share an edge.  
To see this call a  vacuum's side
boundary any edge  of a triangle in the vacuum 
facing a zero.  
Now simply note if the intersection of  two vacuums contains an edge then
it contains a first edge $e_i$ with respect to one of the vacuums.
There are two possibilities for this edge. One is that $t_{i+1}$ has two zeros 
and $k^t(x) =0$, which we showed was impossible in the previous paragraph.
The other is that $e_i$ is a side boundary of both  vacuums.
In this case we have an edge facing zero angles in both directions in triangle with zero curvature,  so this would force $\psi^e(x) = \pi$, a contradiction. 
So either case 
is impossible, and indeed no distinct vacuums share an edge. 
 
Let 
$S$ be the removal
from $Z$ of all these vacuums.  First I'd like to note that 
$S$ is non-empty.  Note every vacuum has  
side boundary.  Since vacuums cannot intersect themselves or 
share edges with distinct 
vacuums,  $S$ 
would be nonempty if side boundary had to be 
in $Z$'s interior.  
Look at any side boundary edge $e$ of a fixed vacuum. Note $e$ 
cannot be on $\partial Z - \partial M$ since then the vacuum 
triangle it belonged to 
would have at least two zeros  and $k^t(x) =0$.  
Furthermore 
$e$ cannot be on $\partial M$ 
since then $\psi^e(x) = \frac{\pi}{2}$.  
So indeed $S$ is nonempty.

Now lets observe the following formula.
\begin{formula}
Given a set of triangles $S$ 
\[ \sum_{\{e \in S \}} \theta^e(x)  = \sum_{e \in \partial S-\partial M}  
\left(\frac{\pi}{2} - \psi^e_{t}(x)\right)  = 
\sum_{t \in S} \left( \pi  - \frac{k^t(x)}{2}\right) 
+ \sum_{e \in \partial S-\partial M}  
\left(\frac{\pi}{2} - \psi^e_{t}(x)\right) , \]
with the   $t$ in $\psi^e_{t}(x)$ term being the triangle on 
the  non-$S$ side of $e$.   
\end{formula}
{\bf Proof:}
\[ \sum_{\{e \in S \}} \theta^e(x) = 
\sum_{e \in S} (\pi - \psi^e(x)) \]
\[=
\sum_{e \in S-\partial M} \left( \left(\frac{\pi}{2} - \psi_{t_1}^e(x)\right) 
+ \left(\frac{\pi}{2} - \psi_{t_2}^e(x)\right)\right)  
+\sum_{e \in \partial S} 
\left(\frac{\pi}{2} - \psi_{t}^e(x)\right) \] 
\[ = \sum_{t \in S} \left( \pi  + \frac{ \pi - l^t(x)}{2}\right) 
+ \sum_{e \in \partial S -\partial M}  
\left(\frac{\pi}{2} - \psi^e_{t}(x)\right)  \]
\[ =\sum_{t \in S} \left( \pi  - \frac{k^t(x)}{2}\right) 
+ \sum_{e \in \partial S-\partial M}  
\left(\frac{\pi}{2} - \psi^e_{t}(x)\right). \]

 \qed

Now every edge in $\partial S -\partial M$ faces a zero on its $S^c$ side
in a triangle with $k^t(x) =0$, so
\[ \sum_{e \in \partial S-\partial M}  
\left(\frac{\pi}{2} - \psi^e_{t}(x)\right) =0.\]
Similarly each triangle has zero curvature so from the above 
formula we have 
\[   \sum_{\{e \in S \}} \theta^e(x) = |S| \pi \]
violating condition $(n_2)$.  So we have constructed a violation to 
$(n_2)$ and $\Psi(x)$ cannot be in $D$ as need.

\qed

\end{subsection}

It is worth noting that  nothing prevents us from 
extending the main theorem of this section  (and the next) 
to the case were $p^v(x) \neq 2 \pi $ and in particular 
to the version of this theorem where $p^v(x) =0$ and the resulting
uniform 
surface is a finite area hyperbolic surface with cusps.

\end{section}
\begin{section}{The Ideal Thurston-Andreev Theorem}\label{itat}
In this section we prove and sate the general ideal convex 
Thurston-Andreev theorem, which is simply the 
polygonal decomposition case of the theorem of the previous section.  
The proof here is done in  
detail when the surface has no boundary, and dealing with the boundary 
can be accomplished exactly as in the previous section. 
 
\begin{subsection}{The Statement and Reduction to Linear Algebra}\label{srla}

Corollary   \ref{lemer} tells us that an ideal disk pattern
is always unique when it exists, and we are now left to deal with the 
dilemma of finding good existence 
criteria. In fact in the convex case, 
(when $p \in (0,\pi)^E$ relative to $\frak{P}$)  
we will produce linear conditions on the possible $p$   
which are  necessary and 
sufficient for the $p$ relative to $\frak{P}$ to be the data 
of an ideal disk pattern.

The first example of a necessary condition on $p$ 
is the condition 
related to the fact that the angles at the vertex in a geometric 
triangulation sum to $2 \pi$, namely 
\begin{eqnarray*}
(n_1)  & & \sum_{e \in v} (\psi^e(p)) = 2 \pi   . \\
\end{eqnarray*}
This and all the mentioned necessary condition 
will be demonstrated as such in the next section. 
Another necessary condition is a global one 
(though often localizable), namely 
an  insistence that  for every set $S$ of
$|S|$  polygons  in $\frak{P}$ that 
\begin{eqnarray*}
(n_2)  & &   \sum_{e \in S}  \theta^e(p) > \pi |S|.  \\
\end{eqnarray*}

With these two necessary condition 
 we have our first pattern existence theorem:

\begin{theorem}\label{cir1}
That $(n_1)$, $(n_2)$, and  $p \in \Bbb (0,\pi)^E$ hold 
is  necessary and sufficient 
for $p$ relative to $\frak{P}$  to be associated to a unique 
convex  ideal disk pattern.
\end{theorem}

\begin{note} \label{bound4}
This theorem works equally well for surfaces with boundary.
The only modifications is the obvious one that at a boundary vertex 
 $\sum_{e \in v} (\psi^e(p)) =  \pi$. 
\end{note}
Such pattern data is extremely common; for example 
the data determined by  
the circumscribing circles in any random Delaunay triangulation 
of a varying negative curvature surface 
(see section \ref{conu}).
%(see \cite{Le}).
This fact is the main reason I've chosen to use the terminology of  
the disk pattern construction  rather than polyhedra construction 
in this section.

One way to prove theorem \ref{cir1} allows some  understanding 
of the non-convex case as well.  Namely we will find some 
necessary criteria on
angles when $p \in (-\pi, \pi)^E$. 
In order to articulate these conditions we need certain 
 snake and a loop concepts in a triangular decomposition.

\begin{figure*}
\vspace{.01in}
\hspace*{\fill}
\epsfysize = 2in 
\epsfbox{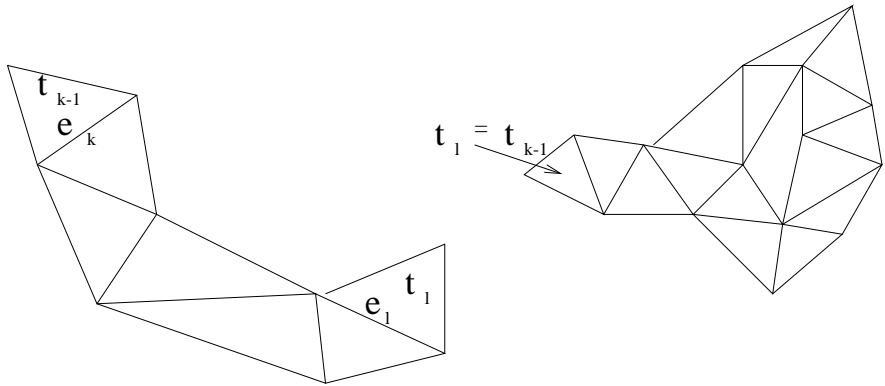}
\hspace*{\fill}
\vspace{.01in}
\caption{\label{the2} A snake and a balloon}
\end{figure*}

\begin{definition}\label{snakes}
  A snake is a finite
directed sequence of edges $\{e_i\}_{i=k}^{l}$
directed in the following sense: if $k <l$ we 
start with the edge $e_k$ between $t_{k-1}$ and $t_k$,  
then we require $e_{k+1}$ to be one of the remaining edges on $t_k$.  
Then letting 
$t_{k+1}$ be the other face associated to $e_{k+1}$ we require $e_{k+2}$ 
to be 
one of the other edges of $t_{k+1}$ and so on until some 
tail edge $e_l$ 
and tail face $t_l$ are  reached, and 
if $l< k$ we reverse the procedure and add rather than subtract from 
the index.  
See figure \ref{the2} for  examples. 
%but note over laps would be indeed allowed. 
A loop is a  snake $\{e_i\}_{i=l}^{k}$ 
where   
$e_k = e_l$ and $t_k = t_l$, see figure \ref{the3} 
for a pair of examples.
%once again keeping in mind the possibilities of over laps.
\end{definition}

\begin{figure*}
\vspace{.01in}
\hspace*{\fill}
\epsfysize = 2in 
\epsfbox{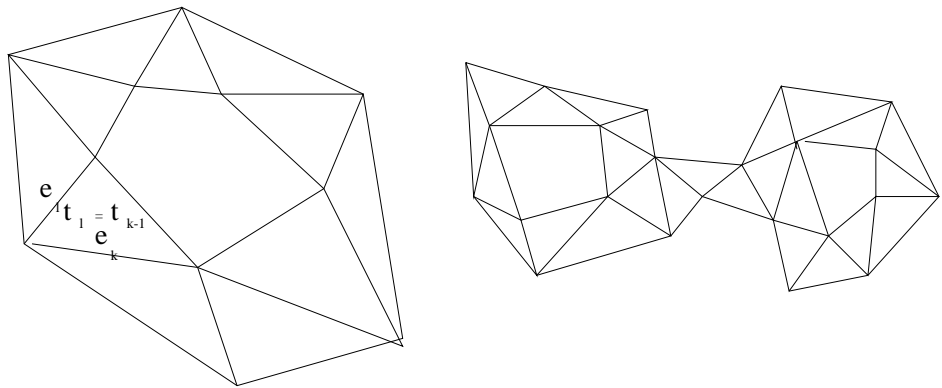}
\hspace*{\fill}
\vspace{.01in}
\caption{\label{the3} A loop and a barbell}
\end{figure*}

It is a condition on 
snakes and loops 
which allows  one to articulate
the remaining necessary conditions.  However as defined 
there are then  an infinite number of such objects 
and it nice to first isolate a 
finite sub-set that does the job.

\begin{definition}\label{barbell}
A set of edges $\{e_i\}_{i=k}^l$  is called embedded if 
$e_i \neq e_j$. A snake $\{e_i\}^l_k$  is said to double back on itself 
if we have a pair of non-empty sub-snakes with
$\{e_i\}_{m}^{n}$ and $\{e_i\}_{k-m}^{k-n}$ containing the same edges.  
A barbell is a loop which doubles back on itself and such that 
$\{e_i\}_{i=k}^l / \{e_i\}_{i=m}^n $ is embedded.  A balloon  is a snake 
which doubles back on itself with 
$\{e_i\}_{i=k}^l / \{e_i\}_{i=m}^n $ embedded and such that $e_l = e_k$.
\end{definition} 

With this terminology the remaining necessary conditions are  
\begin{eqnarray*}
(n_3) & \sum_{i=k}^{l-1} \theta^e_i(p) < |k-l| \pi & 
\mbox{ when } \{e_i\}_{i=k}^l \mbox{ is an embedded loop or barbell,}
\end{eqnarray*}  
and 
\begin{eqnarray*}
(n_4) & \sum_{i =k}^{l} \theta^e_i(p) < (|k-l| +1)  \pi & 
\mbox{ when } \{e_i\}_{i = k}^l \mbox{ is an embedded sake or balloon.}
\end{eqnarray*}

With these conditions let
\[ N = \{p \in (-\pi,\pi)^E \mid  p \mbox{ satisfies } 
(n_i) \mbox{ for each i} \}. \]

To each  ideal disk pattern   we may produce 
an uniform element of a $\frak{N}$ by 
choosing a geodesic $\frak{P}^t$ 
associated to 
the patterns circumscribing $\frak{P}$. 
So the above necessary conditions would follow if 
$\Psi(\frak{N}) \subset N$.
In fact in the next section we shall prove.... 
\begin{theorem}\label{main}
\[\Psi(\frak{N}) = N\] 
\end{theorem}

With this result we are in a position to prove theorem  \ref{cir1}.

{\bf Poof of theorem \ref{cir1} from theorem \ref{main}:}
First we will show that  
$\hat{p} \in N$ relative to  a chosen $\frak{P}^t$ (where we assume polygons have been triangulated as in figure \ref{fan}). 
To do this we need that conditions 
$(n_3)$ and $(n_4)$ are satisfied.  Since 
$\hat{\theta} \in (0,\pi] \mbox{  }(n_4)$ 
is automatic and  $(n_3)$ can only be false if there is a 
loop or barbell $\{e_i\}$ on which 
$\hat{\theta}^{e_i} \equiv \pi$.  
Note that by our choice of $\frak{P}^t$ 
(though any other choice would in fact still work with a 
slight modification) we see that the only snakes $\{e_i\}$ 
containing all $\hat{\theta}^{e_i} = \pi$ edges 
are snakes with edges contained 
in some polygon's fan and in particular can never loop up or form a barbell.   
So by theorem \ref{main} we see the point $\hat{p}$
described in theorem \ref{cir1}
has a preimage  which intersects  $\frak{N}$ non-trivially.

At this point all we need is that a point $y$ in this 
preimage $y+C$ satisfies the conditions of theorem \ref{tri1}.  
Namely we will suppose that $(y + C) \bigcap B$ is not empty, i.e. 
$(y+w) \in B$ with $w \in C$, and produce a contradiction. 
In particular  this assumption gaurentees there is   
some triangle $t_0$ with   
$d^{t_0}(y + w) = \{\pi,0,0\}$.  
Let the edge $e_1$ of $t_0$ be the edge with 
the $\{0,0\}$ of $t_0$ on it.  The fact $d^{t_0}(y + w) = \{\pi,0,0\}$ 
allows us 
to control the $d^{t_1}(y +w)$
data in the other triangle containing the edge $e_1$.
This  follows by observing the natural and 
what will prove be useful decomposition 
of $\psi^{e_1} = \psi^{e_1}_{t_0} + \psi^{e_1}_{t_1}$ as in 
figure \ref{new1}, and the following trivial but useful fact...   
\begin{figure*}
\vspace{.01in}
\hspace*{\fill}
\epsfysize = 2in  
\epsfbox{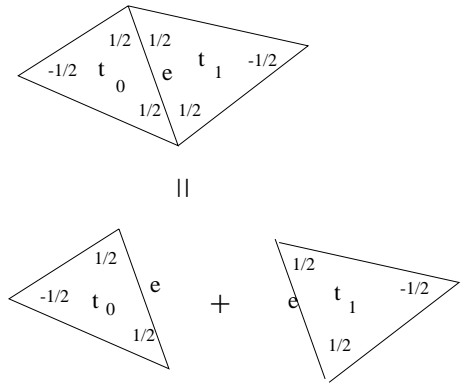}
\hspace*{\fill}
\vspace{.01in}
\caption{\label{new1} The  decomposition $\psi^e = 
\psi^e_{t_0} + \psi^e_{t_1}$}
\end{figure*}
\begin{fact}\label{ob}
When $x \in \frak{N}$ we have
$\psi_{t_i}^e(x) \in \left(\frac{-\pi}{2},\frac{\pi}{2} \right)$, and 
when $x \in \partial \frak{N}$ we have
$\psi_{t_i}^e(x) \in \left[\frac{-\pi}{2},\frac{\pi}{2} \right] $ .
\end{fact}
{ \bf reason for the fact:}

Let  $d^{t_i}(x) = \{A,B,C\} $ and note since 
$B +C \leq A+B+C = l^t(x) < \pi$ and $A<\pi$ we have  
\[ -\frac{\pi}{2} < -\frac{A}{2} \leq  \psi^e_{t_i} = \frac{B+C-A}{2}
\leq \frac{B +C}{2} < \frac{\pi}{2} . \]
The second statement follows from the possibility 
of these inequalities becoming equalities.

Back to the proof. From this fact  we have that  
\[\theta^{e_1}(y + w) = 
\pi + \frac{\pi}{2} - \psi^{t_1}(y +w),\]
forces $\psi^{t_1}(y +w) = \frac{\pi}{2}$,  and hence 
$d^{t_1}(y+w) =\{A,\pi -A,0\}$ with the zero opposite to $e_1$.  
Now the assumption we are in $B$ 
implies $d^{t_1}(x)$ is not legal and $A = \pi$ or $0$.
Let $e_2$ be the $\{0,0\}$  edge $t_1$ and continue 
this argument hence forming a snake of edges with $\theta^{e_i} = \pi$.  
Note by finiteness of the triangulation there must be a 
first $l$ and $k <l$ where  $t_l = t_k$.  
Note when this happens that $\{e_i\}_{k+1}^l$ 
forms an embedded loop with all its $\theta^{e_i} = \pi$,  
contradicting condition $(n_3)$.   
  
So $(y+C) \bigcap B$ is indeed empty and we have our 
need triangulation, hence 
our needed ideal disk pattern.   

\qed

%With this concept we have the another necessary condition
%\begin{eqnarray*}
%(\hat{\hat{n}}_3) & \mbox{ } & \mbox{ There is no loop } \{e_i\}_{i=k}^l \mbox%{ where } 
%\theta(e_i) \equiv  \pi.  
%\end{eqnarray*}

%This  condition allows us to extend $\theta(e)$'s possible values to abstractly to $(0.\pi]$ with...
%\begin{theorem}\label{cir2}
%If a topological triangular decomposition and $\theta$ function with values in$(0,\pi]$ satisfies $(n_1)$, $(n_2)$, and $( \hat{\hat{n}}_3)$ then there exits a unique  solution to the $\theta$-triangular circle pattern with this data. 
%\end{theorem}

%\begin{figure*}
%\vspace{.01in}
%\hspace*{\fill}
%\epsfysize = 1.5in 
%\epsfbox{the1.eps}
%\hspace*{\fill}
%\vspace{.01in}
%\caption{\label{the1} Angle Determination}
%\end{figure*}

%In this section we articulate particularly nice conditions for a polygonal decomposition $\frak{P}$ and angle discrepancy information $p \in \Bb R^E$to be the data $\theta$-ideal disk pattern. 

\end{subsection}
\begin{subsection}{Proof of Theorem \ref{main}}\label{pitat}
\begin{subsubsection}{Injectivity}

Our first goal in proving theorem \ref{main} is to show
\[\Psi(\frak{N}) \subset N.\]

For starters note the fact  
$p \in (-\pi,\pi)^E$ follows immediately 
from fact 
\ref{ob1} 
in the previous section.   Now we need to verify the conditions 
$(n_1)$ through $(n_4)$ hold in $\Psi(\frak{N})$.   
Condition $(n_1)$ is equivalent to the following simple lemma.

\begin{lemma}\label{n1}
\[ \Psi(V) = \{p \in \Bbb R^{E} \mid \sum_{\{ e_i \in v\}} 
e^j(p) = 2 \pi \}.\]
\end{lemma}
{\bf Proof:}
First recall that we know $\Psi$ is surjective, so 
we may express any $p \in N$ as $p = \Psi(x)$.
L`So letting $\{e \in v\}$ denote the set of edges at a vertex $v$ we have 
 \[ \sum_{e_i \in v} e^j(p) =  \sum_{e_i \in v} \psi^{e_i}(x) =
p^{v} (x).\]
So in particular the affine flat  
\[W = \{p \in \Bbb R^{E} \mid \sum_{\{ e_i \in v\}} 
e^j(p) = 2 \pi \} \]
is precisely $\Psi(V)$.  
\qed

Verifying condition $(n_2)$ relies on the following formula.
\begin{formula}\label{set}
Given a set of triangles $S$ 
\[ \sum_{\{e \in S \}} \theta^e(x)   = 
\sum_{t \in S} \left( \pi  - \frac{k^t(x)}{2}\right) 
+ \sum_{e \in \partial S}  
\left(\frac{\pi}{2} - \psi^e_{t}(x)\right) , \]
with the   $t$ in $\psi^e_{t}(x)$ term being the triangle on 
the  non-$S$ side of $e$.   
\end{formula}
{\bf Proof:}
\[ \sum_{\{e \in S \}} \theta^e(x) = 
\sum_{e \in S} (\pi - \psi^e(x)) =
\sum_{e \in S} \left( \left(\frac{\pi}{2} - \psi_{t_1}^e(x)\right) 
+ \left(\frac{\pi}{2} - \psi_{t_2}^e(x)\right)\right)  \] 
\[ = \sum_{t \in S} \left( \pi  + \frac{ \pi - l^t(x)}{2}\right) 
+ \sum_{e \in \partial S}  
\left(\frac{\pi}{2} - \psi^e_{t}(x)\right)  \]
\[ =\sum_{t \in S} \left( \pi  - \frac{k^t(x)}{2}\right) 
+ \sum_{e \in \partial S}  
\left(\frac{\pi}{2} - \psi^e_{t}(x)\right). \]

 \qed

Note for any point $x \in \frak{N} $ that 
 $k^t(x) < 0$ and so with this and observation \ref{ob1}  we  have 
$\sum_{e \in S} (\frac{\pi}{2} - \psi^e(x)) > \pi|S| $  
and in particular condition $(n_2)$ is  necessary. 

$(n_3)$ and $(n_4)$  rely  on certain a pair of related 
formulae.
\begin{formula}\label{loop}
Let $A^{i,i+1}$ be the angle slot between $e_i$ and 
$e_{i+1}$ in a snake $\{e_i\}_{i = k}^l$.  We have 
\[ \sum_{i = k}^l \theta^ {e_i}(x)
= |l-k| \pi - \sum_{i=k}^{l-1} A^{i,i+1}(x) -\psi^{e_k}_{t_{k-1}} - 
\psi^{e_l}_{t_l} \]
and if $\{e_i\}$ is a loop
\[ \sum_{i = k}^{l-1} \theta^{e_i}(x)
= |l-k|\pi - \sum_{i=k}^{l-1} A^{i,i+1}(x) . \]
\end{formula}
{\bf Proof: }
Simply note both  that 
\[\theta^{e_i}(x) = \left(\frac{\pi}{2} - \psi_{t_{i-1}}^{e_i}(x)\right) 
+ \left(\frac{\pi}{2} - \psi_{t_i}^{e_i}(x)\right), \]
and that 
\[ \left(\frac{\pi}{2} - \psi_{t_{i}}^{e_{i}}(x)\right) 
+ \left(\frac{\pi}{2} - \psi_{t_{i}}^{e_{i+1}}(x)\right) 
= \pi - A^{i,i+1} (x), \] 
and sum up.
\qed

The second formula immediately   
implies $(n_3)$ immediately, and from 
fact \ref{ob1} we have  $\psi^{e}_{t} > \frac{-\pi}{2}$
allowing  the first formula to demonstrate $(n_4)$.

\end{subsubsection}

\begin{subsubsection}{Surjectivity}\label{surjectivity}

In this section we will finish the proof of theorem \ref{main} by showing 
$\Psi$ maps $\frak{N}$ onto $N$. 
To do it let's assume the contrary allowing that that $\Psi(\frak{N})$ 
is strictly contained in $N$ and produce a contradiction.
With this assumption we have
a  point $p$ on the boundary of $\Psi(\frak{N})$  
inside $N$.
Note $p = \Psi(y)$ for some 
$y \in \partial \frak{N}$ and that 
$(C+y) \bigcap \frak{N}$ empty, since other wise for some 
$w \in C$ we would have $(y + w) \in \frak{N}$ hence forcing by the 
openness of $\Psi$ $\mbox{ } p = \Psi(y) =\Psi(y +w)$ to be 
in the interior of $\psi(N)$.

At this point we need to choose a particularly nice conformal version of 
$y$, which requires the notion of a 
stable boundary point of $\frak{N}$.  
Before defining stability note since  $\frak{N}$ is a convex 
set with hyper plane boundary if 
$x \in \partial \frak{N}$ such that  
$(x + C) \bigcap \frak{N} \neq \phi$, then 
$(x + C) \bigcap \partial \frak{N} $
is its self a convex $k$ dimensional set.    

\begin{definition} \label{stable}
A point in $x \in \partial \frak{N}$ is stable if 
$(x + C) \bigcap \frak{N} = \phi$  and $x$ is in the interior of 
$(x +C) \bigcap \partial \frak{N} $ as a $k$ dimensional set.
Any inequality forming  $\frak{N}$ violated in order to
make a stable $x$ a boundary point  will be 
called a violation.
\end{definition} 

The key property of a stable point is that 
a conformal change $w \in C$ 
has $x+\epsilon w \in \bar{\frak{N}}^c$  for all $\epsilon > 0$ 
or for some  sufficiently small  $\epsilon > 0$ we have 
$x + \epsilon w$ must still be on $\partial{\frak{N}}$ and 
experience exactly the same violations as $x$.
  The impossibility of any other
phenomena when conformally changing a stable point is at the heart of the 
arguments in lemma \ref{sub} and lemma \ref{sub2} below.
At this point subjectivity would follow if
for a stable  $x \in \partial {\frak{N}}$ 
we knew that  $\Psi(x)$ could not be in  $N$, contradicting the choice of 
$p = \Psi(x)$ as needed. 
 
We will prove this by splitting  
up the possibilities into the two cases in
lemma \ref{sub} and lemma  \ref{sub2}. 

Before starting lets define ...

\begin{definition}\label{head}
The end of a snake $\{e_i\}_{i =k}^l$ is said to have 
a head with respect to $x$ if 
$d^{t_l}(x) = \{0,0,\pi\}$ with the pair of zeros located at the angle 
slots of $t_l$  along $e_l$. 
\end{definition}

\begin{lemma}\label{sub}
If $x \in \partial {\frak{N}}$ is stable
and $\alpha^i(x) =0$ for $\alpha^i$ in a triangle where 
$k^t(x) <0$, then $\Psi(x)$ is not in $N$.
\end{lemma}
{\bf Proof:}
I will  suppose that  $\Psi(x) \in N$ and produce a violation to 
the $(n_3)$ or $(n_4)$ conditions.

%e $(y+C) \bigcap \frak{N}$ is empty  yet  
%(y+C) \bigcap \bar{\frak{N}} \neq \phi$ we have 
%(y+C) \bigcap \partial \frak{N}$ is a convex subset.  
%hose a conformally equivalent point $x$ in the 
%interior of this set in its subspace topology.
%now with such a choice that either
 %of the boundary conditions forming 
%$\frak{N}$ without sending $x$ into $\bar{\frak{N}}^c$.  
%Such an $x$ choice is particularly stable.
%Now we'll use the assumption that 
%there were no triangles with $k^t(x) = 0$, and hence  
%$x$ violates and only violates  the $\alpha^i(x) = 0$ condition 
%for various $\alpha^i$. 

Now look at an  angle slot which is zero in  triangle $t_0$ satisfying  
$k^{t_0}(x) < 0$.
View this angle as living  
between the edges $e_0$ and $e_1$.  
Note that in order for $x$ to be stable that  
the $\epsilon w_{e_1}$ transformation 
(with its positive side in $t_0$) 
 must be protected 
by  a zero on the  $-\epsilon$   side forcing the condition that 
$\epsilon w_{e_1} \in \bar{\frak{N}}^c$ , or else for small enough 
$\epsilon$ we have 
$x + \epsilon w_{e_1}$  would be a conformally 
equivalent point on $\partial{\frak{N}}$  with
fewer violations. Call this neighboring triangle $t_1$.
If we see a pair of zeros and a $\pi$ facing the 
$t_0$ from $t_1$ we stop.   Otherwise let $e_2$ be another 
edge bounding the zero angle slot in $t_1$ and repeat the above procedure if
$k^{t_1} < 0$.  If $k^{t_1}(x) = 0$ form 
\[x +  \epsilon w_{e_1} + \epsilon w_{e_2}.\]
Note that the only way this construction could have difficulty is 
precisely the case in which $t_1$ was a head - in which 
case we already stopped.

So we may continue this process  
forming a snake $\{e_i\}$ with $A^{i,i+1} =0$ until we hit a head.  
Notice we can 
also make the same construction in the other direction.

%RIDDDDDD!!!!! look at  
% $\epsilon v_{e_{2}}$ 
%%either or both of the  $\epsilon v_{e_2}$ and $\epsilon v_{e_2}$ so that  
%If we see a single zero note 
%\[x +  \epsilon v_{e_1} + \epsilon v_{e_2}\]
%needs once again a zero on the $-\epsilon$ sides of 
%$e_2$.  If this neighbor has zero curvature there still no issue since we will be adding and subtracting the same amount of angle.
%If there are two zeros on the $- \epsilon$ side of $e_1$ and its not a head
%then either of the new zeros must its self have a protection on the $-\epsilon$ side of $e_3$ and $e_3$.  Choose one of these to continue the zero chain.
%(there is one case in which there is choice of $e_2$, namely when  
%$e_1$ faces a trianlg with $d^t(x) = \{0,0,0\}$ in which case 
%at least one of the possible choices of $e_2$ 
%must have  this property).  

Note as such if  his snake formed a embedded 
loop, a barbell, an embedded snake with two heads or 
a balloon with a head 
then the $A^{i,i+1} = 0$ condition  would 
contradict one of the formulas in formula \ref{loop}, 
hence violating $(n_4)$ or $(n_5)$, and we would be done.  
 It will be shown that one of these cases   must occur.

%If this snake formed a loop  by formula \ref{loop}
%\[ \sum_{\{e_i\}} \theta^{e_i}(x) = card\{e_i\} \pi, \]
%violating condition $(n_3)$ and producing our needed contradiction.

To produce the needed snakes note that by finiteness in the positive direction 
there is a first time when some $t_k = t_l$ and $k < l$ or we 
terminated at a head before such an over lap.
If this sequence  
terminated in a head look at the snake in the negative direction 
and if it terminates in a head then we are done.   If not we have the same 
situation as the positive snake not terminating in a head, i.e. there is 
a first  
$e_k$ when the one headed snake hits itself.  
If we hit the head we have our needed embedded loop.  If not  
we have two possibilities  either that $e_{k+1}$ can be chosen to be $e_{l+1}$
in which case we have our need embedded loop or $e_{k+1}$ must be
$e_{l}$.  In this case we can reverse the construction 
going form $e_0$ to $e_l$to form the needed  barbell.

 Suppose in the positive direction we experienced our first moment when
for $k <l$ $t_k = t_l$.  Then as above if it hinges so   
$e_{k+1}$ can be chosen to be $e_{l+1}$ we have our embedded loop.  
If not we will form a chain in the opposite direction starting with 
 $e_{l} = \hat{e}_0$.  Now we are searching for the first $ m < n$ when 
$\hat{e}_n$ hits $\{\hat{e}_{j}\} $ at $\hat{t}_m$  
or hits some  or $\{e_i\}_{i=k}^{l}$
at $t_p$ or terminates in a head.  
If it terminates in a head we can form our needed  balloon as   
\[\{\hat{e}_i\}_{i=n}^{0} \bigcup \{ {e}_i\}_{i = k}^{l} \bigcup
\{\hat{e}_i\}_{i=0}^{n}. \]
If our snake hits itself at $\hat{t}_m$  and $\hat{e}_{n+1}$ can
be chosen as $\hat{e}_{m+1}$ we once again get our needed embedded  loop.  
If it hinges such that
 $\hat{e}_{n+1}$ must $e_{m}$
then we can double back to form our needed barbell
\[\{\hat{e}_i\}_{i=0}^{e_n} \bigcup \{\hat{e}_i\}_{i = m-1}^{0} \bigcup
\{e_i\}_{i=l}^{k}. \]

The other possibility is that $\{\hat{e_i}\}$ hits first $\{e_i\}$ 
at $t_p$ with $\hat{e}_n$. 
%So there is a zero  angle between $e_p$ or $e_{p-1}$ and $\hat{e}_{m}$. 
As always there are two possibilities for how  they hinge and in either
case we can form one of the following embedded loops  
\[\{\hat{e}_i\}_{i=0}^{n}  \bigcup
\{e_i\}_{i=p}^{l} \]
or
\[\{\hat{e}_i\}_{i=0}^{n} \bigcup 
\{e_i\}_{i=p-1}^{k}. \] 

So in any case violation to $(n_4)$ or $(n_5)$ can be produced.  

\qed

%It is worth noting as a schollium that the condition $(n_3)$ really 
%only needs to be verified on embedded loops and ``bar bells'' lops - 
%a finite set of objects in a triangulation; and similarly 
%$(n_4)$ only needs verification on embed sakes and half barbells. 
%Now suppose we are the compliment of this case and...

\begin{lemma}\label{sub2}
If a stable $x$ satisfies the condition that 
if $\alpha^i(x) =0$ then $\alpha_i$ is in a 
triangle $t$ 
with
$k^t(x) = 0$, 
then $\Psi(x)$ is not in $N$. 
\end{lemma}
{\bf Proof:}
In order to be a boundary point of $\frak{N}$ 
for some $t$ we have that 
$k^t=0$. We will be looking at the set of 
all triangles with $k^t = 0$, $Z$,  
which is not all of $M$ and has a non-empty boundary 
To see this note 
\[\sum_{t \in \frak{P}} k^t(x) = \sum_{e \in v} A^i -\pi F 
 = 2 \pi V - 3 \pi F + 2 \pi F \]
\[  =  
2 \pi V - 2 \pi E + 2 \pi F  = 2 \pi \chi(M) < 0,\]
so there is negative curvature somewhere.

%(In the manifold with boundary case one similarly gets that the region of zero curvature is forced to have boundary in the interior of the triangulation.  This because the assumption that  $\sum_v \alpha^i(x) = \pi$ at a boundary vertex preserves the above fact that   $\sum_{t \in \frak{P}} k^t(x) = 2 \pi \chi(M)$.)

By the stability  
of $x$ once again 
there can be no conformal transformation capable of moving negative 
curvature into this set. Suppose we are at a boundary 
edge of $Z$, call the triangle  on the $Z$ side of 
the boundary edge $t_0$ and the triangle on the
non-boundary edge  $t_{-1}$.      Since  $t_{-1}$ has negative 
curvature (and hence   no $\pi$ angles) the obstruction to 
the $\epsilon w_{e_0}$ transformation
being able to move 
curvature out of $Z$ must be due to  $t_0$.  
In order for $t_0$ to protect against 
this there must be zero along $e_0$ in the $t_0$ side.

Now we will continue the attempt to suck curvature out with a curvature vacuum.
Such a vacuum is an element of $C$ indexed by a snake.  The first edge in 
the snake is the boundary edge $e_0$.  If in $t_0$ $e_0$ faces a 
$\pi$ we stop.  We say that we stopped at a head.
If not let $e_1$ be the other edge sharing the unique zero angle along $e_0$ in $t_0$ and if 
$e_1$ is another boundary edge we stop. If $e_1$ is not a boundary edge
 use $w_{e_1}$ 
to continue the effort to remove curvature.  Continuing this 
process $n$ steps forms  a snake $\{e_i\}_{i =0}^{m \leq n}$ and  
$x + \epsilon \sum_{i =0}^{m \leq n} w_{e_i} \in C$.

Suppose a   vacuum hits its self and  $t_n = t_m$. 
Then we must have an extra zero in 
$t_m$ in which case we have a $d^{t_m}(x )= \{0,0,\pi\}$.  Note it 
 is not a head 
with respect to either direction and it fact now form a vacuum 
loop.  In the conformal change associated to this
loop  consistently changes the angles with value  $\pi$ to having value 
$\pi - 2 \epsilon$, 
a contradiction to stability.  

So any vacuum in fact pokes  through $Z$ into $Z^c$. In fact this 
argument shows us 
something slightly stronger, namely  if a vacuum hits a triangle with 
$d^{t_m}(x )= \{0,0,\pi\}$ then it is a head.  If not after we poke through
we could still reduce the $\pi$ to $\pi -\epsilon$ conformally.  
This because under the lemma's hypothesis,\
there can be no zeros in $Z^c$ protecting the 
vacuum form consistently sucking.

From this note one vacuum can never pass through another since this would
force a 
$d^t(x) = \{0,0,\pi\}$ triangle which is not a head for 
at least one of the vacuums.
Similarly the outside edges of a vacuum always face zero angles in the vacuums
and to be in $N$ two zeros can never face each other since then 
\[\theta_e = \frac{\pi - \pi}{2} +  \frac{\pi - \pi}{2} = 0.\] 
So {\bf all} the edges associated to 
distinct vacuums are distinct.

Now simply let 
$S$ be the removal
from $Z$ of all these vacuums. If $S$  is non-empty 
then every boundary edge of this set faces a zero in a 
triangle of zero curvature so formula \ref{set} receives all zeros from 
the boundary terms.  
Similarly each triangle having exactly zero curvature gives us exactly 
a $\pi$ for each internal triangle in  formula \ref{set}, so 
$\sum_{e \in S} \theta^e(x) = \pi |S| $ as needed to violate $(n_2)$. 
 
So we are reduced to seeing that $S$ is nonempty.  Since two vacuums can never 
border each other, this is reduced to seeing that every vacuum has a $Z$ 
internal 
edge.  Well suppose not then our vacuum would be an embedded snake with all 
boundary edges having a zero along them  and all internal angles being 
zero and zero curvature.   This forces our vacuum to 
have only triangles $t$ with 
$d^t(x) = \{\pi,0,0\}$,  so could only be a pair of heads.  
When two heads face each other at edge $e$ we have $\theta^e = 2 \pi$ 
contradicting that
fact we are in $N$.    
So $S$ must be  non-empty and we are done.

\qed

\end{subsubsection}

\end{subsection}

\end{section}

\end{chapter}

\begin{chapter}{From the Discrete to the Continuous}\label{dis2con}

This chapter is dedicated to setting up the geometry 
and probability needed to
compute the random variables discussed in section 1.2.2 of the 
introduction.  In particular 
we fill in the details to all the steps in the probabilistic 
proof of the Gauss-Bonnet theorem sketch in section 1.1.  In section 
3.1 we develop all the geometric tools necessary to prove theorem
\ref{dtri} from the introduction.  This includes section 3.3 
where we examine some properties of Delaunay triangulations on surface
independent of the rest of the thesis but of interest to anyone 
wanting get a feel for these triangulations.

In section 3.2 we explore random Delaunay triangulation.  In 3.2.1  we 
develop all the ideas need to prove the theorem \ref{size} and
Euler-Delaunay-Poisson formula.  In sections \ref{coords} and 
 \ref{dercomp} 
we develop a formula for 
computing random variables on the space of random Delaunay triangulations 
(or complexes), and in particular prove the  Euler-Gauss-Bonnet-Delaunay
formula form the introduction.  Section 3.2.4 contains some 
particularly boring 
facts concerning the measurability of certain function and sets
which arise in the first two sections of this chapter. 

\begin{section}{Delaunay Triangulations}\label{triangus}

%This section is dedicated to proving the geometric facts  leading 
%up to  the Euler-Delaunay-Poisson formula  from the introduction.  
This section is dedicated to the exploration of Delaunay triangulations.
The technical backbone for all that occurs in this proof is theorem 
\ref{lit} from the introduction, which is dealt with in section 
3.1.2. 
The geometry continues in section  3.1.3 where  theorem \ref{dtri} is 
proved (with the help of certain "inflating families" also 
dealt with in section 3.1.2).  

Then in section \ref{prop} we prove some interesting 
properties about Delaunay triangulations.  These facts will 
not  needed in the rest of thesis 
but are of  interest in showing how certain facts about 
Euclidean Delaunay triangulations carry over to surfaces.  
 The facts explored  include that Delaunay triangulations 
are local,  several local facts, justification of algorithm constructions, 
as well as results showing  that the Delaunay triangulation 
of a dense ${\bf{p}}$ 
is in several ways optimal amongst 
dense triangulations; where a triangulation $T$ is called dense is  
if each  $t \in T$  has its 
vertices on a ball of radius less than 
$\delta$.  For the Euclidean versions of essentially 
all these facts with quite different proofs see
\cite{Fo}. Throughout this section I'll assume $\delta$ is simply 
$\min\{\frac{i}{8},\tau\}$.
To entice the reader perhaps I'll mention now what these 
optimality properties are.  For our  first optimality 
property we have. 

\begin{property}\label{delen}
Among dense triangulations
associated to a dense ${\bf p}$ 
the energy
\[E(T) = \sum_{\{e \in T \}} length(e) \]
is minimized precisely  
at the Delaunay triangulation. 
\end{property}

In section 4 we will see that the  ``gradient flow'' of this energy
tells us how to deform 
a dense triangulation associated to a dense ${\bf{p}}$ into its 
Delaunay triangulation. 

Notice the smaller this energy the squatter the triangles.  That Delaunay 
triangles minimize this energy is one reflection of  the fact that 
they prefer fat triangles. There are many realizations  of this fact, 
another is the  fact that the Delaunay triangulation attempts to 
minimize the sizes of associated spheres.  
To articulate this we first must acknowledge that 
by  lemma \ref{lit} of the next section the vertices of  
a triangle $t$ in a dense triangulation lie on 
a uniquely  associated ball which will be denoted  $B_t$.  

\begin{property} \label{maxrad}
Among all dense triangulations $T$ 
associated to a dense ${\bf{p}}$  the 
Delaunay triangulation minimizes 
\[ maxrad(T) = max\{radius(B_t) \mid t \in T \}.\]
\end{property}

In section \ref{prop} we will also see the 
sense in which this is locally true.

\begin{subsection}{Some Geometric Reminders and Notation}\label{georem}
%This proof uses a bit surface geometry and not much else. In this section I'll discuss  the needed surface geometry results and notation for future reference(see  \cite{Do} for the proofs). 

To understand this paper 
one  must be aware of geodesics, the exponential map 
$exp_p: T_pM \rightarrow M$, and this map's 
implicit interaction with balls and spheres.  Let the ball of radius $r$ at $p$ be $exp_p(B_r(0))$ - where $B_r(0)$ is the  open  ball of radius $r$ in $\Bbb E^2 \cong T_p M$, and denote it $B_r(p)$.  Let the sphere $S_r(p)$ be its boundary.  The first half of the needed results can be summed up in the following lemma.

\begin{lemma}[Geometric Reminders]
\label{normal}
 Assume $M$ is a compact Riemannian surface then: 
\begin{enumerate}
\item (Normal and Convex Neighborhoods)
$i$ and $\tau$ are greater than zero.
\item
For any $p$ we have that  $exp_p$ is diffeomorphism of $B_i(0)$ onto $B_i(p)$; and if $p$ and $q$ satisfy $d(p,q) < i$, then there is a unique geodesic of length less than $ i$ between them. 
\item (Gauss's Lemma)
The unique unit speed geodesics from $p$ to points in $B_{i}(p)$ are given by $\gamma(r) = exp_p \left(r \frac{v}{||v||}\right)$ for some $v$; and any such geodesic is  orthogonal to $S_r(p)$.
 \end{enumerate} 
\end{lemma}

It is worth explicitly reminding the reader that given an orthonormal basis $\{ e_1, e_2 \}$  at  $p \in M$ we have the lovely normal coordinates:
\[ N(p,z_1,z_2) = \exp_p(z_1 e_1 + z_2 e_2):B_i(0) \subset \Bbb R^2 \rightarrow M. \] 
Sometimes it is useful to think in terms of angular and radial coordinates.  Let $v(\theta) = \cos(\theta)e_1 + \sin(\theta)e_2$, then as alternate coordinates we have the  geodesic polar coordinates: 
 \[ G(p,r , \theta) = \exp_p(r v(\theta)): (0,i) \times S^1 \rightarrow M. \] 
 
In order to vary $p$ in the above it is necessary to have smoothly varying orthonormal frames and they will be denoted $f = \{e_1,e_2\}$.  We can always construct one on, say, a convex set; and can even globally have one on $M - \{points\}_f$.  This is accomplished by Graham-Schmidting a pair of generic vector fields, where  $\{points\}_f$ is the finite set of points where the vector fields are not linearly independent.  In the the presence of a frame we have a canonical choice for a  a $\frac{\pi}{2}$ rotation field   $\Theta$;  by using the fact $v(\theta)$ parameterizes the tangent spaces we can let  \[ \Theta(v(\theta)) = v(\theta)^{\perp}=-\sin(\theta) e_1 + \cos(\theta) e_2. \]

The last frame idea used is that of a {\bf geodesic frame} at $p$.  
Fixing an orthonormal basis  $\{ e_1,e_2 \}$ at $p$, let the geodesic frame be the the frame given by the parallel transport  of this orthonormal basis of $T_p M$ along the geodesics spitting out from $p$.  

The other bit of geometry used are some basic Jacobi field results.   Recall that a Jacobi field is a vector field along a geodesic $\gamma(r)$ satisfying 
\[ \frac{D^2 J}{{dr}^2}  = -R( \dot{\gamma}, J) \dot{\gamma} ,\]
with the initial condition $J(0) = V$ and $\frac{D J}{dr}(0) =W$. (The choice of $r$ here,  as opposed to the usual $t$,  stems from the fact that our Jacobi fields will be thought of as  living along geodesics parameterized by $r$ in some geodesic polar coordinates.)  

Jacobi fields have the wonderful property of being in 1-1 correspondence with smooth one parameter families of geodesics (in the standard notation $J = \frac{\partial \Gamma}{\partial s} $ where $\Gamma(r,s)$ is a geodesic for each $s \in [s_0 - \epsilon, s_0 + \epsilon ]$) and  $\Gamma(r,s_0) = \gamma(r)$ .
On a surface they come in four flavors.  To taste these flavors first one notes Jacobi fields with initial conditions perpendicular or parallel to $\frac{d\gamma}{dr}(0)$ remain as such for all time.   Also the equation is a second order O.D.E., hence linear in its initial conditions - so a Jacobi field can be decomposed into its component along $\gamma$ and its perpendicular component simply by decomposing its initial conditions as such.  

So we have all Jacobi field are linear combinations of the following types (and some examples of corresponding $\Gamma(s,r)$):
\begin{itemize}
\item
$J(0) = v(\theta_0)$ and  $\frac{DJ}{dr}(0) = 0$ (The two parameter family corresponding to this case is the friendly $\Gamma(s,r) = exp_p((s + r)v(\theta_0))$.)
\item
$J(0) = v^{\perp}(\theta_0)$ and  $\frac{DJ}{dr}(0) = 0$ 
\item
$J(0) = 0$ and  $\frac{DJ}{dr}(0) = v(\theta_0)$ 
\item
$J(0) = 0$ and  $\frac{DJ}{dt}(0) = v^{\perp}(\theta_0)$   (The two parameter family corresponding to this case is the friendly $\Gamma(s,r) =exp_p(r v(s))$.  In particular this Jacobi field is precisely  $ {N}_{*} \left( r v^{\perp}(\theta_0) \right) =G_{*} \left( \frac{\partial }{\partial \theta} \right) $.
\end{itemize}

Sometimes one starts with a $\Gamma(s,r) $ and wants to understand the associated field - a well known example, that will prove relevant to us, can be constructed by fixing  a geodesic $\alpha(s)$ and a frame giving the vector field $v(\theta)(s)$ along  $\alpha(s)$. The Jacobi field along  $exp_{\alpha(s_0)}{rv(\theta)}$  corresponding  to 
$exp_{\alpha(s)}{\left( r v(\theta) \right)}$ is the  the  one with initial conditions $J(0) =  \frac{d \alpha}{ds}(0)$, $\frac{D J}{dr}(0) = \frac{D v(\theta)}{ds}$.

Here are some facts we will be needing about Jacobi fields, and in particular what they look like in normal coordinates.  

 \begin{lemma}
\label{jac}
\begin{enumerate}
\item
The Jacobi field along $ exp_{p}({rv(\theta_0)})$ with initial conditions $J(0) = v(\theta_0)$ and  $\frac{DJ}{dr}(0) = 0$ in normal coordinates is $J(r) = v(\theta_0)$.
\item
The Jacobi field along  $exp_{p}({r v(\theta_0)})$ corresponding to $J(0) = 0$ and 
$\frac{d J}{dr} = v^{\perp}(\theta_0)$ in normal coordinates is  $r v^{\perp}(\theta_0)$.  Calling $ ||r v^{\perp}(\theta_0)||_M = j_{\theta_0}(r)$,  we have $j_{\theta_0}$'s  Taylor expansion is   $j_{\theta_0}(r) =  r(1  - r^2 \frac{k}{6} + O(r^3)) $.

\item
The Jacobi field along $ exp_{p}({r v(\theta_0)})$ with initial conditions $J(0) = v^{\perp}(\theta_0)$ and  $\frac{DJ}{dr}(0) = 0$ in normal coordinates can be written  $J(r) = h_{\theta_0}(r) \frac{r}{j_{\theta_0}(r)}v^{\perp}(\theta_0)$. 
\item
Here are a few immediate consequences of part two - the area of a ball at $p$ function, $a(r)$, satisfies $a(r) = r^2(\pi - \frac{\pi k}{12} r^2 + O(r^3))$; and the product   $j_{\theta_1} j_{\theta_2}j_{\theta_3} =  r^3(\pi - \frac{ k}{2} r^2 + o(r^3))$.
\end{enumerate}
\end{lemma}

It is necessary to have certain global estimates of the above $o(r^3)$ functions, resulting from the fact we are on a compact surface. 
%(This result is not in  \cite{Do},b ut is very basic.  A proof can be found in\cite{Le}.) 

\begin{lemma}\label{dippy}
\label{dual}
If $M$ is compact, then
there is a $C_M >0$ such that when $r < i $ we have the  $O(r^3)$ functions in 3 and 4 above all satisfying 
$|O(r^3)| < C_M r^3$ globally.
\end{lemma}

{\bf Proof of lemma \ref{dippy}:}
First I'll find the constant related to 4(c) above. The Jacobi fields are the solutions to an o.d.e. - so vary continuously with initial data - which is indexed by $ (p,v) \in  UTM$ (the unit tangent bundle).  In particular $j_{\theta_0}(r) = j(p,v,r)$ varies continuously with initial data.  Now recalling from above the Taylor expansion in the radial variable at $r = 0$ we have   $j(p,v,r) = r + \frac{r^3 k(p) }{6} + r O(r^3)(p,v,r)$; and form Taylor's formula the third term in this sum can in fact be represented as (assuming the metric is smooth) 
\[ r_4(p,v,r) = \frac{r^4}{3!} \int_{0}^{1}(1-t)^{3} \frac{d j}{d r^4} (v,p,tr)  dt .\]
In particular this term  is continuous even after dividing  by $r^4$.  Now 
$UTM \times [-i,i]$ is compact so we can feel free to take $| \frac{r_4(r,v,p)}{r^4} |$'s  maximum over this set for our $C_M$. 

Now observe that the other $O(r^3)$ are directly  related to this one and the curvature function, and since the manifold is compact $sup|k|$ exists. Choose the $C_M$ in the lemma to be the biggest of the constructed  constants among these $O(r^3)$ functions.

\qed

\end{subsection}

Now we will prove some potentially less familiar geometric facts.
The first of which will 
be lemma  \ref{lit} from the introduction.  My original proof 
of this fact was quite inelegant and can be found the the appendix.  
I'd like to thank Albert Nijenhuis for sharing his beautiful proof with me.

%\begin{subsection}{Small circles and Inflating families}

%\begin{lemma}\label{lit}If $\{u,z,w\}$ lie on the boundary of a disk of radius less than $\frac{i}{8}$  disk, then this disk is the unique such disk among disks of radius less than $\frac{i}{8}$.  \end{lemma}
 
The next facts  we will need  concerns  the  notion of an 
inflating family of circles through a pair of points $\{ p, q \}$.  
As a set, this family will be all circles of radius 
$r < \delta$, 
passing through both the points $p$ and $q$. The following 
lemmas justify 
the fact that this set 
can be thought of as the continuous inflating family of  circular 
balloons to the left or right "sides" of the geodesic through $p$ and $q$,  
as in figure $\ref{pen}$. 

Too articulate this given a 
continuous curve $c_{pq}(t)$ let $D_{pq}(t)$ be the continuous family of 
closed disks centered at $c_{pq}(t)$ of radius $d(p,q) + |t|$, and let 
$\partial D_{pq}(t)$ be the corresponding family of circles.  For 
starters we have a lemma gaurenteeing the existence of inflating families.

\begin{figure*}
\vspace{.01in}
\hspace*{\fill}
\epsfysize = 2in 
\epsfbox{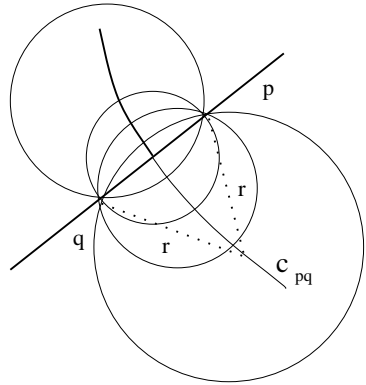}
\hspace*{\fill}
\vspace{.01in}
\caption{\label{pen} The Inflating Family}
\end{figure*}

\begin{lemma}\label{famex}
For each pair of points $p$ and $q$ such that $d(p,q) < \delta$
there is a curve $c_{pq} : (d(p,q) - \delta, \delta -d(p,q)) \rightarrow M$ 
such that every circles of radius less $\delta$ going through $p$ and $q$ is 
$\partial D_{pq}(t)$ for some $t$.    
\end{lemma}

Notice that the radius increases monotonically as $|t|$ does.  We need a lemma 
giving us another 
sense of monotonicity.  To articulate it 
we first develop a little notation.
If $d(p,q) < i$ let 
$pq$ be the unique shortest length geodesic connecting $p$ and $q$, let 
$mid(pq)$ be its midpoint, let $B_r(p)$ be the ball of radius $r$ at $p$,
and let $\partial B_r(p) $ be its boundary. To articulate the next lemma 
note that the geodesic 
connecting $p$ and  $q$ removed form $B_{i}(mid(pq))$ decomposes 
$B_i(mid(pq))$
into two open sets, which will be referred as the decomposition determined by 
$pq$.  Furthermore since $\delta < \tau$
any ball of 
radius less than $\delta$  with $p$ and $q$ on its boundary 
is also divided into two such pieces.

\begin{lemma}\label{mono}
$c_{pq}(-(\delta - d(p.q)),0)$ is in one half of the decomposition 
determine by $pq$ (call it  $H^-$) and $c_{pq}(0,\delta -d(pq)$  
in the other ($H^+$).  If $d(p,q) - \delta < c < d < \delta -d(q,q)$ 
then $D_{pq}(c) \bigcap H^+ \subset D_{pq}(d) \bigcap H^+$ and
$D_{pq}(d) \bigcap H^- \subset D_{pq}(c) \bigcap H^-$ with all 
the subsets proper.
\end{lemma}

The essence of this lemma is that figure \ref{poop} is accurate.  

\begin{figure*}
\vspace{.01in}
\hspace*{\fill}
\epsfysize = 2in 
\epsfbox{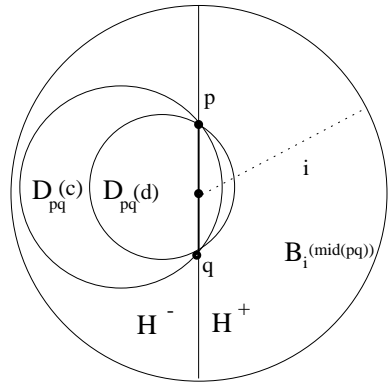}
\hspace*{\fill}
\vspace{.01in}
\caption{\label{poop} The Monotonicity of Area}
\end{figure*}

\begin{subsection}{Small Circle Intersection}\label{samllcir}

My original 
proof of lemma \ref{lit}  was a bit long winded 
 and 
I would like to thank Albert Nijenhuis for showing me the 
elegant proof presented here.  Actually the proof presented here 
is a much less elegant 
modification of the one Nijenhuis showed me
where the needed injectivity radius is 
chased through the argument.  Any errors or any realization  
that the 
injectivity bound is not as sharp as possible is solely my fault.
Note that the lemma will follow if one could show that 
two circles of radius less than $\frac{i}{8}$ intersected 
in at most two points, which is proposition \ref{nij5} below.
Another proof of this can be found in  
\cite{Le}  where the sharper bound of  $\frac{i}{6}$ is also demonstrated.

If $a$ and $b$ have distance between them less  
than the injectivity radius, 
 $d(a,b) < i$, denote as $ab$ the unique minimal length geodesic connecting 
them.
Let $x$ and $y$ be the centers of two intersecting 
circles of radius less than 
$\frac{i}{8}$. Let $C_y$ be the circle centered 
at $y$, and  $xy_i$ be the open geodesic segment 
containing $xy$ with midpoint  $x$ of 
length $i$.

\begin{lemma} \label{nij4} 
 $xy_i$ intersects $C_y$ in exactly two points.
\end{lemma}
{\bf Proof:}
Note $d(x,y) < \frac{i}{4}$ and any point on $C_y$ has a distance less than
$\frac{3i}{8} < \frac{i}{2}$ form $x$, so the diameter of $C_y$ 
intersecting $xy_i$ is in fact included in it.  
So $xy_i$ intersects $C_y$ at least twice.
However  
a geodesic ray form  $y$ of length less than $\frac{3i}{4} < i$ can hit each 
circle of radius $r < i$ centered at  $y$ only once, so $xy_i$ 
can only hit $C_y$ twice.  
\qed

Consider  
$f:C_y \rightarrow \Bbb R$, $f(z) = d(x,z)$, so $f$  measures the distance 
between $x$ and the points of $C_y$.

\begin{lemma} \label{nij2} 
If $z,w \in C_y$, $z \neq w$, and $f(z) = f(w)$, then $f$ has (at least)
one critical point on each of the two circular arcs zw.   
\end{lemma}   

{\bf Proof:}
This follows  from a standard  min-max argument.
\qed

%\begin{lemma} \label{nij3} 
%If $p$ is a critical point of $f$  then $d(x,p) < \frac{3i}{4}$
%and $xp$ either contains $xy$ or an be continued in 
%the $z$ lies on a geodesic segment in
%$V$ that contains $x$ and $y$.  
%\end{lemma}  

\begin{lemma} \label{nij3} 
If $p$ is a critical point of $f$, then $p \in xy_i$.  
\end{lemma}  

{\bf proof:}
Since $p$ is a critical point and has a distance at most 
$\frac{3i}{8}$ from $x$, the tangent to the segment $xz$ 
at $z$  is perpendicular to $C_y$.  So by Gauss's lemma 
the distance to $y$ along 
this geodesic or its potential $\frac{i}{8}$ or less 
length extension  
is less than $\frac{i}{2}$ form $x$, and we have  
$p \subset xy_i$. 
\qed

\begin{proposition} \label{nij5} 
Any two distinct circles of radius less than $\frac{i}{8}$ 
have at most 2 points in common.
\end{proposition}

{\bf proof:}
Suppose not.  
By lemma \ref{nij2}  there would be at least 3 distinct critical points, 
which would all lie on $xy_i$ by lemma \ref{nij3}.  But lemma \ref{nij4}  
assures there are only two such critical points, the need contradiction.
\qed

\end{subsection}

\begin{subsection}{Inflating Families}

 I will construct the  family by first describing the point set  
along which the centers of the circles in the family 
live with a different parameterization than that of lemma \ref{famex}.  
Note  a point set this curve consists of points satisfying 
$d(p,z) - d(q,z) = 0$ - with $d(p,z)$  with less than $\delta$.
We will always denote as  $\bar{pq}$ the most sensible
connected  extension of $pq$.  For the following lemma 
let it be the connected 
extension of $pq$ 
in  $B_{\frac{i}{4}}(p) \bigcup  B_{\frac{i}{4}}(q)$.  
With this notation we have:

\begin{lemma}\label{curv} 
If $d(p,q) < \frac{i}{6}$ the  point set described by $d(p,z) - d(q,z) = 0$ 
in $B_{\frac{i}{4}}(p) \bigcup  B_{\frac{i}{4}}(q)$ can be described
by a curve 
 $c_{p,q}(t)$ with $t \in (c,d)$ satisfying  
\begin{enumerate}
\item
 $c < 0 < d$, $c_{p,q}(0) = mid(pq)$,  and  $mid(pq)$ is
 the unique point of $\bar{pq}$ on $c_{p,q}(t)$
\item
$d(c_{p,q}(t),p)$ and $d(c_{p,q}(t),q) $ strictly increase as the 
parameter  $|t|$ increases.
\end{enumerate}
 
\end{lemma} 

{\bf Proof:}
 To see that the point set  is a nicely parameterized curve it is 
useful to note that  it can be described as 
the integral curve of a vector field.  
Let $D_p$ denote $ d(p,x)$ and let $\nabla D_p$ denote its gradient.
Note the solution to the equation $D_p - D_q=0$
are integral curves of  the vector field $ \Theta (\nabla D_p  - \nabla D_q) $
where $\Theta$ is a  $\frac{\pi}{2}$ rotation field.
 
To understand these integral curves we will first deal with the 
the uniqueness of $mid(pq)$: suppose a point $l \neq mid(pq)$ is 
on $pq$.  Then $l$ is within $\frac{i}{4}$ of $p$; hence 
$d(p,l)$ is determined by the length of the segment of  $pq$ 
from $p$ to $l$, similarly for $q$ 
(using $\frac{i}{2}$).  Now note that as we move from $mid(pq)$ toward, 
say, $p$ that $D_p$ decreases while $D_q$ increase - 
so $D_p - D_q \neq 0$ at another point of $pq$.  
When $l$ is on $ \bar{pq}  / pq$, say above $p$, the segment 
of $\bar{pq}$ from $q$ to $l$ in fact covers the shorter segment 
from $p$ to $l$ - forcing  $D_p - D_q \neq 0$ once again.  
So $l$ cannot satisfy $D_p - D_q =0$, forcing $mid(pq)$ 
to indeed be the unique point of $\bar{pq}$ on $D_p - D_q = 0$.

Now we we will see that we indeed get a union of curves by noting that 
 that the vector field has no zeros in this set.  
In fact the triangle inequality  tells 
$B_{\frac{i}{4}}(p) \subset B_{\frac{i}{2}}(p) \bigcap B_{\frac{i}{2}}(q)$ 
and in this  region we will prove the stronger fact that  
$\nabla D_p \neq c \nabla D_q $ for any $c$. First note that 
$\nabla D_p$ is unit length, with integral curves the 
geodesics emanating from $p$. So at $mid(pq)$ we have 
$ \Theta (\nabla D_p  - \nabla D_q) $ is length $2$.   
To finish the assertion assume at some point $p$ not on $\gamma_{p,q}$  
that  we have  $\nabla D_p  =  c\nabla D_q$.  First note from 
the that fact that $\nabla D_p$ is unit length we have $c = \pm 1$. 
There are two cases, 
first we'll deal with $c =1$.   Since the geodesics satisfy a second order 
O.D.E they are uniquely determined by their position and tangent vector, 
so when $c = 1$ we have  both the geodesic from $p$ and the geodesic form 
$q$ are the same curves.  Without loss of generality $p$ is further away than 
$q$ and this point lies along the same minimal length geodesic 
(of length less than $\frac{i}{2}$) which connects $p$ and $q$, i.e. 
$pq$.  In the case  $c = -1$ we can follow the geodesic form 
$p$ to the point and then from the point back to $q$ forming a geodesic 
of length less than $i$ - which then must by the 
definition of the injectivity radius 
be the unique such one, i.e. $pq$.

 To finish off the first part we need that our curve has only one component.  
This is intimately related to the second part. To see why we first 
look at the component of $c_{p,q}(t)$ in $B_{\frac{i}{4}}(p)$ and 
note any component of  $D_p - D_q = 0$  would have to 
have a point closest to $p$.   
 This closest point  is tangent to a sphere emanating form $p$.  
The same sort of phenomena must take place for the distance 
function to have a critical point; namely 
  if a point $z$ along any integral curve  of 
$\Theta (\nabla D_p  - \nabla D_q)$ is  a critical  point of 
the distance function $D(p, \cdot)$ then either 
$\nabla D_p  =  \nabla D_q$ or a circle is tangent to the solution curve. 
In the tangent case  $\Theta(\nabla D_p  - \nabla D_q) = c \Theta \nabla D_p$,
 or rather $\nabla D_p  - \nabla D_q = c  \nabla D_p$;  so both these 
situation have forced the case  $\nabla D_p  =c\nabla D_q$.  so we may use 
the above observation 
to note that the point where this occurs is on  
$\bar{pq}$; but from above to be on $D_p - D_q = 0$ and  
$\bar{pq}$ means you must be exactly $mid(pq)$.  
So we have both that every  component of   $D_p - D_q = 0$ 
in $B_{\frac{i}{4}}(p)$ contains $mid(pq)$ , and that the distance to $p$ 
parameterized by $t$ can have no critical points except at $mid(pq)$ 
(similarly for $q$).  

\qed

We can now prove lemma \ref{famex} by noting by part (b) of the above lemma we may reparameterize as claimed and that any circle of radius less 
than $\delta$ has its center  contained in 
 $B_{\frac{i}{4}}(p) \bigcup  B_{\frac{i}{4}}(q)$ so this new parameterization can be chosen on (and beyond) $ (d(p,q) - \delta, \delta -d(q,q))$.

Now we shall prove lemma \ref{mono}.

{\bf  Proof of lemma \ref{mono}:}
Notice the first part follows from 
part (a) of the above  lemma. To prove the rest of it 
it is useful to isolate a sub-lemma.

\begin{sub-lemma}\label{sir}
The containments between the halves must switch.
Precisely if 
$D_{pq}(c) \bigcap H^+ \subset D_{pq}(d) \bigcap H^+$ then
$D_{pq}(d) \bigcap H^- \subset D_{pq}(c) \bigcap H^-$ and visa versa. 
\end{sub-lemma}

{\bf Proof:}
Suppose that $D_{pq}(c) \bigcap H^+ \subset D_{pq}(d) \bigcap H^+$.
Note that either $D_{pq}(d) \bigcap H^- \subset D_{pq}(c) \bigcap H^-$
or $D_{pq}(c) \bigcap H^- \subset D_{pq}(d) \bigcap H^-$ since a 
violation of this inclusion would result in a third intersection of 
two circles of radius less than  $\delta$ - 
contradicting lemma  $\ref{lit}$.
From this observation, to violate the 
above choices would mean that 
$D_{pq}(c) \bigcap H^- \subset D_{pq}(d) \bigcap H^-$.  
Such an inclusion would force the circles to be tangent 
at there intersection points - and hence the centers of both disks to be on   
$pq$ (since the curve orthogonal to the tangent is a 
geodesic heading to the circle's center by Gauss's lemma).  
But then $d(p,c_1) = d(p,c_2)$  forcing both the centers  and the 
radii to be the same.  So the disks would be identical 
contradicting distinctness. 

\qed

We may finish the proof of  
the monotonicity of inflation half of lemma \ref{mono}. 
By the above sub-lemma we are left to explore  three cases.

The first case is where $c$ or $d$ is zero.  If $c = 0$ then the 
geodesic from $mid(pq)$ to $c_{p,q}(d)$ to the boundary of  $S_{pq}(d)$  
is strictly larger than $d(mid(pq),p)$  (by part (b) of 
lemma \ref{curv}).  
So a switch of containment is impossible in this case. For the $d=0$ 
case note by the above sub-lemma one is contained in the other, and as 
just noted 
$D_{pq}( 0) \bigcap H^- \subset D_{pq}(c) \bigcap H^-$.  
So it must be the case  
$D_{pq}(c) \bigcap H^+  \subset D_{pq}(0) \bigcap H^+$ 
as needed.

If $0 < d$ then from the above we know 
$D_{pq}(0) \bigcap H^+ \subset  D_{pq}(d) \bigcap H^+ $.
From this the fact    $D_{pq}(d) \bigcap H^+ \subset  D_{pq}(c)\bigcap H^+$
would by the continuity of $D_{pq}(t)$ and the uniqueness of the 
circles in the pencil  (lemma  $\ref{curv}$) force there to be a 
$0 < f < c$ (or $c <  f \leq 0$)   such that  
$\partial D_{pq}(f) \bigcap  \partial D_{pq}(d)$ 
contains a pont in $H^+$.  
This produces extra intersections of distinct circles 
contradicting lemma $\ref{lit}$.

To handle $d < 0$ recall form the first case that we have both 
$D_{pq}(0) \bigcap H^- \subset  D_{pq}(c) \bigcap H^-$ and  
$D_{pq}(0) \bigcap H^- \subset  D_{pq}(d)  \bigcap H^-$.  
So by the sub-lemma  we must have   
$D_{pq}(c) \bigcap H^+ \subset  D_{pq}(0)  \bigcap H^+$ and  
$D_{pq}(d)  \bigcap H^+ \subset  D_{pq}(0)  \bigcap H^+$.
 Now we can contradict  
$D_{pq}(d)  \bigcap H^+\subset  D_{pq}(c)  \bigcap H^+$ using same 
argument as in the second case, with the one variation 
 being that this time we construct an $f$ such that  $d \leq f \leq 0$ where 
 $\partial D_{pq}(f)  \bigcap H^+ \bigcap  \partial D_{pq}(c)  \bigcap H^+$  
contains other points.

\qed

\end{subsection}
%\end{subsuction}
\begin{subsection}{The Existence of Delaunay Triangulations}\label{deltri}

%{ \bf Proof of Theorem \ref{dtri}:}

Here we prove theorem \ref{dtri} from the introduction.  
To do so  first we needed to actually construct $R$.  
Part of $R$'s construction 
 is canonical;  namely the mapping of the 1-skeleton.  
This because  $|K_{\{p_1, \dots , p_n\}}|$ 
is in an affine space and we can let the edges of 
$|K_{\{p_1, \dots , p_n\}}|$ map onto there corresponding unique geodesic 
segments by factoring with an affine map
through the geodesic's unit speed parameterization. Now we need to extend 
this continuous mapping of the 1-skeleton to a continuous map of 
$|K_{\{p_1, \dots , p_n\}}|$.

Let ${\bf q} = \{q_i\}_{i=1}^{3}$  lie on a disk $D$ of radius less 
than $\delta$ and note the existence of $R$ requires  only an
identification of the correct triangle $q_1 q_2 q_3$. Throughout this 
proof a 
bold face letter will always denote a triple forming a face.

Since $D$'s radius is less than $\tau$
for each pair of distinct points
$a$ and $b$  on $\partial D$  we have that 
$D$ is split into its two halves.
If ${\bf q}$ lies in one half let $D^{ab}_{q}$ denote the closure of 
that half, the region enclosed by the bold dashed line in figure \ref{faces}.
Notice since the radius is less than $i$ 
that by the uniqueness of small geodesic this set is  convex.  

 \begin{figure*}
\vspace{.01in}
\hspace*{\fill}
\epsfysize = 2in 
\epsfbox{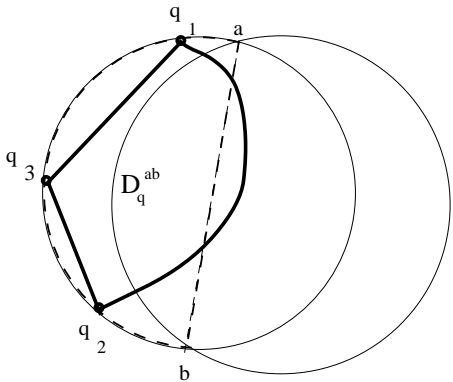}
\hspace*{\fill}
\vspace{.01in}
\caption{\label{faces} That Which Cannot Occur}
\end{figure*}

\begin{lemma}\label{contri}
If {$\bf q$} lies on a disk $D$ of radius less than 
 $min\{\tau, \frac{i}{6}\}$, 
then the triangle $q_1 q_2 q_3 = D^{q_1q_2}_q \bigcap D^{q_1q_3}_q 
\bigcap D^{q_2q_3}_q$ is a 
convex topological disk bounded by $q_1q_2,$ $q_1 q_3$ and $q_2 q_3$. 
Further more for any $a$ and $b$ on $\partial D$ if 
${\bf q} \subset  D^{ab}_{q}$ then $q_1 q_2 q_3 \subset D^{ab}_{q}$ . 
\end{lemma} 

{\bf Proof:}
Since $q_1 q_2 q_3$ is the intersection of  convex sets it is convex.
Each    $D^{q_iq_j}_q$ contains 
${\bf q}$ so by convexity  each contains all the $q_i q_j$.  
Now the only possible 
boundary of this set is the $q_i q_j$ or the boundary of $D$, but 
other than ${\bf q}$ all the points on $D$'s boundary  have 
been eliminated. So $q_1 q_2 q_3$ is a convex set with boundary $\{ q_i q_j\}$ 
contained in the convex disk $D$,  hence itself a topological 
disk as needed.

Now for the  second assertion. $q_1q_2q_3$
is an embedded disk 
so if the assertion fails  then 
one of the $q_iq_j$   crosses $ab$.    By continuity it will hit $ab$ in at
least 
a pair of points, as  in figure \ref{faces}.
So this pair of points 
is connected by two distinct geodesics of length less than  $i$.  
Hence we have a contradiction regarding the the uniqueness 
of geodesics of length less than the injectivity radius.

\qed

Now that we have the mapping $R$ let us prove   $R$ onto.  
Since $R(|K_{\{ p_1 \dots p_n \}}|)$ is closed 
(a union of closed sets) if it is not onto it must miss an 
open set, and the boundary of this open set must be composed of the
edges of triangles in $|K_{\{ p_1 \dots p_n \}}|$ 
(see bold region in figure $\ref{dink}$). So we would be done if the notion 
of edge and the notion of an edge belonging to a pair of faces 
becomes interchangeable in the presence of $\delta$-density.

 \begin{figure*}
\vspace{.01in}
\hspace*{\fill}
\epsfysize = 2in 
\epsfbox{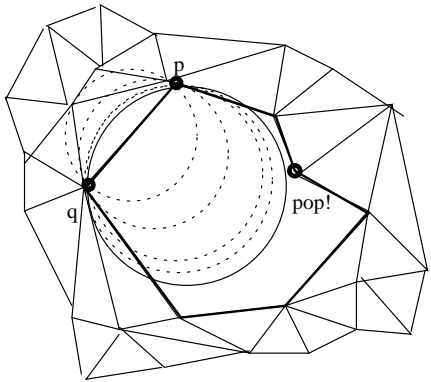}
\hspace*{\fill}
\vspace{.01in}
\caption{\label{dink} A Missing Region and Balloon Popping}
\end{figure*}

  To see this suppose we have an edge connecting a pair $p$ and $q$, 
then by its very definition there is a 
$k \in (d(p,q) - \delta,\delta -d(p,q))$ with the property 
$int D_{pq} D(k)$ is empty of points.  We may right off the bat use our 
assumption of $\delta$-density to see that this empty disk must correspond 
to a radius $d(p,c_{p,q}(k)) < \delta$.  Now start inflating to the left of 
$k$.  When $t = d(p,q) -  \delta $  the radius is $\delta$    
so $int D_{pq}(d(p,q) - \delta)$ contains 
a point.  By the monotonicity lemma (lemma \ref{mono}) 
once you hit a point moving leftward 
you cover for it all future time, so there is a  unique 
$c \in (d(p,q) - \delta,\delta -d(p,q))$ 
such that $\partial D_{pq}(c)$  
first contains a third point to the left of $k$. 
You may view this process as the blowing up of a circular balloon 
as dynamically represented in  figure \ref{dink}, 
and note we  proved that it must pop.  
Similarly there is a $d \geq k$ to the right  were the balloon pops.  
Since we are assuming that there are never four points on a 
circle  we have $c < d$.  So    $\partial D_{pq}(c)$ and 
$\partial D_{pq}(d)$ correspond 
to  the unique left and right faces  as needed.

%\begin{figure*}
%\vspace{.01in}
%\hspace*{\fill}
%\epsfysize = 2in 
%\epsfbox{face.eps}
%\hspace*{\fill}
%\vspace{.01in}
%\caption{\label{face} Balloon  Popping }
%\end{figure*}

Now let us deal with $R$ being 1-1.  
Since the individual faces are embedded this would follow immediately from the following lemma.

\begin{lemma}
Two distinct faces $p_1 p_2 p_3$ and $q_1 q_2 q_3$ can intersect 
only in a vertex or an edge.
\end{lemma}
{\bf Proof:}
Suppose $p_1 p_2 p_3$ and $q_1 q_2 q_3$ 
intersect and  have associated now intersecting circles 
$X$ and $Y$ respectively.  If $X$ and $Y$ intersect in one point then 
$\bf p$ and $\bf q$ being the only points of $p_1p_2p_3$ and $q_1q_2q_3$ on
$\partial X$ and $\partial Y$ respectively forces  
$p_1 p_2 p_3 \bigcap q_1 q_2 q_3 $ to be a vertex.
So we may assume by lemma \ref{lit}  that 
$X$ and $Y$ intersect in precisely 
two distinct  $a$ and $b$. 
Furthermore as in lemma \ref{mono} 
$X$ and $Y$ are both split by $ab$ and since no disk contains 
four points  the interiors of $X^{ab}_p$ and $Y^{ab}_q$ 
are on opposite sides of $ab$. 
So by  the second part of lemma \ref{contri} 
$p_1p_2p_3$ and $q_1 q_2 q_3$ can only intersect  along $ab$.  Since 
the boundary of the say $q_1q_2q_3$ is  $\{q_iq_j\}$  
for this to occur either $q_1 q_2 q_3 \bigcap ab$ is precisely $a$ or $b$, 
some $q_i q_j$ is tangent to $ab$, or some $q_iq_j$ 
contains two points of $ab$.  
By the uniqueness of small geodesics in the last case $q_iq_j =ab$.
In the tangent case the fact that geodesic are the solutions 
to a second order O.D.E. once again gives us $q_iq_j = ab$. 
So in any case distinct faces can only intersect by 
sharing a vertex or 
an edge as needed.

\qed

So we have $R$ is a bijective continuous map from a compact space
hence a homeomorphism, and theorem \ref{dtri} has been proved.

\qed (theorem \ref{dtri})

%Suppose not and $R(x) =R(y)$ with $x \neq y$. Since no disk ever contains a fourth vertex  so we know neither $x$ nor $y$ can be vertices.Triangles are embedded so $x$ and $y$ must come from distinct faces $p_1 p_2 p_3$ and $q_1 q_2 q_3$ on the boundary of distinct but intersecting disks $X$ and $Y$.  Since neither $x$ nor $y$ can be vertices and the triangles are contained in their associated disks lemma \ref{lit} assures us that $X$ and $Y$ intersect in precisely two distinct  $a$ and $b$. Furthermore as in lemma \ref{mono} $X$ and $Y$ are both split by $ab$ and since no disk contains four points  the interiors of $X^{ab_p$ and $Y^{ab}_q$ are on opposite sides of $ab$.  So by  the second part of lemma \ref{contri} $p_1p_2p_3$ and $q_1 q_2 q_3$ can only intersect  along $ab$.  Since the boundary of the say $q_1q_2q_3$ is  $\{q_iq_j\}$  for this to occur either $q_1 q_2 q_3 \bigcap ab$ is precisely $a$ or $b$, some $q_i q_j$ is tangent to $ab$, or some $q_iq_j$ contains two points of $ab$.  By the uniqueness of small geodesics in the last case $q_iq_j =ab$.In the tangent case the fact that geodesic are the solutions to a second order O.D.E. once again gives us $q_iq_j = ab$. So in any case distinct faces can only intersect by sharing a vertex or an edge, and  no such $x$ and $y$ can exist. 

\end{subsection}

\begin{subsection}{Basic Properties of Delaunay Triangulations}\label{prop}
The first property worth exploring is 
the fact that a Delaunay triangulation 
is a local phenomena. 

\begin{definition}  
Call a  dense triangulation Delaunay at an edge $e$ 
if the  vertex forming the face to left  side of $e$ is   
out side the circle associated the right side's face.
\end{definition}

Notice by the monotonicity lemma (lemma \ref{mono}) 
that the property of 
being Delaunay at $e$ is equivalent to the same property 
with the sides reversed.

\begin{lemma} \label{locy}
A triangulation is Delaunay at each edge if and only if it is Delaunay.   
\end{lemma} 
{\bf Proof:}
Clearly Delaunay implies locally Delaunay.  To see the converse 
suppose it were not true and there is a triangle $t$  with an extra 
point $p$ in its associated disk.   Then $p$ is in a region of the disk to one 
side of an edge of $t$.  However since each edge is Delaunay  
this vertex cannot be the third vertex of the face on this side 
or be in its associated face.  Further more by the monotonicity 
lemma it is also in this new triangle's associated disk.
So the same situation persists for this new face.  Using this observation
one can now construct 
a sequence of  such triangles and edges with each new edge clearly closer 
to the point.  So we produce an infinite sequence of distinct edges containing 
points 
in a neighborhood  of $p$  
contradicting even local finiteness of the triangulation (let alone 
the fact it is globally finite.)

\qed 

There is another very basic local property of the dense triangulations.
I will call it the switching property.  

\begin{lemma} \label{switchy}
If two triangles  $abc$ and $abd$ of a dense triangulation 
are not Delaunay at $ab$ then $abcd$ is convex.  The triangulation formed
by switching the diagonal inside $abcd$ is Delaunay at 
$cd$ and $length(cd) < length(ab)$.
\end{lemma}
{\bf Proof:}
First we prove the convexity assertion. Note the triangles 
both are inscribed  in 
$B_{abc} \bigcap B_{abd}$ which is convex and have interiors 
on opposite sides of $ab$.  The key fact is that a shortest length 
geodesic connecting two points in  $abc$ and $abd$ must cross $ab$. 
From this observation 
$abcd$ is convex, since if it were not then then we could construct a 
shortest length geodesic form an interior point of a triangle to $ab$
 which must cross through another triangle side, contradicting the fact that 
the triangle is  convex.  

From the monotonicity lemma (lemma \ref{mono})
the ball $B$  with diagonal $ab$ is contained in 
$B_{abc} \bigcup B_{abd}$ and contains $B_{abc} \bigcap B_{abd}$, 
see figure \ref{switch}. 
Furthermore lemma \ref{mono} gives us  that
the closures of $B_{abc} \bigcap B_{abd}$ and $B$ intersect only 
at $a$ and $b$.  In particular the ball with diagonal $cd$
is strictly contained in 
$B$ with $a$ and $b$ now out side and on different sides of $cd$. 
So we immediately see the the length of $cd$ is less than the length of 
$ab$, and we may use the inflating family 
associated to $c$ and $d$ to hunt down 
$a$ and $b$ demonstrating that  the triangulation is 
indeed Delaunay at $cd$.

 \begin{figure*}
\vspace{.01in}
\hspace*{\fill}
\epsfysize = 1.5in 
\epsfbox{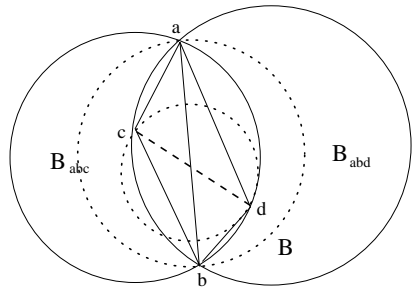}
\hspace*{\fill}
\vspace{.01in}
\caption{\label{switch}The Diagonal Switch }
\end{figure*}

\qed

It worth noting that the converse holds among convex Delaunay 
parallelograms, although the convexity assumption is now necessary.

The first corollary of the locality and switching lemmas is the proof of property \ref{delen}.

{\bf Proof of Property \ref{delen}:}
Note the set of possible triangulations is finite (it is certainly 
less than or equal to the cardinality of the vertices choose 3),
so there is at least one minimizer of $E$.  If the minimizer is not $D$ then 
by lemma \ref{locy} there is an edge which is not Delaunay, so by lemma 
\ref{switchy}  we may rearrange the triangulation at the parallelogram 
living at this edge.  Note the only edge which is changed is the diagonal 
and its length is strictly decreased hence the energy is decreased contradicting $T$'s minimality.  So the minimizer must be $D$.

\qed
 
To see property \ref{maxrad} it is useful to first 
prove the following local lemma.
 
\begin{lemma} \label{loccir}
Let $D$ be  $t_d$ a triangle in the Delaunay triangulation and let   
$c_d$ be the center of $t_d$'s associated circle.   
If $T$ is a dense triangulation relative to the same points
and $t$ is a triangle in $T$ 
such that both $c_d \in t $, then 
$rad(B_t) \geq rad(B_d)$.
\end{lemma}  

{\bf Proof:}
Assume $t_d \neq t$, 
since the result is trivial otherwise. First note 
all of the points including the 
vertices of $t$ are in $B_d^c$.  By the triangle inequality if 
$B_d$ is a proper 
subset of $B_t$ we are done. 
So we may assume the bounding circles,  $\partial B_{t}$ 
and $\partial B_{t_d}$,
intersect  in two points (lemma \ref{lit})
with lemma \ref{sir} assuring us the vertices of 
$t$ all on one side of the decomposition determined by these 
circles (or this side's closure). 
From lemma \ref{contri} $t$ and hence $c_d$ will also lies 
on this same side side of the decomposition (or its closure).  
Now look at the  inflating family 
determined by the intersection points of the 
bounding circles.
By the monotonicity lemma (lemma \ref{mono})
in order to contain greater or equal 
area on the side of the pencil where the vertices of $t$ and $c_d$ live
 $rad(B_t)$ must be greater than  $rad(B_d)$ as needed.

% since both 
%circles lie on the same side of the Circe of radius one half the 
%distance between the intersection points (the monotonicity lemma).
%So we indeed arrive at the needed inequality.   

\qed

{\bf Proof of Property \ref{maxrad}:}
First take a triangle $t_d \in D$ with $B_d$ having maximal radius in $D$.
For any dense triangulation $T$  there is a  
triangle $t$  of $T$ containing the center of $B_d$.  
Now lemma \ref{loccir} guarantees that the radius of $B_t$ is at 
least that of $B_d$ proving the corollary.   
\qed

Its also worth noting that the geometrically minimal spanning 
trees (GMST) make sense on a 
surface as well as relative neighborhood graphs (RN) and
Gabriel graphs (GG) and for the same exact reasons as in 
Euclidean space we have 
\[ GMST \subset RNG \subset GG \subset D. \]
So to construct all these things it would be 
nice to see that we can construct $D$ computationally,
and that in fact  the basic algorithm for 
changing a triangulation $T$ into $D$ works just fine.

The algorithm is essentially gradient flow of the energy $E$; namely
switch edges which are
not Delaunay as in lemma \ref{switch}.  
This procedure terminates in certainly fewer than the number 
possible triangulation steps since each triangulation has 
an associated energy and by lemma \ref{switch} the energy 
is always decreasing under this ``flow''. 
By
lemma \ref{locy} this flow must terminate in the Delaunay triangulation.

In fact this procedure must 
end before the number of vertices choose two  
steps.  This follows immediately from the following 
lemma telling us that an edge once switched 
can never be formed  again. 

\begin{lemma}
If one dense triangulation contains $ab$ and fails to be 
Delaunay at $ab$ 
then there is no dense triangulation 
containing $ab$ which is Delaunay at $ab$.
\end{lemma}  
{\bf Proof:}
This follows from the monotonicity lemma (lemma \ref{mono}) 
which guarantees that no disks in 
$ab$'s inflating family  can ever be empty.  
\qed
\end{subsection}

\end{section}

\begin{section}{Random Delaunay Triangulations}\label{euler}

As described in the chapters introduction here we tackle a the 
details needed about random Delaunay triangulations and the needed 
random variable computations.

\begin{subsection}{Basic Facts About Random Delaunay  
Triangulations}\label{randtri}
Poisson Point Process Reminders: In order to set up some  notation and 
convince the unfamiliar that nothing deep is occurring, we will now construct 
from scratch what little we need of the Poisson point process, i.e. the
probability space from the introduction.
% it is completely standard (see \cite{Du} for example), but worth doing here to set up notation and . 

 To get started it is useful to look at the space of 
sets of ordered sets of $n$ points, with a small set conveniently removed.  
This small set is  the  union of three closed measure zero sets: the set where some $x_i = x_j$ for $i \neq j$, the set where four points land on a circle of radius $r \leq \min\{ \frac{i}{6}, \tau \}$, and the set were three points land on a circle of radius exactly $\delta$.    (These removals are merely a 
technical convenience - so I will not index $\frak P$ with a $\delta$ - 
though a bit of $\delta$ has been programmed into it.) When I say 
measure zero I mean using the Riemannian product measure $dA^n$; 
that these sets  are closed and measure zero is straight forward though  
bit boring; detailed proofs can be found in 
section \ref{secsilly}.   
Let $\times^n M_-$ 
denote this full measure open subset of $\times^n M$, and let $dA^3$ be 
the Riemannian volume element restricted to this open set.  
The sets of points of interest to us   
can now be expressed as points in  
$\frak P \equiv \bigcup_{n \in \Bbb Z^{+}} \times^n M_- $.

 The measure, $\Bbb P^{\lambda}$, on this space is given by weighting  $dA^n$ on each component by $ \frac{\lambda^n}{n!} e^{{- A \lambda}} $.  This is a probability measure since the measure of  $\times^n M_-$ under $dA^n$ is $A^n$,  and so the size of $\frak P$ is $\sum_{n=0}^{\infty} \frac{ A^n \lambda^n}{n!} e^{- A \lambda} = e^{- A \lambda}  e^{ A \lambda} = 1 $.  Lastly, the measurable sets $\bf B$  will be the Borel $\sigma$-algebra.

On the $\Bbb Z^{+}$ index we have what is usually referred to as the Poisson distribution, namely $\Bbb P^\lambda( \times^n M _-) = e^{-\lambda A} \frac{A^n \lambda^n }{n!}$; which justifies (in this model) the computation of $\Bbb E_{\lambda}(F) = \lambda A$ from the introduction.   As further warm up from this view point it is useful to  explore the probability that some chunk $U$ of area $A(U)$ in $M$ is empty of points.  Strictly speaking of course $x \in \frak P$ is not a set of points in $M$.  Being careful about is perhaps less burden then any confusion resulting form it, so I will introduce the mapping $set:\times^n M_- \rightarrow 2^M$ defined as $set(p_1, \dots , p_n) = \{ p_1 , \dots ,p_n \} \subset M$.    So we are trying to find the size of the set where $set(x) \cap B = \emptyset $,or rather $\bigcup_{n \in \Bbb Z^{+}} \left( \times^n U^c \bigcap \frak P \right)$.  
% (Form now on I will suppress the intersection with $\frak P$ when discussing a sets whose measure I will finding, since the measure all that has been removed are measure zero sets.)
 For each for each $n$ using the Riemannian volume element this set has precisely the measure of $ \times^n U^c $, which  has size $(A - A(U))^n \frac{\lambda^n e^{-A \lambda}}{n!}$; and now we can sum them up to find the needed  measure is
 \[ e^{- A \lambda}  \sum_{n=0}^{\infty}  \frac{(A - A(U))^n \lambda^n}{n!} = e^{{- A \lambda}}  e^{{ (A - A(U)) \lambda}} =  e^{{- A(U) \lambda}} .\]
So we arrive at the last property of this model need for this proof.

Now we would like to beginning  proving theorem $\ref{size}$ from the introduction.  We start with: 

  \begin{proposition}  \label{shrimp}
 $\Bbb E_{\lambda}  (1_{\frak T_{\delta}} ) > 1 - c e^{d  \lambda}$, with $c$ and $d$  greater than zero.
\end{proposition}

{\bf Proof:}   The heart of the proof  is that points when distributed as above  land very densely which is what being in $\frak D_{\delta}$ means; and we know by theorem $\ref{dtri}$ that $\frak D_{\delta} \subset \frak T_{\delta}$.  Note by the triangle inequality any ball of radius $\frac{\delta}{2}$ covering the center of a ball of radius  $\delta$ must be contained in it.  So to force every ball of radius  $ \delta$  to contain a point it is sufficient to cover the surface with balls of radius $\frac{\delta}{2}$,
and then force this finite set of balls  to 
contain points.  So cover the surface with the $\frac{\delta}{2}$ 
balls and pick out a finite sub-cover with say $c$ elements.  Let $d$ be 
the minimum area amongst these $c$ balls. Now let    $\frak E_{\delta} \subset \frak P$ be the set where none of these $c$ balls is empty.  By the above observation  $\frak E_{\delta} \subset \frak D_{\delta}$, so we have $\frak T^c_{\delta} \subset \frak D^c_{\delta} \subset \frak E^c_{\delta}$.   Now $\frak E^c_{\delta}$  is precisely the union of the sets were some individual of the $c$ balls is empty - so clearly measurable. By the sub-additivity of measures, the warm up computation,  and the choices of $c$ and $d$ we now have  $\frak E^c_{\delta}$ has size less than   $   c  e^{- d A \lambda }$.  So assuming $\frak T_{\delta}$ is measurable (in fact it is open, to be seen in the final section), we have $\Bbb E_{\lambda}(1_{\frak T^c_{\delta}}) \leq  c  e^{- d A \lambda }$, as needed. 

\qed

From this we can get the needed estimates regarding the  $F$, $V$, and $E$ finishing of theorem \ref{size} from the introduction 
(up the measurability which is in fact trivial since the 
functions are continuous on $\frak P$ - 
a fact whose details 
can be found in section \ref{secsilly}). 
%the interested reader can find in \cite{Le}).

{\bf Proof of Theorem \ref{size}:}
  First let's do it for the faces. If there are $n$ points then the number of faces is usually much less than $\binom{n}{3}$.  So by the above proposition even the worse case for the area of  $\frak T_{\delta}$ is better than when  $  c  e^{-  d A \lambda } $ worth of area is crammed into the part of the space where $\binom{n}{3}$ is the largest (or $\binom{n}{2}$ for the edges or $n$ for the vertices).  In other words  letting  $S = \{ B \subset \bf{B} : \Bbb E_{\lambda}(1_{B}) = ce^{-d \lambda}\}$ we have 
\[ \Bbb E_{\lambda} \left(  \binom{n}{3} 1_{\frak T_{\delta}^c} \right) \leq  \sup_{ S} \Bbb E_{\lambda} \left(   \binom{n}{3} 1_{B} \right) = s. \] 
So really it is this quantity we estimate.  The monotonicity  of $\binom{n}{3}$ in $n$ indicates that we can realize this $sup$ with any set size   $ce^{-d \lambda} $ which fills up all the $\times^n M$ for $n \geq N_\lambda$ along with some subset of size $ ce^{-d \lambda} -  \sum_{n=N_{\lambda}+1 }^{\infty} e^{-A \lambda} \frac{A^n \lambda^n}{n!} $  in  $\times^{N_{\lambda}} M$.  To make sure we can perform this construction it is useful to observes as a scholium to the above proposition that we have $d < A$. So we naturally may stay away from the $e^{-\lambda A}$  point mass at $\times^0 M$, and can indeed realize a set of the needed size.   So for this $N_{\lambda}$ we have 

 \[ \sum_{n=N_{\lambda}+1 }^{\infty} e^{-A \lambda} \frac{A^n \lambda^n}{n!} \leq  c e^{- d  \lambda }\] and
 \[ \sum_{n=N_{\lambda} }^{\infty} e^{-A \lambda} \frac{A^n \lambda^n}{n!} >  c  e^{-d  \lambda }\]

So we now have

\[ \Bbb E(F |\frak C_\lambda^c ) \leq s \leq    \sum_{n=N_\lambda }^{\infty} \binom{n}{3} e^{-A \lambda} \frac{A^n \lambda^n}{n!} \]
\[ \leq  \frac{A^3 \lambda^3}{3!} \sum_{n = N_{\lambda} }^{\infty}  e^{-A \lambda} \frac{A^{n-3} \lambda^{n-3}}{{(n-3)}!} = \frac{A^3 \lambda^3}{6} \sum_{n = N_\lambda - 3}^{\infty}  e^{-A \lambda} \frac{A^{n} \lambda^{n}}{n!} \] 
\[\leq   \frac{A^3 \lambda^3}{3!} \sum_{n = N_\lambda - 3}^{N_\lambda }  e^{-A \lambda} \frac{A^{n} \lambda^{n}}{n!}  +  \frac{c A^3 \lambda^{3}}{6}  e^{- d \lambda}. \]

The second term decays faster than any polynomial, so  we are reduced to seeing the terms in the form $ {(A \lambda)}^k e^{-A \lambda} \frac{{(A \lambda)}^{N_{\lambda}-l}}{(N_\lambda -l)!} $  decay quickly, with $l$ and $k$  non-negative integers.  

To approach these terms first we use that our upper bound gives us  $  \frac{{(A \lambda)}^{N_\lambda + 1} }{(N_{\lambda + 1})!} <  c  e^{(A - d)  \lambda } $, or rather  $A \lambda <  {((N_{\lambda} +1)! c  e^{(A - d)  \lambda } )}^{\frac{1}{k+1}}$.

In particular  
\[  \frac{{A \lambda}^{N_{\lambda}-l}}{(N_{\lambda}-l)!}  <  \frac{{((N_{\lambda} + 1)! c  e^{(A - d)  \lambda } )}^{\frac{N_{\lambda}-l}{N_{\lambda}+1}}}{(N_{\lambda}-l)!} \]
\[ < \frac{ (N_{\lambda} + 1)!  c  e^{(A - d)  \lambda }}{(N_{\lambda}-l)!} 
< {(N_{\lambda}+1)}^{l+1} c  e^{(A - d)  \lambda }. \]

Applying this to our term gives the quite manageable estimate
 \[ {(A \lambda)}^k e^{-A \lambda} \frac{{(A \lambda)}^{N_\lambda -l}}{(N_\lambda -l)!} 
< (N_\lambda+1)^{l+1} {(A \lambda)}^k   e^{- d \lambda }. \]

Now we see these terms would decay faster than any polynomial if the $N_{\lambda}$ 
grew (with $A\lambda$) no faster than a polynomial.  To see this is so, note that our estimate in the other direction gives us  
$ \sum_{n=N_{\lambda} }^{\infty}  \frac{A^n \lambda^n}{n!} >  c  e^{(A-d)  \lambda }.$
Now $d < A$ so we might hope $N_\lambda$ must remain quite small for this to be so large.
 In fact this is easy to show; we may even observe a rather extreme fact that  $N_{\lambda}$ cannot have a subsequence grow even as fast as  ${(A \lambda)}^3$  with

\[ \lim_{\lambda = \frac{m}{A} \rightarrow \infty}\sum_{n = {(\frac{A m}{A})}^3}^{\infty}   \frac{A^n \lambda^n}{n!} =  \lim_{m \rightarrow \infty}\sum_{n = {m}^3}^{\infty}   \frac{m^n}{n!} =  0 \]

 To see this note that Stirling's formula tells us that $n! > \frac{1}{C} n^n e^{-n} {(2 \pi n)}^{\frac{1}{2}}$, so assuming $m > e$ we have

\[ \lim_{m \rightarrow \infty} \sum_{n = m^3}^{\infty}   \frac{m^n}{n!}
\leq C \lim_{m \rightarrow \infty} \sum_{n = m^3}^{\infty} {\left( \frac{e m }{n} \right)}^n  {(2 \pi n)}^{-\frac{1}{2}} \]
\[ \leq  C \lim_{m \rightarrow \infty}\sum_{n = {m}^3}^{\infty} {\left( \frac{1 }{m} \right)}^n \leq C \lim_{m \rightarrow \infty}\sum_{n = {m}^3}^{\infty} {\left( \frac{1 }{e} \right)}^n .   \]

This goes to zero since the geometric series at $\frac{1}{e}$ converges.
 
The $E$ and $V$ cases are virtually identical - simply do the same thing with $  \binom{n}{2}$ and $n$ rather than $\binom{n}{3}$.   So once again up to the measurability we are done.

\qed

\end{subsection}

  Let $V_{\delta}$ be the set of
ordered triples in $\times^3 M$  which are 
on circles of radius less than $\delta$ 
%(note the irrelevant measure zero change of this set from its definition in th%e introduction)
, and let $r$ be a  measurable function on $\times^3 M$ with its  restriction 
to $V_{\delta}$ in 
  $ L^1(V_{\delta})$ and which  is  
symmetric under any  permutation of the coordinates.
 Let $R$ be a random variable on $\Bbb P_{\lambda}$ which is 
given as a $R({\bf p})  = \sum_{t \in K_{\bf p}} r(t) $ where by $t$ 
we mean the ordered triple in $\bf{p}$  corresponding to the face 
$t \in K_{\bf p}$.
(that this function and the functions introduced below  are  measurable can be easily seen; see section \ref{secsilly} if there is any confusion).
%\cite{Le} if there is any confusion).  
For such a random variable we have its expect value given by: 

\begin{theorem} \label{ranvar}
 With the above notation we have 
\[ \Bbb E_{\lambda} (R) =  
\frac{\lambda^3}{6} \int_{V_{\delta}} r(y) e^{- a(y) \lambda}  dA^3(y). \]
\end{theorem}

{\bf Proof:}
%With the probabilistic context set up as in section 2.3 equation  
%(\ref{ex1}) from the introduction
%  can be most easily handled  by  producing a handy 
%formula for $F$.
Given a set $s = \{i_1,i_2, i_3 :i_1 < i_2 < i_3   \} \subset \{1, ..., n \} $
 let $\pi_s$ be the projection mapping of $M_1 \times \dots \times M_n$ 
onto  $ M_{i_1} \times M_{i_2} \times   M_{i_3}$. Now note that every 
triple will be uniquely represented as $set \cdot \pi_s({\bf p})$ by one 
of these $s$.  Define a function $f_s(x)$ 
to be one if  ${\bf p} \in \times^3 M_{-}$ and the triple 
$  set \cdot \pi_s({\bf p}) $ is on a disk of radius less than $\delta$ 
and has its uniquely associated open disk is empty of 
points in the configuration $set({\bf p})$; and zero otherwise.     
Let  $R_n({\bf p}) = \sum_{s=1}^{\binom{n}{3}}r(\pi_s({\bf p})) f_s({\bf p})$ 
on  ${\times}^n M$, and  
note the above random variable $R$ is precisely the function defined on  
$\frak{P} = \bigcup_{n \in \Bbb Z^{+}} \times^n M_- $ 
which is $R_n$ on each $\times^n M_- $.  
%(Also on $\frak P$ the $r_s$ are, like $R$, are measurable
%(which is easily seen and a carefully treatment of can be found in \cite{Le}).

This formula buys us a better look at $\Bbb E_{\lambda} (R)$.
 First noting that $\times^n M_-$ and $\times^n M $ 
differ by a measure zero set
and breaking up the integral into the pieces over the disjoint 
pieces of the space we have

%previously $M_-$ was always used, why?
   
 \[ \Bbb E_{\lambda} (R) = \sum_{n=3}^{\infty} \int_{\times^n M}
 R_n({\bf p}) dA^n \frac{\lambda^n}{n!} e^{{- A \lambda}}. \]
Using our formula for  $R_n$ and the linearity of the integral we have
  
\[  \Bbb E (R) =  \sum_{n=3}^{\infty} 
\sum_{s=1}^{\binom{n}{3}}  \int_{\times^n M} 
r(\pi_s({\bf p})) f_s({\bf p}) dA^n \frac{\lambda^n}{n!} e^{{- A \lambda}}. \]
Now by symmetry each $r(\pi_s({\bf p})) f_s({\bf p})$
 (for a fixed $n$) has the 
same integral so we may write $\Bbb E(F)$ as 
  
\[  \sum_{n=3}^{\infty} \binom{n}{3}  
\int_{\times^n M} r(\pi_s({\bf p})) f_s({\bf p}) dA^n \frac{\lambda^n}{n!} 
e^{{- A \lambda}} \] 
\[ =  e^{{- A \lambda}}\frac{ \lambda^3}{6} \sum_{n=3}^{\infty}  
\frac{\lambda^{n-3}}{(n-3)!} \int_{\times^n M} 
r(\pi_s({\bf p})) f_s({\bf p})  dA^n. \]
 
To do this integral note $ r(\pi_s({\bf p})) f_s({\bf p}) $ 
is a function which is zero off the set  with two properties.  
First of all, the triple $set(\pi_s({\bf p}))$ is on a ball of radius less 
than $\delta$ - call it 
$B(set(\pi_s({\bf p})))$; and, secondly, all 
of the other points in $set({\bf p})$ are 
in $B(set(\pi_s({\bf p})))^c$.  So this set can be given by 
$\bigcup_{y \in V_{\delta}} 
\left( y \times \left(\times^{n-3} B(y)^c \right) \right)$ 
with  $y$ denoting the coordinates onto which 
 $\pi_s$ projects. 
%(think about this $\frak{P}$).  
Note  that on this set our random variable is 
 $ r(\pi_s(x))$ (at least off a measure zero subset where it is zero) and so 
is constant on each of the disjoint  
$ y \times \left(\times^{n-3} B(y)^c\right) $ pieces of the above set.
%(This set is measurable since $f_s$ is see \cite{Le}.)
So Fubini's theorem (with the $r \in L^1(V_{\delta})$ assumption) tells us  
\[ \int_{\times^n M}  r(\pi_s({\bf p})) f_s({\bf p}) dA^n
 =  \int_{ V_{\delta}} r(y) 
\left( \int_{\times^{n-3} M} 1_{\times^{n-3} B(y)^c} dA^{n-3} \right) dA^3(y)
 \]
\[ =   \int_{V_{\delta}} r(y)  (A -a(y))^{n-3} dA^3(y), \]
where $a(y)$ is the area of the ball function for $y$ with 
$set(y)$ containing three distinct points and zero other wise.

%the indicator function on the set with two properties.  So $f_s $ is the indicator%ator of   $\bigcup_{y \in V_{\delta}} \left( y \times \left(\times^{n-3} B(y)^%c \right) \right) \bigcap \frak P$.  Since $f_s$ is measurable this set is  me%asurable in $\bf B$, and  so (at the very least)  $\bigcup_{y \in V_{\delta}} %\left( y \times \left(\times^{n-3} B(y)^c \right) \right)$ is in the completion%n of the Borels of $\times^n M$ with the same measure. 

Plugging this into the above computation, and using Fubini's theorem once again, we have
\[  \Bbb E (F) = e^{{- A \lambda}} \frac{\lambda^3}{6} \sum_{n=3}^{\infty}  
\frac{\lambda^{n-3}}{(n-3)!}\int_{V_{\delta}}  
r(y) (A -a(y))^{n-3} dA^3(y) = \] 
\[\frac{e^{{- A \lambda}} \lambda^3}{6}  \int_{V_{\delta}} r(y)
 e^{(A- a(y)) \lambda}  dA^3(y)  = 
\frac{\lambda^3}{6} \int_{V_{\delta}} r(y) e^{- a(y) \lambda}  dA^3(y). \]

{\bf Q.E.D.}

Note in particular when $r = 1$ we have equation (\ref{ex1}) 
from the introduction.

%To further compute this we need to coordinatize  $V_{\delta}$, which we do in the next section.  Then in section 3.2 we take a look in these coordinates and note that we have in fact found the Gauss-Bonnet theorem! 

%, where $ y \times \left(\times^{n-3} B(s)^c \right)$ really means to place $y$ in the coordinates specified by $s$

\begin{subsection}{The Geometric Coordinates}\label{coords}

Now its time to carefully set up the coordinates form the introduction, 
and to derive equation (\ref{ex2}) from the introduction.  
Here we will parameterize a  full measure subset 
$(V_{\delta})_- \subset V_{\delta}$  via the set  
\[ L_{\delta} = \times^3 S^1_- \times (0,\delta) \times \{ M -\{points\} \}. \]   This mapping will be  based on the
exponential mapping $exp_p : T_pM \rightarrow M $ 
(see \cite{Do}), and its definition requires a  orthonormal frame 
$f = \{e_1, e_2 \}$ - which can be defined at all but a finite 
number of points ($\{ points \}_f$) of $M$.
%(see section 1.1).  
Identifying  $\times^3 S^1$ with the triple of variables
  $(\theta_1  \mod 2 \pi,\theta_2 \mod 2 \pi,\theta_3 \mod 2 \pi)$,  
denoted as $\vec{\theta}$, we can explicitly define our mapping 
$\Phi_f:L_{\delta} \rightarrow (V_{\delta})_-$ as

\[ \Phi_f(\vec{\theta}, r ,z) = (exp_z(r v(\theta_1)), 
exp_z(r v(\theta_2)),exp_z(r v(\theta_3)) ), \]
where the points excluded from $M$ are precisely $\{ points \}_f$.

Clearly we land inside the triples on circles of radius 
less than $\delta$ and so are well defined. The differentiability of 
the exponential map guarantees 
$\Phi_f$ is a differentiable on $L_{i}$.  
We would of course like to say more than that we have a differentiable map. 
As we know by lemma $\ref{lit}$ a triple uniquely determines its 
disk when the radius is less than $\delta$, so the $z$ and $r$ 
coordinates are uniquely determined.    Also $set^{-1} 
(set \cdot \Phi_f (\vec{\theta}, r ,z))$ corresponding to the 
6 permutation of a distinct triple is hit precisely by the 6 distinct 
images under  $\Phi_f$ of $ (set^{-1}\{ \theta_1 , 
\theta_1, \theta_1 \},r,z)$ with $r$ and $z$ fixed.  
So $\Phi_f$ is injective.  Also $\Phi_f$'s image is almost everything 
since  any triple in $(V_{\delta})_-$ not centered at $z \in \{ points \}_f$ 
is hit; a set which his of course measure zero (see section \ref{secsilly}).
%(see \cite{Le} for the details).

Being bijective (to a full measure set) and differentiable, $\Phi_f$ 
forms a re-parameterization of  a full measure subset of 
$V_{\delta}$; and with it we  may continue 
the computation from the previous section  this section finding  
\[  \Bbb E (R) = \frac{\lambda^3}{6} \int_{V_{\delta}} 
r(y) e^{- a(y) \lambda}  dA^3(y) \]
\[ = \frac{\lambda^3}{6} \int_{\Phi_f(L_{\delta})}r(y)  e^{- a(y) \lambda}  dA^3
=  \frac{\lambda^3}{6} \int_{L_{\delta}} 
r(\phi_f(\vec{\theta},r,p))e^{- \lambda 
a(\phi_f(\vec{\theta},r,p))}\Phi_f^*dA^3
. \]

The remainder of this sub-section will be dedicated to  
finding an explicit formula for 
$\Phi_f^*(dA^3) $.  Since we are pulling back a Riemannian volume 
form, understanding how the metric pulls back 
will determine how the the volume pulls back; and this will be our strategy.
%So we need  to compute the inner-products
%between elements of the tangent space in $L_{\delta}$ using the $M$ metric.   
%If one 
Pulling back the metric in an arbitrary 
frame  is fortunately not needed; 
in fact the expression for the volume form at $(\vec{\theta},p,r)$ 
when pulled back by 
$\Phi_f$ depends only on the orthonormal basis choice at $p$. Letting 
$dA^3$ 
be the the volume form on $ V_{\delta}$ (since $V_{\delta}$ is an open 
set in $\times^3 M$) we can express this as:
\begin{lemma} \label{orth}
Suppose the frames $f$ and $g$ agree at $p \in M$, then  
\[ \Phi_f^{*} \left(dA^3 \right)(\vec{\theta},r, p)  = 
\Phi_g^{*} \left(dA^3 \right)(\vec{\theta},r, p). \] 
\end{lemma}
{\bf Proof:}
First notice that given two smooth  orthonormal frames $f = \{e_1, e_2 \}$ and $g = \{d_1, d_2 \}$   defined on a simply connected open set of $M$ (call it $U$), they can be orthogonal compared via a smoothly varying element of $O(1)$.   $U$ being simply connected allows us to lift this mapping from $O(1)$ to its cover $\Bbb R^1$.  
 Otherwise said: there is a differentiable function $\theta_{f,g}(z)$ such that 
\[ \left[ \begin{array}{l} e_1 \\ e_2 \\ \end{array} \right] = \left[ \begin{array}{ll} \cos(\theta_{fg}) &  -\sin(\theta_{fg}) \\  \sin(\theta_{fg}) & \cos(\theta_{fg}) \\  \end{array} \right] \left[ \begin{array}{l} d_1 \\ d_2 \\ \end{array} \right]. \]

Let $F$ be the mapping 
\[ F:\times^3 S^1_- \times (0,\delta) \times \{ U -\{points\}  \} \rightarrow \times^3 S^1_- \times (0,\delta) \times \{ U -\{points\} \} \]  such that  
\[ F:(\vec{\theta}, r ,z) = (\theta_1 + \theta_{fg}, \theta_2+ \theta_{fg}, \theta_3+ \theta_{fg}, r ,z). \]
Using this mapping note that our mapping satisfies  $ \Phi_f = \Phi_g \cdot F $
on $\times^3 S^1_- \times (0,\delta) \times \{ U -\{points\}  \}$.

Now simply note that  using $d(\theta_{fg}(z)) \wedge dA = 0 $ we have
\[F^{*}( v(\vec{\theta}, r ,z) dA \wedge dr  \wedge d \theta_1 \wedge d \theta_2  \wedge  d \theta_3) \]
\[ =   v(\theta_1 + \theta_{fg}, \theta_2+ \theta_{fg}, \theta_3+ \theta_{fg}, r ,z) dA \wedge dr \]
\[  \wedge d (\theta_1 + \theta_{fg}(z)) \wedge d (\theta_2 + \theta_{fg}(z))  \wedge  d (\theta_3 + f(z)(z)) = \]
\[    v(\theta_1 + \theta_{fg}, \theta_2+ \theta_{fg}, \theta_3+ \theta_{fg}, r ,z) dA \wedge dr  \wedge d \theta_1 \wedge d \theta_2  \wedge  d \theta_3 .\]

Assuming the frames agree at $p$, we have $\theta_{fg}(p) = 0 $. So by the above formula $F^{*}$ acts as the identity at this point - giving 
\[ \Phi_f^{*} \left(dV \right)(\vec{\theta},r, p)  = F^{*}\Phi_g^{*} \left(dV \right)(\vec{\theta},r, p)
=  \Phi_g^{*} \left(dV \right)(\vec{\theta},r, p), \]
as needed.

\qed

To actually compute the inner products it is useful
to canonically identify all the tangent spaces 
near $p$ as as $\vec{\theta}$ and $r$ vary with a 
frame determined by the triple product of the normal coordinates at $p$, 
$N(p,\vec{\theta},r)$, via 
\[L_p(\vec{\theta},r) = \times^3  N_{*}: \times^3 \Bbb E^2 
\rightarrow T_{\Phi_f(\vec{\theta},r,p)}V_{\delta}.\]

While not an orthonormal frame a certain amount of the Euclidean geometry 
in  $\times^3 \Bbb E^2$ 
is preserved by this frame.  For starters each of the $\Bbb E^2$ copies in
$\times^3 \Bbb E^2$ is orthogonal to the others,  since 
$ V_\delta$ is an open subset of $\times^3 \Bbb M$ in its product metric, 
and this mapping is respecting the product structure. 
Denote the vectors in this $ \times^3 \Bbb E^2$ as  $(w_1,w_2,w_3)$. 
Further more by Gauss's lemma 
each of these $w_i$ components decomposes 
orthogonally (when $0< r \neq i$) into $v(\theta_i)(p)$ 
and $v^{\perp}(\theta_i)(p)$ - with respect to the both the Euclidean metric  
and the surface metric.
  So for each triple of angles we may in fact represent 
the vectors with the following orthogonal decomposition 
 \[ ( a_1 v(\theta_1)(p) + b_1r v^{\perp}(\theta_1)(p), a_2 v(\theta_2)(p) + b_2 r v^{\perp}(\theta_2)(p), \] 
\[ a_3 v(\theta_3)(p) + b_3 r v^{\perp}(\theta_3)(p)) . \] 
To understand the lengths of these vectors in this decomposition
it is necessary to remind our selves about Jacobi Fields.

We can now get a grip on our needed vector lengths.

\begin{lemma}
\label{vec}
Suppose $f$ is a geodesic frame at $p$ then using the above notation we have:  
\begin{itemize}
\item
\[ L^{-1} \cdot {(\Phi_f)}_{*} \left( \frac{\partial }{\partial \theta_1} \right) = (r v(\theta_1)^{\perp}(p), 0 , 0). \] 
 Similarly for $\theta_2$ and $\theta_3$. 
\item
 \[ L^{-1} \cdot {(\Phi_f)}_{*} \left( \frac{\partial }{\partial r}\right)  = ( v(\theta_1)(p),  v(\theta_2)(p) ,  v(\theta_3)(p)).\]
\item
%Using coordinates $(z_1,z_2)$ on $M$ such that  $\frac{\partial }{\partial z_i} = e_i$ we have 

\[ L^{-1} \cdot {(\Phi_f)}_{*} \left( e_1 \right) =  (\cos(\theta_1) v(\theta_1),\cos(\theta_2) v(\theta_2), \cos(\theta_3) v(\theta_3)) + \]
\[   (-\frac{r h_{\theta_1}}{j_{\theta_0}} \sin(\theta_1) v(\theta_1)^{\perp},- \frac{r h_{\theta_1}}{j_{\theta_0}} \sin(\theta_2) v(\theta_2)^{\perp}, - \frac{r h_{\theta_1}}{j_{\theta_0}} \sin(\theta_3) v(\theta_3)^{\perp}).\]
Similarly for $e_2$. 
\end{itemize} 
\end{lemma} 
$\bf{Proof:}$
The idea for all these computations is identical.  For each component one finds a curve $\Gamma(s)$ such that $\frac{d}{ds} \Gamma(s) = \pi_{\star} v$ and then one notes that in fact $\Gamma(s) = \Gamma(r,s)$ for some two parameter family of geodesics - forcing our vector, $v$ , to be Jacobi fields. At this point one uses the Jacobi lemma (lemma $\ref{jac}$ in section \ref{georem}) to observe the above formulas.

For example for ${(\Phi_f)}_{*} \frac{\partial }{\partial \theta_1} $ is by our Jacobi field observations is  
  \[ \frac{d}{d s} (exp_{p}(rv(s)),exp_{p}(rv(\theta_2)),exp_{p}(rv(\theta_3))) |_{\theta_1} . \] 
 By the discussion preceding lemma $\ref{jac}$ this   is the Jacobi field with initial conditions $J(0) =0$ and $\frac{d J}{dr}(0) = v^{\perp}(\theta_1)$ - located in the in the first component of $M^3$. So we have a description of it from part 3 of lemma $\ref{jac}$ as $(r v(\theta_1)^{\perp}(p), 0 , 0)$.

%\begin{figure*}
%\vspace{.01in}
%\hspace*{\fill}
%\epsfysize = 1 in 
%\epsfbox{angle.eps}
%\hspace*{\fill}
%\vspace{.01in}
%\caption{\label{ang} The Angle Fields}
%\end{figure*}

  Now let's  play the same game for the radial direction and   note that the image of ${(\Phi_f)}_{*} \frac{\partial}{\partial r}$ is 
 \[ \frac{d}{d s} (exp_{p}((r + s) v(\theta_1),exp_{p}((r + s) v(\theta_2) , exp_{p}((r + s) v(\theta_3))|_{s =0}, \] which by part 1 of lemma $\ref{jac}$ is  $J(r) = (v(\theta_1)(p),(v(\theta_2)(p),(v(\theta_3)(p))$.

% \begin{figure*}
%\vspace{.01in}
%\hspace*{\fill}
%\epsfysize = 1 in 
%\epsfbox{radiu.eps}
%\hspace*{\fill}
%\vspace{.01in}
%\caption{\label{rad} The Radius Field}
%\end{figure*}

 Finally to compute the vector in the  
$e_1$ direction  use the geodesic pointing in that direction 
and vary along it. This is the situation of the example preceding 
lemma $\ref{jac}$, and we have ${(\Phi_f)}_{*} e_1 $ 
described as 
 \[ \frac{d}{ds}( exp_{exp_p(s e_i)}(r v(\theta_1)), exp_{exp_p(s e_i)}(r v(\theta_2)), exp_{exp_p(s e_i)}(r v(\theta_3))). \] 

Since the frame we are using is geodesic at $p$ we have $\frac{Dv(\theta)}{dr} = 0$ so  $\frac{D J}{dr}(0) = 0$ (by the example preceding lemma $\ref{jac}$ - once again). By the discussion preceding the Jacobi lemma $\ref{jac}$, we may decompose this field by decomposing $J(0) = e_1(p)$ into $v(\theta_i)$ and $v(\theta_i)^{\perp}$ in each component; so we have  this Jacobi field has one summand corresponding to  
\[ J(0) = (<v(\theta_1) ,e_1>  v(\theta_1) \]
\[ ,<v(\theta_2) ,e_1> v(\theta_2) , <v(\theta_3) , e_1> v(\theta_3) )(p)\] and $\frac{D J}{dr}(0) = 0$.  This field is dealt with by the first part of lemma $\ref{jac}$, and gives the first set of vectors. The other summand in the decomposition corresponds to   
\[ J(0) = ( <v(\theta_1)^{\perp} , e_1> v(\theta_1)^{\perp}\] 
\[, < v(\theta_2)^{\perp} ,e_1>  v(\theta_2)^{\perp}, <v(\theta_3)^{\perp} ,e_1>  v(\theta_3)^{\perp})(p) \]
 with $\frac{D J}{dr}(0) = 0$ - which is dealt with in by the second part of lemma $\ref{jac}$.  These formulas are precisely the last of the needed inner-product relationships.
% \begin{figure*}
%\vspace{.01in}
%\hspace*{\fill}
%\epsfysize = 1 in 
%\epsfbox{center.eps}
%\hspace*{\fill}
%\vspace{.01in}
%\caption{\label{cent} The center Fields}
%\end{figure*}

\qed

With the use of these lemmas  we  accomplish our goal of computing the pullback of the volume form in $V_{\delta}$ coordinates and derive  equation (\ref{ex2}). 

\begin{proposition}
\label{calc}
\[ \Phi_f^*(dA^3) = \nu( \vec{\theta})j_{\theta_1}j_{\theta_2}j_{\theta_3} d \vec{\theta}  \wedge dr \wedge dA \]
with  $\nu( \vec{\theta_1})$ the area of a triangle with vertices at $\{ v(\theta_i) \}$ in the Euclidean plane.

\end{proposition}
$\bf{Proof:}$
Let's compute the form at a point $(\vec{\theta},r,p)$.  
Using the normal coordinates for the $M$ variables note we have  
$\frac{\partial }{\partial z_i} = e_i$.  We need to find 
$ {(det(g_{ij})}^{\frac{1}{2}}  dz_1 \wedge d z_2 \wedge d \vec{\theta}  
\wedge dr$.  To do so first note by lemma $\ref{orth}$ that we may use 
the geodesic frame at $p$ with the initial vectors $\{e_1(p), e_2(p) \}$.   
Recalling that   $< v(\theta_i),  v(\theta_i)>= 1$, and that from lemma 
$\ref{jac}$  in section 1.1 $<r v^{\perp}(\theta_i), 
r v^{\perp}(\theta_i)>= j^2_{\theta_i}$.  Using the above 
observation we now  have can easily compute the determinant.

Letting
\[  \nu( \vec{\theta_1}) = |\sin(\theta_2 - \theta_1) + \sin(\theta_3 - \theta_2)+\sin(\theta_1 - \theta_3)|\] 
\[= \left|\frac{1}{4} \sin\left(\frac{\theta_2 - \theta_1}{2}\right)  \sin\left(\frac{\theta_3 - \theta_2}{2}\right) \sin\left(\frac{\theta_3 - \theta_1}{2}\right)\right|, \]
 we find
\[ det(g_{ij})= 
\nu^2( \vec{\theta})j^2_{\theta_1}j^2_{\theta_2}j^2_{\theta_3}. \]
Take its square root to find  
\[ \Phi^*(dA^3) =   \nu(\vec{\theta})j_{\theta_1}j_{\theta_2}j_{\theta_3}   d \vec{\theta}  \wedge dr \wedge dz_1 \wedge d z_2.   \]

Now $dz_1 \wedge d z_2$ is precisely $dA$ at in these coordinates 
(at $p$) - so we have our needed equality. 
%The resulting inequality can be found in the lemma 
%$\ref{jac}$ in section 1.2.

\qed

Actually it is  useful to witness the matrix forming 
the determinant in the above proposition and 
hence find that it is relatively easy to compute. 
Note 
{\tiny \[ \left[ \begin{array}{llllll}
< \frac{\partial}{\partial \theta_1} ,  \frac{\partial}{\partial \theta_1} > &
< \frac{\partial}{\partial \theta_1} ,  \frac{\partial}{\partial \theta_2} > &< \frac{\partial}{\partial \theta_1} ,  \frac{\partial}{\partial \theta_3 }> &< \frac{\partial}{\partial \theta_1 },  \frac{\partial}{\partial r} > &
< \frac{\partial}{\partial \theta_1} ,  \frac{\partial}{\partial z_1} > &
< \frac{\partial}{\partial \theta_1 },  \frac{\partial}{\partial z_2 }> \\
< \frac{\partial}{\partial \theta_2 },  \frac{\partial}{\partial \theta_1} > &
< \frac{\partial}{\partial \theta_2} ,  \frac{\partial}{\partial \theta_2} > &
< \frac{\partial}{\partial \theta_2} ,  \frac{\partial}{\partial \theta_3} > &
< \frac{\partial}{\partial \theta_2} ,  \frac{\partial}{\partial r} > &
< \frac{\partial}{\partial \theta_2} ,  \frac{\partial}{\partial z_1} > &
< \frac{\partial}{\partial \theta_2} ,  \frac{\partial}{\partial z_2} > \\
< \frac{\partial}{\partial \theta_3 },  \frac{\partial}{\partial \theta_1} > &
< \frac{\partial}{\partial \theta_3} ,  \frac{\partial}{\partial \theta_2} > &
< \frac{\partial}{\partial \theta_3} ,  \frac{\partial}{\partial \theta_3} > &
< \frac{\partial}{\partial \theta_3} ,  \frac{\partial}{\partial r} > &
< \frac{\partial}{\partial \theta_3} ,  \frac{\partial}{\partial z_1} > &
< \frac{\partial}{\partial \theta_3} ,  \frac{\partial}{\partial z_2} > \\
< \frac{\partial}{\partial r},  \frac{\partial}{\partial \theta_1} > &
< \frac{\partial}{\partial r} ,  \frac{\partial}{\partial \theta_2} > &
< \frac{\partial}{\partial r} ,  \frac{\partial}{\partial \theta_3} > &
< \frac{\partial}{\partial r} ,  \frac{\partial}{\partial r} > &
< \frac{\partial}{\partial r} ,  \frac{\partial}{\partial z_1} > &
< \frac{\partial}{\partial r} ,  \frac{\partial}{\partial z_2} > \\
< \frac{\partial}{\partial z_1},  \frac{\partial}{\partial \theta_1} > &
< \frac{\partial}{\partial z_1} ,  \frac{\partial}{\partial \theta_2} > &
< \frac{\partial}{\partial z_1} ,  \frac{\partial}{\partial \theta_3} > &
< \frac{\partial}{\partial z_1} ,  \frac{\partial}{\partial r} > &
< \frac{\partial}{\partial z_1} ,  \frac{\partial}{\partial z_1} > &
< \frac{\partial}{\partial z_1} ,  \frac{\partial}{\partial z_2} > \\
< \frac{\partial}{\partial z_2},  \frac{\partial}{\partial \theta_1} > &
< \frac{\partial}{\partial z_2} ,  \frac{\partial}{\partial \theta_2} > &
< \frac{\partial}{\partial z_2} ,  \frac{\partial}{\partial \theta_3} > &
< \frac{\partial}{\partial z_2} ,  \frac{\partial}{\partial r} > &
< \frac{\partial}{\partial z^2} ,  \frac{\partial}{\partial z_1} > &
< \frac{\partial}{\partial z^2} ,  \frac{\partial}{\partial z_2} > \\
\end{array} \right] \] }
{\tiny \[ = \left[ \begin{array}{llll}
j^2_{\theta_1} &
0 & 0 & 0 
 \\
0 & j^2_{\theta_2}& 0 & 0 
  \\0& 0& j^2_{\theta_3} & 0 
 \\
0 & 0 & 0 & 3  
 \\
-j_{\theta_1} h_{\theta_1} \sin(\theta_1) &
-j_{\theta_2} h_{\theta_2} \sin(\theta_2)  &
-j_{\theta_3} h_{\theta_3} \sin(\theta_3)  &
\sum_{k=1}^{3} \cos(\theta_k)   \\
j_{\theta_1} h_{\theta_1} \cos(\theta_1) &
j_{\theta_2} h_{\theta_2} \cos(\theta_2)  &
j_{\theta_2} h_{\theta_2} \cos(\theta_2)  &
\sum_{k=1}^{3} \sin(\theta_k)  \\ \end{array} \right. \]

\[  \left.  \begin{array}{ll}
-j_{\theta_1} h_{\theta_1} \sin(\theta_1)  &
j_{\theta_1} h_{\theta_1} \cos(\theta_1) \\
-j_{\theta_2} h_{\theta_2} \sin(\theta_2)  &
j_{\theta_2} h_{\theta_2} \cos(\theta_2) \\
-j_{\theta_3} h_{\theta_3} \sin(\theta_3) &
j_{\theta_3} h_{\theta_3} \cos(\theta_3) \\
\sum_{k=1}^{3} \cos(\theta_k) &
\sum_{k=1}^{3} \sin(\theta_k)  \\
 \sum_{k=1}^{3}( h_{\theta_k}^2 \sin^2(\theta_k) + \cos^2 (\theta_k)) &
\sum_{k=1}^3 (\sin(\theta_k) \cos(\theta_k) - h_{\theta_k}^2 \sin(\theta_k) \cos(\theta_k))  \\
\sum_{k=1}^3 (\sin(\theta_k) \cos(\theta_k) -  h_{\theta_k}^2 \sin(\theta_k) \cos(\theta_k)) &
 \sum_{k=1}^{3}( h_{\theta_k}^2 \cos^2(\theta_k) +  \sin^2 (\theta_k)) \\
\end{array} \right] \] }
so one can pull the $j_{\theta_k}$ out of the matrix and 
the symmetries of the remaining matrix make the computation 
of the determinant relatively transparent.   (Trying to see 
this first hand proves useful when one begins 
exploring the three 
dimensional computation.)
\end{subsection}
  
\begin{subsection}{The Euler-Gauss-Bonnet-Delaunay Formula}\label{dercomp}
It is now time to compute the needed random variables. 
 All the random variables encountered have expected values expressible as
\[ \Bbb E_{\lambda} (R) = \frac{\lambda^3}{6} \int_{L_{\delta}} r^3 
f(r,p,\vec{\theta})
  e^{- a(r,p,\vec{\theta}) \lambda} dr  d \vec{\theta} 
dA  \]
where 
\[f(r,p,\vec{\theta}) =  \left(f_0(p,\theta)
 + f_1(p,\theta) r + \right. \left. f_2(p,\theta) r^2 \ln(r) +  f_3(p,\theta) r^2 + 
O(\ln(r) r^3) \right)\]
with the $f_i$ bounded and measurable functions on $L_{\delta}$.
In this expression and throughout out this paper remainder 
a function denoted 
$O(f(r))$ will mean a  
function  $O(p,\vec{\theta},r)$ on 
 $M \times [0,2 \pi]^3 \times (0,\delta)$ satisfying 
$|O(p,\vec{\theta},r)| \leq
C f(r)$ for some constant $C$.

To compute the above  expected value first note that one can Taylor 
expand in $r$ the area of a ball of radius $r$ at $p$ as 
$a(r,p) = 
\pi r^2 -\frac{k(z) \pi}{12}r^4 +O(r^5)$ (in fact this follows immediately 
form the o.d.e. describing the $v^{\perp}(\theta)$ Jacobi fields and the fact the surface is compact). 
This formula 
allows us to choose  constants $0 < \epsilon < 1$ and $\rho_M >0$
such that 
\[ e^{- \lambda a(p,\theta,r)} < e^{- \lambda \pi r^2 (1 - \epsilon)} \]  
for all $\rho_M  > r  > 0$ and at all $p \in M$. 
%Using this $\rho_M$ as our decision radius we have: 

\begin{theorem}\label{exp} 
Using the above notation,
letting $c = \int_{0}^{\infty} \ln(r) r^2 e^{-r^2} dr $,
 and using a decision radius smaller than $\rho_M$ 
we have $ \Bbb E_{\lambda} (R)$ equals

\[ \frac{\lambda}{6}\int_{M} \int_{S^1 \times S^1 \times S^1}
\frac{f_0(p,\theta)}{2 \pi^2}
 d \vec{\theta} dA  + \frac{\ln(\lambda) \lambda^{\frac{1}{2}}}{6}\int_{M} 
\int_{S^1 \times S^1 \times S^1} 
\frac{ 15 f_2(p,\theta) }{16 \pi^2} 
d \vec{\theta} dA \]
\[ + \frac{\lambda^{\frac{1}{2}}}{6}\int_{M} \int_{S^1 \times S^1 \times S^1}
\left( \frac{ 15 f_1(p,\theta) }{16 \pi^2} +
  \frac{c  f_2(p,\theta)}{\pi^{\frac{5}{2}}}+
 \frac{105 k(p) f_1(p,\theta)}{12 \cdot 32 \pi^3} \right)
 d \vec{\theta} dA \]

\[ \frac{1}{6}\int_{M} \int_{S^1 \times S^1 \times S^1}
\left(  \frac{f_3(p,\theta)}{\pi^{3}} +  \frac{k(p) f_0(p,\theta) }{4 \pi^3}
\right)
 d \vec{\theta} dA  + O(\ln(\lambda)\lambda^{-\frac{1}{2}}) .  \] 
\end{theorem}
{\bf Proof:}
As a preliminary observation, the above fact about $a(r,p)$  and the 
mean value theorem gives us
\[  e^{- \lambda a(p,\theta,r)} = e^{- \lambda\pi r^2}
\left(1 + \frac{k(z) \pi \lambda}{12}r^4 + O(\lambda r^5)\right) + 
O\left(\lambda r^8 e^{\lambda \pi r^2 (1 - \epsilon)}\right). \]

Using this fact
$ \Bbb E_{\lambda}(R)$ equals

 \begin{eqnarray}
  \frac{\lambda^3}{6} 
\int_{L_{\rho_M}} \left( f_0(p,\theta)
 + f_1(p,\theta) r +f_2(p,\theta) r^2 \ln(r) +  f_3(p,\theta) r^2 \right)   
  e^{- \pi r^2 \lambda} dr d \vec{\theta} 
dA ,
\\
 +  \frac{\lambda^3}{6} 
\int_{L_{\rho_M}}  r^3  \left( f_0(p,\theta)
 + f_1(p,\theta) r\right) \frac{\lambda \pi k(p) r^4}{12}   
  e^{- \pi r^2 \lambda} dr d \vec{\theta} 
dA,  
\\
+  \frac{\lambda^3}{6} 
\int_{L_{\rho_M}}  r^3  
\left(f_2(p,\theta) r^2 \ln(r) +  f_3(p,\theta) r^2  \right) 
\frac{\lambda \pi k(p) r^4}{12}   
  e^{- \pi r^2 \lambda} dr d \vec{\theta} 
dA  ,
\\+
 \frac{\lambda^3}{6} 
\int_{L_{\rho_M}}  r^3   
\frac{\lambda \pi k(p) r^4}{12}   
    f(r,p,\theta) O(\lambda r^5 e^{- \pi r^2 \lambda})
 dr d \vec{\theta} 
dA  ,
\\
+ \frac{\lambda^3}{6} 
\int_{L_{\rho_M}}  r^3   
\frac{\lambda \pi k(p) r^4}{12}   
   f(r,p,\theta) O(\lambda^2 r^8 e^{- \pi r^2 \lambda (1 - \epsilon)})
 dr d \vec{\theta} 
dA  ,
\\
+ \frac{\lambda^3}{6} 
\int_{L_{\rho_M}}  r^3   
\frac{\lambda \pi k(p) r^4}{12}   
   O(r^6 \ln(r)) O(\lambda^2 r^8 e^{- \pi r^2 \lambda (1 - \epsilon)})
 dr d \vec{\theta} 
dA 
\end{eqnarray}

The first thing to observe at this point is that for any $\delta$

 \[ \int_{0}^{\delta} e^{-\lambda \pi r^2} r^k dr =
\frac{m_k}{{(\pi \lambda)}^{\frac{k+1}{2}}}   + O(\lambda^{-\infty}) \]
and similarly with the $\ln(r) r^k $ integral.  This is because 

\[ \int_{0}^{\delta} e^{-\lambda \pi r^2} r^k dr =  \int_{0}^{\delta {(\pi \lambda)}^{\frac{1}{2}}} e^{-s^2} \frac{1}{ {(\pi \lambda)}^{\frac{k+1}{2}}}  s^{k }  ds \]
\[ = \int_{0}^{\infty} e^{-s^2}   \frac{1}{ {(\pi \lambda)}^{\frac{k+1}{2}}} s^{k }  ds  - \int_{\delta {(\pi \lambda)}^{\frac{1}{2}}}^{\infty} e^{-s^2}  \frac{1}{ {(\pi \lambda)}^{\frac{k+1}{2}}} s^{k}  ds \] 
\[ = \frac{m_k}{{(\pi \lambda)}^{\frac{k+1}{2}}} +  \frac{1}{ {(\pi \lambda)}^{\frac{k+1}{2}}} d^k\left(\delta {(\pi \lambda)}^{\frac{1}{2}}\right) \]

Where   $d^k(x) = \int_{x}^{\infty} e^{-s^2} s^k ds$ - which by 
l'hoptial's rule decays faster than any polynomial.

Armed with this observation we may explicitly integrate
 the first two terms in 
the above function with the only expense being a $O(\lambda^{-\infty})$ term.
 After doing so we  find the non $O\left(\ln(\lambda)
\lambda^{-\frac{1}{2}}\right)$ terms in 
$  \Bbb E_{\lambda}(R)$ to equal equal 
 \[  \frac{1}{6} 
\int_{M} \int_{S^1 \times S^1 \times S^1}  \left( f_0(p,\theta) 
\frac{\lambda m_3}{\pi^2}
 + f_1(p,\theta) \frac{\lambda^{\frac{1}{2}} m_4}{\pi^{\frac{5}{2}}} + 
\right.\]
\[ \left.
f_2(p,\theta) 
\left( \frac{\lambda^{\frac{1}{2}}  m_4^1}{\pi^{\frac{5}{2}}}  +
 \frac{\ln(\pi \lambda) \lambda^{\frac{1}{2}} m_4}{\pi^{\frac{5}{2}}}   
\right) +  f_3(p,\theta) \frac{m_5}{\pi^{3}}
  \right)   
d \vec{\theta} 
dA \]
\[+  \frac{1}{6} 
\int_{M} \int_{S^1 \times S^1 \times S^1} \frac{k(p)}{12}  
\left( f_0(p,\theta) \frac{m_7}{\pi^3}
 + f_1(p,\theta) \frac{\lambda^{\frac{1}{2}}m_8}{\pi^{\frac{7}{2}}} \right)
   d \vec{\theta} 
dA,\] 
where     
$m_k = \int_{0}^{\infty} r^k e^{-r^2} dr $.  Plugging in 
$m_3 = \frac{1}{2}$,$m_2 = \frac{3 \sqrt{\pi}}{8}$,
 $m_5 = 1$, $m_7 =  3$ and $m_8 = \frac{105 \sqrt{\pi}}{32}$ we 
arrive at the above formula.

  Likewise upon integration we see the remaining 
 four terms are indeed $O(\ln(\lambda) 
\lambda^{-\frac{1}{2}})$.

{\bf Q.E.D}

We now can compute this for the random 
variable $F$  counting the number of faces.

%[The Many Faces of Curvature Theorem]
\begin{corollary}\label{gauss}
\label{mfc}
With $\delta < \rho_M$ and calling   
$\int_{\times^3 S^1} \nu(\vec{\theta})  d \vec{\theta}  =  \nu$, 
we have the expected number of faces in the configuration can be expressed as 
\[ \Bbb E_{\lambda} (F)
% =  \frac{\lambda^3}{6} \int_{\Phi(L_{\delta})} e^{-\lambda a(y)} dA^3   
=  A \frac{\nu }{12 \pi^2}   \lambda 
-    \frac{\nu }{24 \pi^3 } \int_{M} k dA  + O\left(\ln(\lambda)
\lambda^{-\frac{1}{2}}\right) . \]
\end{corollary}
$\bf{Proof:}$
From the theorem in the previous section the $f(r,p,\theta)$ 
corresponding to the $F$ random variable is 
\[ j_{\theta_1}j_{\theta_2}j_{\theta_3} \nu(\vec{\theta}).\]
Using the o.d.e. describing the Jacobi field and the fact $M$ is compact 
we now immediately
 have $j_{\theta_0}$'s  Taylor expansion is   $j_{\theta_0}(r) =  
r(1  - r^2 \frac{k}{6} + O(r^3)) $, 
and as such  $j_{\theta_1} j_{\theta_2}j_{\theta_3} 
=  r^3(\pi - \frac{ k}{2} r^2 + o(r^3))$. 
So we have
\[ f(r,p,\vec{\theta}) =  \nu(\vec{\theta}) \left(1 - \frac{kr^2}{2} + O(r^3)\right) ,\]
as needed.
{\bf Q.E.D.}

%\[\Bbb E(F) = 
%A  \frac{ \nu }{12 \pi^2}   \lambda  +  
%\nu \int_{M} k dA \left(  \frac{1}{ 24 {\pi}^3} -  
%\frac{1  }{12 \pi^3 } \right)  + O(\lambda^{-\frac{1}{2}}). \]

%\qed
 
Let us now use the Euler-Delaunay-Poisson formula 
(formula \ref{edpf}) to figure out the constant $\nu$ 
( if you prefer, the integral it represents is easy to 
compute). We know the Euler characteristic  
of the flat tours of area one is zero; so as 
$\lambda \rightarrow \infty $ we have

\[ 0 = \Bbb E(\chi(M)) = \lambda - \frac{1}{2} \Bbb E (F) = (1-
 \frac{ \nu }{24 \pi^2})   \lambda , \]
forcing  $\nu  = 24 \pi^2 $.

Plugging this into the formula in the corollary  we  
have the needed equation (\ref{ex4}).

%$  \Bbb E_{\lambda}(R)$ equals 
% \[  \frac{1}{6} 
%\int_{M} \int_{S^1 \times S^1 \times S^1}  \left( f_0(p,\theta) 
%\frac{\lambda}{2 \pi^2}
% + f_1(p,\theta) \frac{\lambda^{\frac{1}{2}} 15}{16 \pi^2} + 
%\right.\]
%\[ \left.
%f_2(p,\theta) 
%\left( \frac{\lambda^{\frac{1}{2}}  c}{\pi^{\frac{5}{2}}}  +
% \frac{\ln(\pi \lambda) \lambda^{\frac{1}{2}} 15}{16 \pi^2 }   
%\right) +  f_3(p,\theta) \frac{1}{\pi^{3}}
%  \right)   
%d \vec{\theta} 
%dA \]
%\[+  \frac{1}{6} 
%\int_{M} \int_{S^1 \times S^1 \times S^1} \frac{k(p)}{12}  
%\left( f_0(p,\theta) \frac{3}{\pi^3}
% + f_1(p,\theta) \frac{\lambda^{\frac{1}{2}}105}{32 \pi^3} \right)  
% d \vec{\theta} 
%dA\] 
%\[  + O(\ln(\lambda)\lambda^{-\frac{1}{2}}) .  \] 

\end{subsection}

\begin{subsection}{Silly Continuity Results}\label{secsilly}

 We have several functions  which we need continuous,  and  several sets we need either open or  closed and  measure zero; here we deal with these issues.  
 Let us start by collecting and dealing with the needed closed measure zero results.

\begin{fact} \label{silly1}
The following sets are closed and measure zero  where indicated: 
\begin{enumerate}
\item
$\{ (x_1, \cdots , x_n) \in \times^n M : x_i = x_j with i \neq j \}$ in $\times^n M$.  
\item
$\Phi_f(L_{\delta})^c$ in  $V_{\delta}$. 
\item
The set of points $x$ in  $\times^n M$ where four points in $set(x)$ lie upon a circle of radius $r \leq min\{ \frac{i}{6}, \tau \} $.
\item
The set of points $x$ in  $\times^n M$  where three points in $set(x)$ lie upon a circle of radius exactly $\delta$.

\end{enumerate}
\end{fact}  
$\bf{Proof:}$ We will simply realize these sets as finite unions of the differentiable images of lower dimension sets, hence measure zero.  They will all be clearly relatively closed in the specified sets.

Let $D: M  \rightarrow M \times M$ be $D(x) =  (x,x)$.  The set we are removing is precisely the union of the  $n!$ index permutations of the set $D(M) \times \left( \times^{n-2} M  \right)$, as needed.

For the next part I'll prove something a bit stronger.  Let $\bar{V}_{min\{ \frac{i}{6}, \tau \}}$ be the set of all triples on circles of radius less than or equal to $min\{ \frac{i}{6}, \tau \}$ with repeats allowed.  I will show the stronger result that  $\Phi_f$ acting on $\times^3 S^1 \times \left[ 0, min\{ \frac{i}{6}, \tau \} \right] \times M - \{ points \}_{f} $ hits all but a measure zero set of $\bar{V}_{min\{ \frac{i}{6}, \tau \}}$.  To see this  let  $g$ be  a frame defined at each point of  $\{ points \}_{f}$; and then we have missed exactly 
\[  \bigcup_{p \in \{ points \}_f} \Phi_g (\times^3 S^1 \times \left[ 0, min\{ \frac{i}{6}, \tau \} \right]  \times M - \{ points \}_{g})(z = p) \]
This is a finite union of the differentiable images of  4 dimensional sets so measure zero and clearly relatively closed in the 6 dimensional space. 

For the third part
first look at the image of the the 7 dimension space 
\[ \times^3 S^1 \times \left[ 0, \frac{\tau}{3} \right] \times \left( M - \{ ponts \}_f \right) \times S^1 \]
under the mapping given by
\[ (\Phi_f(\vec{\theta},r,z) ,exp_z(r v(\theta_4)). \] 
It is measure zero as a subset of $ \times^4 M$. Union it with the same set using the above $g$.  With this unioning it is clearly relatively closed in $\times^4 M$.  Now product this with the other $n - 4$ components and take the union over the $n!$ permuted copies of this set in  $\times^{n} M$ and remove them; hence removing all the possible quadruples on the same circles.  This amounts to a  finite removal of closed measure zero sets, as needed. 

Similarly for the last part, where we use instead the the image under $\Phi_f$ (and $\Phi_g$)  of $\times^3 S^1 \times \delta \times M$ in $\times^3 M$.

\qed

Now we need some functions to be continuous and some sets to be open.

\begin{fact} \label{silly2}
$ T_{\delta}$ is open in $\frak{P}$; and   $F$ , $f_s$,  $E$ and $V$ are continuous in $\frak P$.
\end{fact}

$\bf{Proof:}$
 All these cases could we handled simultaneously if one could perturb a set of points and not change any of the assignments given to a triple,  a pair, or a singleton.  We may restrict or attention to a set coming form $\times^n M_-$, since $\frak P$ is a disjoint union of such sets. Clearly the assignment is discontinuous on $\times^n M$; and what is being claimed is that all the discontinuities occur on the removed measure zero sets.

That any singleton  is still a distinct vertex is true due to the removal of the set corresponding to part one of the above fact. 
In fact this removal guarantees all points in $set(x)$ can simultaneously  be separated by open sets  (including any pair or triple). 

To deal with the stability of face assignment to a triple, $set \cdot \pi_s(x)$, first observe the removal of the fourth set  guarantees $V_{\delta}$ and $V_{\delta}^c$ are open.  So $set \cdot \pi_s(x)$'s relationship to $\delta$ is open, finishing off the stability argument if the point fails to be on a circle of radius less than $\delta$.  In the case where the triple does lie on such  a circle, note that $\pi_s(x) \in \Phi_f(L_{\delta})$ for a suitable $f$.  So the disk associated to a triple varies continuously with the triples position ($r$ and $z$ in the parameterization).   Now suppose $set \cdot \pi_s (x)$ is not assigned a face due to a point in its associated disk's interior.  This violating point must have a neighborhood in the disks interior.  So by the above continuity of disk position, there is an open set $U_1 \times U_2 \times U_3$ about $ \pi_s(x)$ and  an open set $U_p$ about $p$ such that $U_p$ is  contained in the disks associated to $set(y)$ for each $y \in U_1 \times U_2 \times U_3$.  So the triple having no associated face is indeed an open condition. Virtually this same argument guarantees that in the case when the disk is empty that there is an open set about $set(x) - set \cdot \pi_s(x)$ such that no point in this open set is in the disk associated to $set(y)$ for $y$ in some  $ U_1 \times U_2 \times U_3$.  The only difference is that we must note that we removed the possibility of four points  on a circle (part 3 above) - so indeed each point in $set(x) - set \cdot \pi_s(x)$ is in an open set separating it form the closed disk, as needed. So off the above measure zero sets the notion of face is indeed stable in an open set.

The edges require a small amount more thought.
An edge $pq$ 
existing  implies there is a $k$ such that $D_{pq}(k)$ is empty.   
We may attempt to deform the inflating family
to the left and to the right of $k$.  
If we cannot it is because a third point lies to the left or right side 
of the circle.   Now since we may assume there is not a fourth point on 
the same circle, in one of these directions so we indeed can deform 
our inflating family; and in fact in both directions if the circle's 
radius is $\delta$ - since no triple lives on such a circle.  
In particular, choosing a different $k$ if necessary, we now 
have a circle of radius less than $\delta$ which is empty, and is 
contained in an neighborhood empty of other vertices.  Now one can 
proceed exactly as above to note that $E$ is constant on a neighborhood.    
The notion of no edge requires a more delicate use of our inflating family
ideas.  
The key observation is that from the monotonicity lemma a point is in 
$intD_{pq}(h) \bigcap intD_{pq}(k) $ if and only if it is in each 
$intD_{pq}(c)$ for $c \in [h,k]$. From this the idea is to find 
$a = c_0 < c_1 < \dots < c_n = b$ such that 
$D_{pq}(c_i - \epsilon) \bigcap D_{pq} (c_{i+1} +\epsilon)$ 
contains a point in its interior.  Then each of these sets will, 
as above, satisfy this property when the points are 
perturbed - and the notion of no edge will be stable.  
(Note at $a$ and $b$ we need not use the 
$\epsilon$ - since there is a point in the 
interior of each of these disks which we may use.)  
To construct the $c_i$ start at $a$ and take a point 
inside it.  Move the inflating family
rightward until at some $d_i$ this point fails to be in the disk 
(if such a point does not exist we may use $c_0 = a $ and $c_1 = b$.)  
Now since there is no empty disk $intD_{pq}(d_i)$ contains some other point.
Being in its interior this is in fact true for the parameters in 
$(d_i - \eta, d_i)$. Let $c_i = d_i - \frac{\eta}{2}$ and 
$\epsilon = \frac{\eta}{4}$.  Now continue this process
 making $\epsilon$ smaller if necessary.  There are 
a finite number of distributed points - so eventually one 
must hit $b$ or an empty disk.    An empty disk is impossible, 
since no edge was put in; and we are done.

%\begin{figure*}
%\vspace{.01in}
%\hspace*{\fill}
%\epsfysize = 1in 
%\epsfbox{zero.eps}
%\hspace*{\fill}
%\vspace{.01in}
%\caption{\label{zero} The Open Sets }
%\end{figure*}

%The last thing to deal with is the openness of $\frak D_{\delta}$ . Suppose not and you could not find an open set about $set(x) = \{p_1 , \dots ,p_n \} \in \frak D_{\delta}$ such that the open set remained in $1_{\frak D_{\delta}}$.  Then  there must be a sequence $p_i^j \rightarrow p_i$
%where this violation occurs, hence a sequence of radius $\delta$ violation disks $D^j$.
%Such a disk is completely determined by its center $c^j$ and $M$ is compact so there is subsequence of which converges to say $c$.  The disk of radius $\delta$ at $c$, being a disk of radius $\delta$, must contain a point $p$ of $set(x)$.  Since it is in the disks interior - there is a neighbor hood of the center $c$ such that all the disks in this neighborhood contain $p$ - contradicting the fact $c^j$ converged to $c$. So indeed the set is open. 

\qed
\end{subsection}
\end{section}

\end{chapter}

\begin{chapter}{The Continuous Uniformization Theorem}\label{contunith}
This section accomplishes two things.  For starters in section 
\ref{randomenergies}
 we use theorem \ref{exp} of the 
previous chapter to actually calculate the energies on the 
space of metrics of interest form section 1.2.2 and 1.2.3. 
This involves first a careful  look at  how the angles in a geodesic 
net deform under a conformal breeze followed by 
the actual computations.
Section 4.2 contains an actual proof of the uniformization 
theorem for surfaces with 
$\chi(M) < 0$,  mimicking in a precise sense the 
discrete uniformization proof.   

\begin{section}{The Random Energies} \label{randomenergies}
\begin{subsection}{Triangulation Deformations} \label{triangulardef}
The key to computing the need expected values in section \ref{conu} 
is to compute
how an angle in a triangle deforms under a conformal change of metric.  
The method used here is to solve the 
 the boundary value
 problem 
 for the geodesics in the $e^{2 \phi}g$ metric, and compare the initial 
directions.   Of real interest is the case when we have 
 the geodesics forming the edge of a triangle,
 so I will phrase the results in this language.  
 Using normal coordinates at the point $p$ in the $g$ metric let
  the points on our triangle are labeled 
 $r v(\theta_i) = r v^i 
 =  ( r \cos(\theta_i), r  \sin(\theta_i))$.  Further let 
$v_{ij} = v(\theta_j) - v(\theta_i)$ for any $i$ and $j$
 and let $v_{ij}^{\perp}$ be the left handed 
  $\frac{\pi}{2}$ rotation  of $v_{ij}$ in the Euclidean metric.  Given that the surface is orientable pick an orientation and let $\sigma(i,j,k)$ be one if $(v_i,v_j,v_k)$ is ordered in the "clockwise" direction and $-1$ if not.  Further more let $k$ be the Gaussian curvature at $p$ and ${\bf H}$ be Hessian 
of $\phi$ in normal coordinates (which are denoted via $(x,y)$).

\begin{lemma} \label{wind}
In the above notation the initial direction of a  geodesic from $v^i$ to $v^j$ in the $e^{2\phi}g$ metric is given by

\[ \left( \frac{r}{2}( v_{ij} \cdot \nabla \phi) + O(r^2) \right) v_{ij}  \]
\[ -   \left( \frac{r}{2}( v^{\perp}_{ij} \cdot \nabla \phi) 
+ r^2 \left( \frac{1}{2}
(v^{\perp}_{ij})^{tr} {\bf H} \left(v^i + \frac{1}{3}v_{ij} \right) -
 \frac{k}{3}(v^{\perp}_{ij})^{tr} v^i \right. \right. \]
\[ \left. \left. +
\frac{1}{3} ( \nabla \phi \cdot (v^{\perp}_{ij})^{tr}  )(\nabla \phi \cdot v_{ij}) \right) + O (r^3)\right) v^{\perp}_{ij} . \]
  
\end{lemma}
     
   {\bf Proof:}
 For the 
   purposes of readable notation 
   let $v_{ij} = {\bf v} =  (v,w)$, and note we may  
rewrite our initial starting direction  as ${\bf v +  l}$.
% and further as 
% $ {\bf v +  l} = (1+m) {\bf v} + n \sigma {\bf v}^{\perp} $.

  The new metric in normal coordinate is up to order $O(r^2)$ given by 
  \[ g_{mn} = e^{2 \phi}  \left[ \begin{array}{ll} 1 - \frac{1}{3} k y^2 
  & \frac{1}{3} k xy \\ \frac{1}{3} k xy 
  & 1 - \frac{1}{3} k x^2  \end{array} \right] .\] 
  Recall the geodesic equation is 
  \[ \frac{d^2 x_k}{dt^2} = - \Gamma_{mn}^k \dot{x}_m^2 \dot{x}_n^2 . \]
  
  Using this a metric up to order $O(r)$ 
 the Christoffel symbols are found to be
  
  \[ \Gamma_{mn}^1 = \left[ \begin{array}{ll}
   \phi_{x}  + \phi_{xx}x + \phi_{xy}y
   & \phi_{y}  + \phi_{yy}y + \phi_{xy}x  - \frac{1}{3}ky    \\
  \phi_{y}  + \phi_{yy}y + \phi_{xy}x - \frac{1}{3}ky   & 
  -(\phi_{x}  + \phi_{xx}x + \phi_{xy}y) + \frac{2}{3}kx
  \end{array} \right] \] 
  and
  \[ \Gamma_{mn}^2 =  \left[ \begin{array}{ll}
  -(\phi_{y}  + \phi_{yy}y + \phi_{xy}x) + \frac{2}{3}ky
   & \phi_{x}  + \phi_{xx}x + \phi_{xy}y - \frac{1}{3}kx    \\
  \phi_{x}  + \phi_{xx}x + \phi_{xy} y-  \frac{1}{3}kx  & 
  \phi_{y}  + \phi_{yy}y + \phi_{xy}x
  \end{array} \right] .\]

 All the above $\phi$ derivatives
 are evaluated at the center of the circle, $p$; 
and I am letting $k$ denote  the curvature a $p$.
   We are solving the boundary value problem where we start at 
$rv^i$ and ending a 
  $rv^j$, i.e. introducing the notation that 
 $(p,q) = (\dot{x}, \dot{y}) -\bf  v$ we have 
\[ \dot{{\bf x}} = \left[  \begin{array}{l} \dot{x}\\ \dot{y} \\ \dot{p}
   \\ \dot{q}
  \end{array} \right] = \left[  \begin{array}{l} p+v\\ q+w \\ -\left(
\Gamma_{11}^1 (p+v)^2 + 2 \Gamma_{12}^1 (p+v)(q+w) +\Gamma_{22}^1 (q+w)^2 
\right)
   \\ -\left( \Gamma_{11}^1 (p+v)^2 + 
2 \Gamma_{12}^1 (p+v)(q+w) +\Gamma_{22}^1
 (q+w)^2\right)  \end{array} \right].\] 
 
  Since we know the Christoffel symbols up to $O(r)$, 
 to analyze the solution up to order $r^3$
 we may
   simply linearize the problem form the point of view of the origin 
(this can be seen immediately by looking at the power series expansion).
     So we need to solve the linear o.d.e.
 \[ \dot{{\bf x}} = \left[  \begin{array}{l} \dot{x}\\ \dot{y} \\ \dot{p}
   \\ \dot{q}
  \end{array} \right] ={\bf b} + A x ={\bf b} + 
  \left[  \begin{array}{ll} {\bf 0} & I  
   \\ -C &-2 D  
  \end{array} \right] x ,\]
  
  where {\bf b}
  \[ {\bf b}  = \left[  \begin{array}{l} {\bf v} \\
   {\bf f}  
  \end{array} \right]  = \left[  \begin{array}{l} v \\ w \\
  -( \phi_x(v^2-w^2) + 2 \phi_y vw) 
   \\  -( 2\phi_x vw+  \phi_y ( w^2-v^2))  
  \end{array} \right] = \left[  \begin{array}{l} {\bf v} \\
   {\bf v}^{\perp}( {\bf v}^{\perp} \cdot \nabla \phi) -
{\bf v}( {\bf v}\cdot \nabla \phi) 
  \end{array} \right] \]
  
  and 
  \[ C = \left[  \begin{array}{ll}
  \phi_{xx}(v^2 - w^2) + 2 \phi_{xy} vw + \frac{2}{3}kw^2
     & \phi_{xy}(v^2 - w^2) + 2 \phi_{yy} vw  - \frac{2}{3}kvw  \\ 
     -\phi_{xy}(v^2 - w^2) + 2 \phi_{xx} vw - \frac{2}{3}kvw
     & -\phi_{yy}(v^2 - w^2) + 2 \phi_{xy} vw  + \frac{2}{3}kv^2
  \end{array} \right] \]
\[ =- {\bf v}^{\perp}{( {\bf v}^{\perp})}^{tr} {\bf H}(p) +
 {\bf v} {\bf v}^{tr} {\bf H} (p)  + \frac{2k}{3}  {\bf v}^{\perp}
{( {\bf v}^{\perp})}^{tr} \] 
  
  \[ D =   \left[  \begin{array}{ll}
    \phi_x v +  \phi_y w 
     & -\phi_x w +  \phi_y v   \\ 
    \phi_x w -  \phi_y v 
     & \phi_x v +  \phi_y w 
  \end{array} \right] = {\bf v} {(\nabla \phi)}^{tr} + 
 {\bf v}^{\perp}{(\nabla \phi^{\perp})}^{tr} \]

For future use it is convenient to note 
\[  D {\bf f} = {\bf v}\left( (\nabla \phi \cdot {\bf v}^{\perp})^2 
- (\nabla \phi \cdot {\bf v})^2 \right) +
 2  {\bf v}^{\perp}\left((\nabla \phi 
\cdot {\bf v}^{\perp})(\nabla \phi \cdot {\bf v})\right). \]

  Now this o.d.e can be solved via variation of parameter via
  
  \[ {\bf x} = e^{tA} \left( \int_{0}^{t} e^{-tA} {\bf b}
   dt \right) + e^{ta} {\bf x_0} \] 
   
   with initial condition
   \[ {\bf x_0} = \left[  \begin{array}{l} r v^i \\
   {\bf l}   
  \end{array} \right]= \left[  \begin{array}{l} r \cos(\theta_1) \\ 
  r \sin(\theta_i) \\
   l 
   \\  k  
  \end{array} \right] .\]
  
  So up to order $o(r^3)$ 
%(assuming  ${\bf x_0}$ is $o(r)$) 
we have 
  
  \[ {\bf x} = (t + \frac{t^2}{2} A + \frac{t^3}{6}A^2) {\bf b} +
  ( I + tA + \frac{t^2}{2}A^2) {\bf x_0}. \]
  
  Now we need to find the $l$ and $k$ so that when this 
is evaluated a $r$ it is
  at the point $r v^j$.
  
  So plugging in $r$ we find we need
  \[ (I - r D) {\bf l} = \frac{r}{2} \left(
   \left(C r v^i + \frac{r}{3}C{\bf v}\right)  - \left(I-  \frac{2r}{3}D)
{\bf f}\right)
 \right). \]
   
   Inverting to isolate ${\bf l}$ we find (up to order $O(r^2)$) that
   \[ {\bf l} = \frac{r}{2}
   \left( 
 rC(  v^i + \frac{r}{3} {\bf v})- (I +
 \frac{r}{3}D) {\bf f} \right) = - \frac{r}{2}f + \frac{r^2}{2}  \left( 
 C(  v^i + \frac{1}{3} {\bf v})-  
 \frac{1}{3}D {\bf f} \right). \]
   
 Simply plug into this expression 
the above formulas to arrive at the claimed formula.

{\bf Q.E.D.}

Now we are capable of computing the needed angle.  
It is convenient to denote the angle in the triangle at the point 
$r v^i$ in  the Euclidean  coordinates as ${\Bbb E}_i$.  With this we have...
  
 \begin{lemma}\label{def}
The angle in the triangle at $r v^i$ is
 
\[\psi^i_{\phi} = {\Bbb E}_i +
 \sigma(i,j,k) \left(- \frac{r}{2}({ v_{kj}}^{\perp} \cdot \nabla \phi)
 + \frac{r^2}{2} b_i \right) + O(r^3) , \]
where
\[ b_i =    \frac{k}{2}( { v}_{kj} \cdot(v^i)^{\perp})
- ({ v}_{kj}^{\perp})^{tr} {\bf H}  v^i 
   - \frac{1}{3} \left( ( { v}_{ji}^{\perp})^{tr} {\bf H}{ v}_{ji} 
-  {({ v_{ki}})^{\perp}}^{tr} {\bf H}v_{ki} \right) \]
\[  -
\frac{1}{6}\left(((\nabla \phi
 \cdot { v_{ji}}^{\perp} )(\nabla \phi \cdot { v_{ji}}))
 - ((\nabla \phi
 \cdot { v_{ki}}^{\perp} )(\nabla \phi \cdot { v_{ki}})) \right)  . \]

  \end{lemma}

{\bf Proof:}
For this problem denote the solution to the previous problem as 
${\bf l + v} = (1+m_{ij}) {\bf v} + n_{ij}{\bf v}^{\perp} $.  
First note that the angle, $\eta_{ij}$,  between this
 initial vector and ${\bf v}$   has its cosine
    given by

\[ \cos(\eta_{ij}) = \frac{<({\bf l} + {\bf v}),  {\bf v}>_g}
{(||{\bf v}||_g)(||{\bf l+v}||_g)}.\]

Hence using the notation of the previous problem and letting
 $c_1 = \frac{||{\bf v}||^2_g}{||{\bf v}||^2_{\Bbb E}}$,
 $c_2 =  \frac{||{\bf v}^{\perp}||^2_g}{||{\bf v}||^2_{\Bbb E}}$,
 and
 $c_3 =  \frac{<{\bf v}, {\bf v}^{\perp}>_g}{||{\bf v}||^2_{\Bbb E}} $,
 we have up to cubic order that   

\[ \cos(\eta_{ij}) 
 = \frac{(1+m_{ij} + c_1 + m_{ij} c_1 + n_{ij}c_3 )}{\sqrt{(1 + c_1)(1 
+ 2m_{ij} +m_{ij}^2 + n_{ij}^2
c_1 + 2m_{ij} c_1 + 2n_{ij} c_3) }} . \]

Multiplying out and using $\sqrt{1+x} =
 1 + \frac{1}{_2} x - \frac{1}{8} x^2 + \frac{1}{16} x^3 + \dots$ and
$\frac{1}{1 + x} = 1 -x + x^2 - x^3 \dots $ up to the third order in $r$
we have the nice fact that this expression (up to $r^3$) is  
independent of the $g_0$ metric and equal to...

 \[ \cos(\eta_{ij}) =
 1 + (- \frac{n_{ij}^2}{2} + m_{ij}n_{ij}^2).\]
% = 1 - \frac{r^2}{2}(f_2^2) + r^3(f_2^2 f_1 - f_2g_2.) \]

   Using cosine's  power series 
this allows us to isolate the angle up to second order as...
   \[ |\eta_{ij}| = |n_{ij}| \sqrt{1 - 2m_{ij}} = |n_{ij}| (1-m_{ij}),\]
Further note that if one would like this angle to 
positive if it contributes to the triangle's internal 
angle and negative if not, then 
\[ \eta_{ji} = \sigma(i,j,k)  n_{ij}(1-m_{ij}). \]

With this computation out of the way  we are left needing 
 to measure the angle between $v_{ij}$ and $v_{ik}$ in the $g$  metric.
 The trick will be to measure the angle $a_{ij}$ from $-v^i$ to 
$v_{ij}$ at $r v^i$  and
 the angle $a_{ki}$ from
 $-v^i$ to $v_{ki}$ at $r v^i$, both with 
 the proper signs as contributers  to $a^i_{\phi}$,  
  and then sum them up.
 Let $a_{ij}^{\Bbb E}$ be the angle in Euclidean
coordinates and $(\cdot,\cdot)$ be the Euclidean innerproduct.

To compute these angle note at $r v^i$ we have 
$||v^i||^2_g = 1$, $||(v^i)^{\perp}||^2_g= 1  -
 \frac{kr^2}{3} + o(r^3)$, $||v^i - v^j||_{\Bbb E
} = \sqrt{2} \sqrt{1 - (v^i, v^j).}$ and that

\[\cos( a_{ij}) = 
 \frac{<v^i,v^i - v^j >_g}{(||v^i||_g)(||v^i - v^j||_g)}\]
\[= \frac{1- (v^i,v^j)}{\sqrt{(1-(v^i,v^j))^2  + ((v^i)^{\perp},v^j)^2(1
  -\frac{kr^2}{3})}}  = \frac{1- (v^i,v^j)}{\sqrt{2}
 \sqrt{1- (v^i,v^j)} \sqrt{1  
  -\frac{kr^2((v^i)^{\perp},v^j)^2}{6(1- (v^i,v^j))}}}   \]
\[ =  \frac{\sqrt{1- (v^i,v^j)}}{\sqrt{2}}
 \left(1 + \frac{kr^2}{12}\frac{((v^i)^{\perp},v^j)^2}{ (1- (v^i,v^j))}
 \right) \]
\[  =
   \frac{\sqrt{1- (v^i,v^j)}}{\sqrt{2}}  - \frac{\sqrt{1+ (v^i,v^j)}}{\sqrt{2}} \left(\frac{-kr^2}{12} 
\frac{((v^i)^{\perp}, v^j)}{|((v^i)^{\perp}, v^j)|}
  ((v^i)^{\perp}, v^j) \right) . \]

Now note that
\[ \cos(a^{\Bbb E}_{ij}) =   \frac{\sqrt{1- (v^i,v^j)}}{\sqrt{2}} \]
and from this by keeping track of the  necessary sign we have 
\[ \sin(|a^{\Bbb E}_{ij}|) = \sigma(i,j,k)
 \frac{((v^i)^{\perp}, v^j)}{|((v^i)^{\perp}, v^j)|} 
\frac{\sqrt{1+ (v^i,v^j)}}{\sqrt{2}} .\]

So the above expression is precisely 

\[  
\cos(a^{\Bbb E}_{ij}) -  \sin(|a^{\Bbb E}_{ij}|)\left(\frac{-\sigma(i,j,k)
kr^2}{12}
  ((v^i)^{\perp}, v^j) \right)=  \cos\left (a^{\Bbb E}_{ij}   - \frac{\sigma(i,j,k)k r^2}{12}((v^i)^{\perp}, v^j)\right). \]

So we have that $a_{ij} =a^{\Bbb E}_{ij}  
  -  \frac{\sigma(i,j,k) k r^2}{12}((v^i)^{\perp}, v^j) $, and from this the needed    
\[ a_{\phi}^i = \Bbb E_i + \sigma(i,j,k) \frac{k r^2}{12}((v^i)^{\perp}, v^j)
+   \sigma(i,k,j) \frac{k r^2}{12}((v^i)^{\perp}, v^k) \]
\[  = \Bbb E_i + \sigma(i,j,k) \frac{km r^2}{12} (v_{jk}^{\perp}, v^i)  \]

Now summing up to get  $\psi^i_{\phi} = \eta_{ji} + \eta_{ki} + a_{\phi}^i$,  
and so 
\[ \psi_{\phi}^i =  \Bbb E_i + \sigma(i,j,k) \frac{kr^2}{12} (v_{jk}^{\perp}, v^i)  + \sigma(i,j,k)  n_{ij}(1-m_{ij}) +  \sigma(i,k,j)  
n_{ik}(1-m_{ik})\]
\[ =  \Bbb E_i + \sigma(i,j,k)\left ( \frac{kr^2}{12} (v_{jk}^{\perp}, v^i)
  +   n_{ij}(1-m_{ij}) -  
n_{ik}(1-m_{ik}) \right) .\] 
Plugging into the formula from the previous lemma 
now finishes the computation.

{\bf Q.E.D.}

\end{subsection}
\begin{subsection}{The Energy Computation}\label{encomp}
Using the angle formula from the previous section I will now compute the 
the expected value of the energy, deriving formula
\ref{cen} and theorem \ref{lapy}. 
In the end we will arrive at the formula 
  
\begin{formula} \label{formu}
For a negative curvature metric $g$ and  $h = e^{2 \phi} g$  we have 
\[ \Bbb E^{g}_{\lambda} \left( E_h \right) = 
  D_0 A  \lambda + D_1 \chi(M) \ln(\lambda) \lambda^{\frac{1}{2}}+ D_2 \chi(M) \lambda^{\frac{1}{2}}+ 
D_3 \chi(M) \]
\[  +\int_M || \nabla \phi ||^2 +(\Delta \phi - k)
 \ln(\Delta \phi - k)  dA + o(\ln(\lambda) \lambda^{-\frac{1}{2}}),\]
 with the $D_i$ constants.
\end{formula}

{\bf Proof:}
To do the computation it is necessary to find the 
$f(r,p,\theta)$ function used in  theorem \ref{exp}  
with respect 
to the volume induced random variable.  

By changing to the $\phi$ 
 coordinates  we reduce   $ \Bbb E^{g}_{\lambda} \left( E_h \right)$ to 
\[ 
  \frac{\lambda^3}{6}\int_M \int_{[0,2\pi]^3}
 \int_{0}^{\infty} V(p,r,\vec{\theta}) 
 e^{-\lambda r^2} \nu(\vec{\theta}) r^3 (1 - r^2 \frac{k}{2} + o(r^3))
  dr d\vec{\theta} dA, \] 
where $V(p,v, \vec{\theta})$ is 
the prism volume associated to the triangle 
withe angle data given by the triangle on this surface formed with this data.
So we need to expand 

\[  V(p,r,\vec{\theta})  \nu(\vec{\theta}) 
r^3 (1 - r^2 \frac{k}{2} + o(r^3)) \]
in the $r$ variable.
To do this it is nice to give our small objects some  names; let
 $A = \pi - \sum \psi^i_{\phi}$, let $\delta_i = \psi_i - \Bbb E_i$,and let
$\hat{\delta_i} =  \psi_{\phi}^{i} + \frac{A}{2} - \Bbb E_i$.  The last of which is small since
\[ \pi + \Bbb E_i = \frac{\pi + \Bbb E_i - \Bbb E_j - \Bbb E_k}{2}. \]
Note that at this point 
\[ V(r,p,\theta) = \sum_{i= 1}^3 \Lambda(\Bbb E_i + \delta_i) + 
\Lambda(\Bbb E_i + \hat{\delta}_i)  + \Lambda \left( \frac{A}{2} \right)\]

The power series expansion of minus the Lobacevskii function
 $\Lambda$ about a positive $E$ is 
\[  \Lambda(E + \delta) = \Lambda(E) -
 \ln(2|\sin(E)|) \delta - \cot(E) \delta^2  +
O\left(\frac{\delta^3}{(\sin(E))^2}\right).\]
Note from the formula for $\psi_i$ in lemma 
\ref{def}
 that the  $\delta$ and $\hat{\delta}_i$ functions are divisible by 
$\sin(E_i) =  \sin(\frac{\theta_j - \theta_k}{2})$ so this 
series when applied to the first six terms in the above expansion 
for $V(r,p,\theta)$ gives us a power series in $r$ with 
bounded continuous coefficients and a remainder of order $O(r^3)$. 

%(check!)

%Now note that for $r < \rho_M$ that $\frac{A}{2}$ we have that 
%$\Lambda$ is singular so
% this term in the volume function is singular, yet still easy to
% express up to order $r^2$ via observing that for small $x$

Note that the remaining term is in the form 
\[  \Lambda\left(\frac{A}{2}\right) = 
 - \int_0^{\frac{A}{2}} \ln(2x) dx - \int_{0}^{\frac{A}{2}} 
\ln \left( \frac{\sin(x)}{x}\right) \]
\[ =  \frac{A}{2}  - \frac{A}{2} \ln(A) + O(A^2).\]

Given  this expression it is extremely useful to get a grip on 
the expression $\pi - \sum_{i=1 }^{3} \psi_{\phi}^{i}$.

\begin{lemma} \label{area}
 \[ A = \sum_{ \in t} \psi_{\phi}^i - \pi 
= (- \Delta \phi(p) + k) r^2 \nu(\vec{\theta}) + O(r^3). \]
\end{lemma}

{\bf Proof:}

Let $i<j$ represent a  pair where $\sigma(i,j,k)=1$.

 We can begin with the observations that $\sum_{i<j}v_{ij} =0$, and that 
since both $i<j$ and $j>i$  include all ordered pairs
\[ \sum_{i<j}{ v}_{ij}^{tr} {\bf H} {v}_{ij}^{\perp} 
 - \sum_{i>j}{v}_{ji}^{tr} {\bf H} {v}_{ji}^{\perp}  =0 \]

\[ \sum_{i<j} ((\nabla \phi
 \cdot {\bf v_{ji}}^{\perp} )(\nabla \phi \cdot {\bf v_{ji}}))
 -  \sum_{i>j}((\nabla \phi
 \cdot {\bf v_{ki}}^{\perp} )(\nabla \phi \cdot {\bf v_{ki}}))  = 0\]

and

\[ \sum_{i<j}{v}_{i}^{tr} {\bf H} { \bf v}_{ij}^{\perp} 
 - \sum_{i>j}{v}_{i}^{tr} {\bf H} {\bf v}_{ij}^{\perp}  =
\sum_{i<j}
\left({v}_{i}^{tr} {\bf H} { v}_{j}^{\perp} 
- {v}_{j}^{tr} {\bf H} { v}_{i}^{\perp} \right) \]
\[ + \sum_{i}
\left({v}_{i}^{tr} {\bf H} { v}_{i}^{\perp} -
 {v}_{i}^{tr} {\bf H} { v}_{i}^{\perp} \right) = \sum_{i<j}
\left({v}_{i}^{tr} {\bf H} { v}_{j}^{\perp} -
 {v}_{j}^{tr} {\bf H} { v}_{i}^{\perp} \right) 
.\]

   To compute explicitly the remaining
 terms the following fact  useful:

\begin{fact} \label{areahes}
If $A$ is a symmetric $2\times 2$ matrix then
\[ \sum_{i<j} (v^j)^{tr} A v^i = -\sum_{i<j} (v^i)^{tr} A v^j =
   tr(A)r^2  \nu(\vec{\theta}).\]
%where \[r^2  A_{\Bbb E} = r^2 \sum_{i<j} \sin(\theta_i - \theta_j) \]
%and is the area of the triangle in normal coordinates.

\end{fact}  

%{\bf Proof:}

%{\bf Q.E.D. (sub-lemma)}

Using these observations, the above fact, the fact that the angle in the 
Euclidean triangle sum to zero, and the fact that  
$tr ({\bf H}) =\Delta \phi$, we find the sum is 
\[ \left( - 2\frac{tr({\bf H})}{2} + 4\frac{k}{4}\right)
 r^2 \nu(\vec{\theta}) + O(r^3), \]
as needed.

{\bf Q.E.D. (lemma)}

It is worth noting as a confirmation to the previous sections computation, 
this is lemma is exactly what one expects from
 the the Gauss Bonnet formula when $\phi = 0$.

With this in mind we see that in fact the expansion for 
$- \Lambda \left( \frac{A}{2} \right)$ has terms of the form $f r^i$ for $i >2$
and $g \ln(r) r^i$ for $i>2$ with continuous coefficient functions.  
So  up to order $\ln(r) r^3$  $V(r,p,\theta)$  is of the form 
%satisfies the requirements of theorem \label{exp} and

\[V(r,p,\theta)
= \sum_{i= 1}^3 \left(  \Lambda(\Bbb E_i) -
 \ln(2|\sin(\Bbb E_i)|) (\delta_i + \hat{\delta}_i) - 
\cot(\Bbb E_i) (\delta_i + \hat{\delta}_i)^2 \right) \]
\[ + 
\frac{A}{2} - \frac{A}{2}  + O(\ln(r)r^3)\]

With this we are in position to apply theorem \ref{exp}.
Note that $\delta_i + \hat{\delta_i} = 
2 \psi_{\phi}^{i} + \frac{A}{2} -2 \Bbb E_i$.  
So using the notation of the previous section we have the 
$f_i(p,\theta)$ from theorem and we get 

\begin{eqnarray}
 f_0(p, \vec{\theta}) & = & \sum_{i= 1}^3  \Lambda(\Bbb E_i) 
\nu(\vec{\theta})  \label{zero} \\
 f_1(p, \vec{\theta}) & = &  -\sum_{i} \ln(2|\sin(\Bbb E_i)|) 
(v(\theta_k) - v(\theta_j))
  \cdot \nabla \phi  
       \nu(\theta) \label{one} \\
 f_2(p, \vec{\theta}) & =  &   ( \Delta \phi -  k)(\nu(\vec{\theta}))^2   \label{two}\\
 f_3(p, \vec{\theta}) & = & 
  \frac{k}{2} \Lambda(\Bbb E_i) \nu(\vec{\theta}) + 
\frac{ \Delta \phi -  k}{2}  \ln(2|\sin(\Bbb E_i)|)(\nu(\vec{\theta}))^2 
  \label{threea} \\
&  & -
 2 b_i  \ln(2|\sin(\Bbb E_i)|) \nu(\vec{\theta})  \label{threeb} \\
&  &- \frac{1}{2} ({\bf v}^{\perp}_{jk} \cdot \nabla \phi)^2 
\cot(\Bbb E_i)\nu(\vec{\theta})  \label{threec}  \\
&  &+  \frac{1}{2} ( \Delta \phi -  k) \ln(\Delta \phi -  k)
(\nu(\vec{\theta}))^2   \label{threed}  
 \end{eqnarray}

Now we attempt to compute 
$ I_i = \int_M \int_{[0,2\pi]^3} f_i(p,\theta) d\vec{\theta} dA $.
%The largest term is the $r^3$ given by
%\[ r^3 \sum_{i} \Lambda(\Bbb E_i) \nu(\theta), \]
The first one clearly can be integrated to a constant times the surfaces 
area, $I_0 = C_0 A $, where 
the constant is  independent of any of the geometry or topology. 
%(the finiteness follows form the integrability of the logarithmatic singular
%y). 
%The second largest is $r^4$ which comes
% form lemma ? in the previous section and is given by 
%\[ r^4 \sum_{i} \ln(2|\sin(\Bbb E_i)|) 
%( {\bf v}^{\perp}_{jk} \cdot \nabla \phi )  \nu(\theta). \]

Happily enough  $I_1 = 0$. 
This follows
  immediately form the first part of the following integral  vanishing 
lemma ( after noting 
$\sin(\Bbb E_i) =\sin(\frac{\theta_k - \theta_j}{2})$):

\begin{lemma}\label{van}
Let $f$ and $g$ be either $\cos$ or $\sin$ functions then
\[   \int_{[0,2\pi]^3} \ln\left(2|\sin
\left(\frac{\theta_k - \theta_j}{2}\right)|
\right)
f(\theta_k)  \nu(\vec{\theta}) d\vec{\theta} = 0 \]
\[   \int_{[0,2\pi]^3} \ln\left(2|\sin\left(\frac{\theta_k 
- \theta_j}{2}\right)|\right)
f(\theta_k) f(\theta_j) \nu(\vec{\theta}) d\vec{\theta} = 0 \]
\[   \int_{[0,2\pi]^3} \ln\left(2|\sin\left(\frac{\theta_k -
 \theta_j}{2}\right)|\right)
\sin(\theta_k) \cos(\theta_k) \nu(\vec{\theta}) d\vec{\theta} = 0 \]
\[\int_{[0,2\pi]^3} \ln\left(2|\sin\left(\frac{\theta_k 
- \theta_j}{2}\right)|\right)
((f(\theta_k) g(\theta_j)+f(\theta_k) g(\theta_j))
  \nu(\vec{\theta}) d\vec{\theta}= 0 \]
\[\int_{[0,2\pi]^3} \ln\left(2|\sin\left(\frac{\theta_k 
- \theta_j}{2}\right)|\right)
((f(\theta_k) g(\theta_i)+f(\theta_k) g(\theta_i))
  \nu(\vec{\theta}) d\vec{\theta}= 0 \]

\end{lemma}

{\bf Proof:}
The idea is simply to note in each case that there are symmetric 
regions of the $\theta$ cube  where the function has  opposite signs 
(the finiteness 
once again follows form the integrability of the logarithmatic singularity). 
 I'll simply indicate the symmetries. 

One key observation is that $\nu$ and  
 $\ln \left( 2|\sin \left( \frac{\theta_k 
- \theta_j}{2} \right) | \right)$ are invariant the 
transformation, $T_c$, were for each $i$  $\theta_i$ goes to 
$\theta_i + c $  modulo $ 2 \pi$; and  
 the transformation $N$ 
sending all $\theta_i$ to $-\theta_i$ modulo $2 \pi$ negates $\nu$
 while of course leaving $\ln \left( 2|\sin \left( \frac{\theta_k 
- \theta_j}{2} \right) | \right)$
 invariant.

Now note for the first integral either $N$ of $N \cdot T_{\pi}$ 
will produce the needed symmetry.  For the second integral 
$N$ alone will work.  
 The Remain  integrals  have  the order four transformation 
$T_{\frac{\pi}{2}}$ producing the four points - two of each possible sign.  

{\bf Q.E.D. (lemma)}

%The next largest term is of size $\ln(r)r^5$ and can be easily computed 
%from the above formula for the angle difference telling 
%us that we must integrate 
%\[ \ln(r) r^5 (- \Delta \phi +  k)(\nu(\vec{\theta}))^2. \]

In $I_2$ note that after integrating the Laplacian term in equation \ref{two}
integrates away and the $k$ term integrates 
to an Euler characteristic - so in the end we get a constant times 
the Euler characteristic  $I_2 = C_2 \chi(M)$, 
with $C_2$ depending on none of the geometry of $g$ or topology of $M$.

%The next term is the $r^5$ and is the smallest term of interest. 
% Each term in the volume function's  power-series contributes
% to it.  Letting $b_i$ be the order $r^2$ piece of the formula 
%for $\psi_{\phi}^i$ we have this term is given by 
%\[  r^5 \left( \frac{-k}{2} \Lambda(\Bbb E_i) \nu(\vec{\theta}) + 
%\frac{- \Delta \phi +  k}{2}  \ln(2|\sin(\Bbb E_i)|(\nu(\vec{\theta}))^2 +
%2 b_i  \ln(2|\sin(\Bbb E_i)|) \nu(\vec{\theta}) \right)\]
%\[ +\frac{ r^5}{2} ({\bf v}^{\perp}_{jk} \cdot \nabla \phi)^2 \cot(\Bbb E_i)\n%u(\vec{\theta})   \]
%\[ + \frac{r^5}{2} (- \Delta \phi +  k) \ln(\Delta \phi -  k)
%(\nu(\vec{\theta}))^2. \]

$I_3$ can be broken up into the four pieces as in its formula 
(equation \ref{threea} -- \ref{threed}).
The first piece (equation \ref{threea}) 
can as above be integrated out to give  $C_3 \chi(M)$,  with $C_3$ 
depending only on the Euler characteristic.

 The last two equations (\ref{threec} and \ref{threed}) can 
be explicitly integrated to give 
\[ - 6 \pi^3  \int_M || \nabla \phi ||^2 +(\Delta  \phi - k)
 \ln(\Delta \phi  - k) dA. \] 

This leaves equation \ref{threeb}  involving the  $b_i$. 
 Looking at the expression for $b_i$ we see that the integral vanishing lemma 
immediately gives us that these terms integrate away to zero, 
with exception of the  terms in the form 

\[ (\frac{2}{3} k + \phi_{xx}) \int_{[0,2\pi]^3} 
\ln\left(2|\sin\left(\frac{\theta_k - \theta_j}{2}\right)|\right)
((-\sin(\theta_k) + \sin(\theta_j)) \cos(\theta_i)
  \nu(\vec{\theta}) d\vec{\theta} \]
\[ (\frac{2}{3} k +\phi_{yy}) \int_{[0,2\pi]^3} \ln\left(2|\sin\left(
\frac{\theta_k - \theta_j}{2}\right)|\right)
((\cos(\theta_k) -  \cos(\theta_j)) \sin(\theta_i)
  \nu(\vec{\theta}) d\vec{\theta} = 0. \]

However the last integral in the vanishing lemma tells us that
 these integral are equal.  So as in the previous piece we are left
 with a term which integrates out to $C_3 \chi(M)$.
% a constant - dependent only on the Euler Characteristic.

So using theorem \ref{exp} our energy is now indeed  in the claimed form.

%\[\frac{\lambda^3}{6}\int_{0}^{\infty}e^{-r^2 \pi \lambda} 
%\left(  C_1 A r^3 + C_2 \chi(M)\ln(r) r^5 +  \chi(M) r^5 \right)   dr  \]
%\[ +  \lambda^3 \pi^3 \left( \int_{0}^{\infty} 
%e^{r^2 \pi \lambda} r^5 dr \right) 
% \int_M || \nabla \phi ||^2 +(\Delta \phi - k)
% \ln(\Delta - k)   dA. \]
%Now integrating out the $r$ variable and including the higher 
%order terms and using an estimate as in section ?, we arrive at the needed for%mula. 

{\bf Q.E.D. (formula \ref{formu})}

Now we can  compute the energy in formula \ref{cen} simply by noting
that the terms in the formula with  $D_i$ constants in them 
cancel (in the prelimit even) and that we
are left precisely  with the needed terms.  
Also noter the formula for $E_1$ is immediate.
With these we can verify the formula for $E_2$  and 
prove theorem 
\ref{lapy}.

{\bf Proof (theorem \ref{lapy})}
It is immediate that  amongst metrics of the same area, $A$,  that 
\[ E_2 = E_1 - E = 
- \int_M k_h \ln(-k_h) dA_h + 
 \int_M ||\nabla\phi||^2 + 
 (\Delta \phi - k) \ln(\Delta \phi - k )d A \]
\[ = - \int_M (\Delta \phi  -k )
\ln\left( e^{-2 \phi}(\Delta \phi -k)\right) dA + 
\int_M ||\nabla \phi||^2 + 
 (\Delta \phi - k) \ln(\Delta \phi - k )d A \]
\[ =   \int_M 2 \phi(\Delta \phi  -k )dA + 
\int_M ||\nabla \phi||^2 d A  = - \int_M  ||\nabla \phi||^2  + 2 k \phi dA . \]

We now state the beautiful formula do to 
  Polyakov (see \cite{Po} or \cite{os})
 for the  the determinant of the Laplacian 
    
  \[ \ln(\det(\Delta_h)) = - \frac{1}{6 \pi} \left( \frac{1}{2}
  \int_M {|\nabla \phi|}^2 dA + \int_M k \phi dA \right) + \ln(A) + C.\]

Plugging in above we arrive at the needed
\[  E(h) = 12 \pi \ln(\det(\Delta_h)) - 12 \pi \ln(A)  + C. \]
which is the claimed formula.

{\bf Q.E.D}
\end{subsection}
\end{section}

\begin{section}{A Proof of the Uniformization Theorem}\label{bigproof}
The goal here is to give the proof of the 
metric uniformization
theorem  indicated in section \ref{conu}.  
It is a direct analog in infinite dimensions of the  
proof of the angle system uniformization  
presented in section \ref{prot}.
To see this it is useful to  organize the finite dimensional  proof 
into 4 basic steps, which will be mimicked amongst metrics.

\begin{subsection}{The Discrete Uniformization Proof Reviewed}\label{revisit}

The negative curvature Delaunay angle system uniformization proof:

{\bf Step 0.}
We defined what what we called  the 
negative curvature Delaunay angle systems conformal to a fixed angle system $y$
and called it 
$\frak{N}_y$, 
which turns out to be 
 nice convex set.
Then we  placed upon it an  energy $E$
\[E(x) = \sum_{t \in \frak{P}} V^t(x) \]
 which was continuous on 
$\bar{\frak{N}}_y$.

{\bf Step 1.}
Now we observed that at least one point of maximum energy must exist,
here by the rather trivial observation that space 
$\bar{\frak{N}}_y$ upon which our continuous energy 
lives is in fact  compact.

{\bf  Step 2.} 
At this step we showed that any  $x$ where $E$ assumes its maximal value  
is in fact in   $\frak{N}_y$.  
To accomplish this recall we took  
a point $x$ on the boundary of $\frak{N}_y$ and constructed a line  
$l(s)$ satisfying $l(0)=x$ and that there is an $\epsilon > 0$
such that for all $s \in (0,\epsilon)$ we have that $l(s) \in \frak{N}_y$ and 
$E(l(s))$ is increasing in $s$.  Clearly now the continuity of $E(l(s))$ on 
$[0,s)$ makes it impossible for $E$ to assume its maximum value 
at any  boundary point.

%Recall the actual method used was to compute $\frac{dE(l(s))}{ds}$ on 
%$(0,\epsilon)$ using the explicitly computed  differential of $E$
%having the additional property that 
%there is an $\epsilon >0$ such that 
%the  continuous function $E(L(s))$ on $[0,\epsilon)$ 
%is increasing in $(0,\epsilon)$.  
%It is worth is requires looking at the differential and 
%in lemma ? we saw   that in fact any direction will do the trick.
%So $E(x)$ could not have been maximal at a boundary point. 

{\bf Step 3.}
At this step we verified that  any point $x \in \frak{N}_y$ 
where $E$ is maximal 
  is in fact a  uniform angle system.
This followed  immediately 
by examining what the differential of $E$ vanishing at $x$ 
implied about $x$ (see observation \ref{disobse}).

%$dE^x = \sum_{\alpha_i \in \frak{P}} E_i(x) \alpha^i$ with 
% \[ E_i(x) = \frac{1}{2} \left( 
%   \ln\left(\frac{\cosh(a)-1}{2}\right) - \ln\left(\frac{\cosh(b)-1}{2}\right)
%   -\ln\left(\frac{\cosh(c)-1}{2}\right) \right), \]

%\begin{observation}
%Thus far all the arguments work just as well as in the positive curvature
%case.  In particular there is a conformally equivalent uniform 
%positive curvature Delaunay  angle system on the 
%sphere relative to any fixed one. 
%\end{observation}

{\bf Step 4.}
We proved uniqueness of the point $x$ where $E$ achieves its maximum.  
To do this we took a second point $y$ and connected it with a line $l(s)$ to our 
$x$ and note that 
 \[ \frac{d^2 E}{ds^2}(l(s)) < 0\] 
for all $s$ such that $l(s)$ remains in $\frak{N}$.  
So $y$ could not also be maximal (see lemma \ref{lemer}).

\end{subsection}

\begin{subsection}{The Indiscreet Proof}

The proof in the previous section 
can be carried out in the indiscreet world step by step.  
%Warning: different proofs using $E$ to uniformize metrics utilize 
%different analytical objects and the
%proof I will give replicating the above proof utilizes easy theorem but some 
%potentially unfamiliar objects

{\bf Step 0.}
We need the correct analogs 
of an $\frak{N}$ and $\bar{\frak{N}}$ on which to interpret the energy 
\[ E(\phi) 
 = - \int_M || \nabla \phi||^2 
+ (\Delta \phi -k) \log(\Delta \phi  - k) dA .\] 
Namely which 
$\phi$ do we use to conformally change our initial metric. 
Clearly a rescaling will change nothing both since it effects 
no angles in the 
discrete model and since $E$ is clearly  scale invariant  
on $C^{\infty}$.
In fact with this in mind we should feel free to re-scale the initial metric 
and for convenience let's assume it  area is $- 2 \pi \chi(M)$. 
One way to eliminate the possibility of rescaling is to  demand
\[ \int_M \phi dA = 0. \]
We'd also like to 
restrict to $\phi$ with negative curvature, i.e. 
$k_{\phi} = e^{-2 \phi}(-\Delta \phi +k) < 0$.    So
the first guess at a reasonable function space might be
\[V =\{ \phi \in C^{\infty}(M) \mid \int_M \phi dA = 
0 \mbox{ and } \Delta \phi - k > 0\}.\] 
However the energy on here would make 
using a compactness argument in step 1 rather difficult.
%would interact with only three of the 
%family of $sup$ norms forming $C^{\infty}$, so would give us no control
%over the compactness of sets or the like. 
The trick to producing a place where a compactness argument will 
work is to take the closure of $V$ with respect to a 
norm which $E$ interacts with
in a sensible way.  The energy being convex in fact means we can 
essentially close $V$ under the energy viewed as a norm.  Such  Banach 
space are well studied and called Sobolev-Orlicz spaces.

\begin{subsubsection}{A BRIEF introduction to  Sobolev-Orlicz spaces}
\label{sob}

In this section we recall several well known theorems concerning  
Banach spaces and in particular the 
 Sobolev-Orlicz spaces.   

One of the key uses of Sobolev-Orlicz spaces (introduced below) is to 
produce Banach spaces $B$ where we have no control over $L^p$ 
growth  for $p>1$ 
yet are still able  to represent $B$ 
as the dual of a 
second Banach space, i.e.  $B = D^{\star}$.   
Recall that Gel'fand's theorem  tell us that even $L^1([0,1])$ cannot be 
realized 
a the dual of any Banach space, so $L^1$ would  not do.

The reason we would like $B = D^{\star}$ is that such a 
relation gives us sensible notions of compactness in the weak 
topology on $B$.  

\begin{theorem}[Alaglu's Theorem]\label{alag}
A closed and bounded set in the norm topology is weak compact.
\end{theorem}

To identify such set's we will recall on of the most well known convexity 
theorems, namely...

\begin{theorem}[Mazur's Theorem] \label{maz}
A closed and convex set in $B$ is also closed in the weak topology.
\end{theorem}

For proof's of these results see \cite{La}.
Now let's produce the spaces with the $B = E^{\star}$ property to which these theorems will apply.

An Orlicz space is an $L^p$ type space using a different 
convex function than $\Phi(t) =t^p$. In fact its good to specify 
the class of convex functions on $[0,\infty)$ 
(the Young functions) of use here.
A young function is 
\[ \Phi(t) = \int_0^{t} \phi(s)ds \]
where $\phi(s)$ satisfies
\begin{enumerate}
\item
$\phi(s) \geq 0$ for $s >0$
\item
$\phi(s)$ is left continuous
\item
$\phi(s)$ is non-decreasing on $(0,\infty)$ 
\item
$\phi(\infty) = \infty$.
\end{enumerate}

We would like to form a norm which behaves something like 
\[ \rho_{\Phi}(f) = \int_0^{\infty} \Phi(|f(x)|)dA. \]
In fact the space of function with satisfy $\rho_{\Phi}(f) <\infty$ will
be essential to us, and we will denote it $\tilde{L}_{\Phi}$.

To actually implement 
this we take what you can of the inverse of $\phi$, namely let  
\[ \psi(t) = \sup_{\phi(s) \leq t} s .\]

Then let 
\[ \Psi(t) = \int_0^{t} \psi(s) ds \]
and call it the the Orlicz conjugate of $\Phi(t)$.

Perhaps the most important example in this context   
are  the following conjugate relationships....
\[ \phi(s) = \left\{ 
\begin{array}{ll} 0 & 0< t <1 \\ 
\log(t) + 1 & t \geq 1
\end{array} \right.\]

\[ \Phi(t) = t \log^+(t) \left\{ 
\begin{array}{ll} 0 &  0 \leq t <1 \\ 
t \log(t) & t \geq 1
\end{array} \right.\]

\[ \psi(t) = \left\{ 
\begin{array}{ll} 1 &  0 < t <1 \\ 
e^{t-1} & t \geq 1
\end{array} \right.\]

\[ \Psi(t)= e^{t-1}_+ = \left\{ 
\begin{array}{ll} t  &  0 \leq t <1 \\ 
e^{t-1} & t \geq 1
\end{array} \right.\]

With these notions let 
\[ || f ||_{\Phi} = sup_{v \in \tilde{L}_{\Psi}} \int|f(x) v(x)| dA.\]

Then this forms a norm on the space of measurable function with 
$||f||_{\Phi} < \infty$ (call it $L^{\Phi}$).
% so that 
%\[ ||f||_{\Phi} \leq \rho(f) +1.\]

The $\Phi$ come in two flavors the happy ones which satisfy
the existence  a $T>0$ and $k>0$ such that 
\[ \Phi(2t) \leq k \Phi(t),\]
and the sad ones which don't.  Notice $t \log^+(t)$ 
is happy and it conjugate $e^{t-1}_+$ is not.

Let $E_{\Phi}$ be the closure of $C^{\infty}$ in $L^1$ under this norm.
For happy $\Phi$ we have that $\tilde{L}_{\Phi} = L_{\Phi} = E_{\Phi}$ 
is a separable Banach space  
and that the norm interacts with $\rho_{\Phi}$ nicely.  For example if 
$f_n \rightarrow_{\Phi} f$ then
$\rho_{\Phi}(f_n) \rightarrow \rho_{\Phi}(f)$ 
(this is a special consequence of what is known as mean convergence)
%the mean convergence result is usually proved and this together with 
%$|\Phi(b) -\Phi(a)| \leq \Phi|b-a|$ for any convex function).
For sad  $\phi$ we have $\tilde{L}_{\Phi} \subset  L_{\Phi} \subset  E_{\Phi}$ 
and that 
these inclusions are always proper,  also $L_{\Phi}$ fails to be separable 
and $\tilde{L}_{\phi}$  fails to even be a vector space. 

$L_{\phi}$ and $L_{\Psi}$ are reflexive if and only if $\Phi$ and $\Psi$ 
are happy.  However it is alway the case that $L_{\Phi} = (E_{\Psi})^{\star}$, 
the key property discussed in the first paragraph of this section.  
%Its also worth  noting that a sequence converges 
%in any Orlicz norm implies it $L^1$ converges hence 
%has a subsequence converging almost surely.  

With these space the Sobolev-Orlicz  spaces are easy understood.  
I will only present and need a very special case, but 
everything here works in complete generality (see the very nice \cite{Don}). 
We will 
embed $C^{\infty}$ into 
$L^2(M)  \times L^2(\Gamma(TM)) \times L_{t \log^+(t)}(M)$
via $I(u) = (u, \nabla u, \Delta u)$ 
and take its closure in the Banach norm. We arrive a Banach space 
$B$.  Just as above even though $e^{t-1}_+$ is bad 
we can realize $B$ as $E^{\star}$ were $E$ is the closure of $C^{\infty}(M)$
in $L^2(M)  \times L^2(\Gamma(TM)) \times L_{e^{t-1}_+}(M)$ under the same embedding.

Its worth noting at this point that we in fact have certain obvious continuous 
inclusions of the classical Sobolev spaces.  Let $H^{k,p}$ be the usual 
Sobolev space where we control the $L^p$ norm of the first $k$ derivatives.  
With these we clearly have 
 have the following continuous inclusions
\[ H^{2,2} \subset B \subset H^{1,2}. \]
 %and that in fact the second inclusion is compact by the usual S
This fact  gives us some nice functions in $B$.  Namely 
we have the Fredholm Alternative assuring us that 
\[ \Delta (H^{2,2,}) = \{ f \in L^2 \mid \int_M  f dA = 0\}; \]
so for any mean zero $L^2$ function $f$  we can construct a function in 
$g \in B$ 
such that $\Delta f = g$.

%It worth mentioning that the Sobolev-0Orlicz space clears up nicely certain as aspects of the Sobolev embedding theorem.  Every Young $\Phi$ in fact has a family of duals $\Psi_k$ relative to how many derivatives one wishes to control.  A wondrous has is that on 2 surface the $\Psi_1$ relative to $\Phi = |t|^p$ is\[ \Psi_1  = C_P |t|^{\frac{2p}{2 -p}}\]for $p < 2$  and \[ \Psi_1  =  e^{|t|} -t-1 \]for $p =2 $.With these the Sobolev embedding theorems becomes\[ W^{1,p}\subset L_{\Psi_1}\]with the inclusion being continuous for $1 < p \leq 2$ and\[  W^{1,p} \subset C(M)\]for $p >2$.  Hence we have demystified the role of $2$ in the theorem (something that always bugged me).
%In fact if one does {\bf not} care about continuity 

The last fact is a certain set inclusion, 
namely  
\[ H^{1,2}(M) \subset L_{e^{t^2}-1}(M). \]
In particular  the mapping 
$\eta$ such that $\eta(u) = e^{u}$ takes
$H^{1,2}$ into $L^p$ for all $p$ and 
$H^{k,2}$ into $H^{k-1,2}$.  I will refer to this fact as 
Trudinger's inequality.
    
%For nice discussion of this (along with some curious 
%compact embedding results) see \ref{Ta}.

\end{subsubsection}   

\begin{subsubsection}{Step zero continued}
With our introduction to the needed spaces  
out of the way we may proceed with step zero 
by letting 
$\bar{\frak{N}}$ be
the closure of 
$V$ in $B$, and $\frak{N}$ the subset of this where 
$esssup(\Delta \phi  -k) > 0$.
One  key property of the space $B$ in this context is that  $B = E^{\star}$. 
so by Mazur's theorem 
the convex set $\bar{\frak{N}}$ is closed in both the norm and weak topologies.
The remainder of step one can be summed up in the following lemma 
assuring the continuity of $E$.

\begin{lemma}\label{cont}
$E$ is continuous on $\bar{\frak{N}}$ in both the 
norm topology and weak topologies. 
\end{lemma}
{\bf proof:} 
To see the norm topology case let 
$\phi_n \rightarrow_B  \phi$  implies
$\Delta \phi_n $ converges to $\Delta \phi$ 
in the $|| \cdot ||_{L \log(L)}$ norm.  
So $\Delta \phi_n -k$ will converge to $\Delta \phi -k$   as well, since 
$k \in C^{\infty}(M)$.  In particular since the $L \log^+(L)$ norm is happy 
we have 
$\int_M (\Delta_n \phi -k) \log^{+}|\Delta_n \phi -k| dA   
\rightarrow \int_M  (\Delta \phi -k) \log^{+}|\Delta \phi - k| dA$.  
Now observe that on $\bar{\frak{N}}$  
$\rho_{L \log^+(L)} = \int_M (\Delta \phi -k) \log^+(-\Delta \phi - k) dA $ 
differs from 
$E$ by a continuous and bounded function.  
So by  the fact that the norm convergence implies $L^1$ convergence 
we have from the dominated convergence theorem  that 
$\int_M (\Delta_n \phi -k) \log|\Delta_n \phi -k| dA  
\rightarrow \int (\Delta \phi -k) \log|\Delta \phi - k|dA$ as needed.

The weak topology assertion follows from the fact that $
E$ is convex on $\bar{\frak{N}}$ 
hence a convex function on its closure, so in the norm topology 
$E^{-1}([a,b])$ is closed
and convex hence by Mazur's theorem closed in the weak topology.

\qed

{\bf Step 1.}
Now we would like to proceed as in the discrete case and use a 
compactness arguments to 
assert the existence of a function achieving the maximum.
Since $\bar{\frak{N}}$ is not quite compact in this case,  we 
must do a little work to see that it is compact enough.  The first 
thing to note is 

\begin{lemma} \label{boundd}
$\sup_{\bar{\frak{N}}} E \leq  0$. 
\end{lemma}
 
{\bf Proof:}
The boundedness  of $E$ on $\bar{\frak{N}}$ 
follows form the fact that for any $\phi \in B$ that 
$\int_M \Delta \phi = 0$  since $\Delta \phi$ is the 
$L^1$ limit of $C^{\infty}$ functions with this property.  
So we have that by 
Jensen's inequality that
\[ 0  =  
(\int_M \left(-\Delta \phi + k)\frac{dA}{- 2 \pi \chi(M)} \right) 
\log \left| \int_M(-\Delta \phi  + k)\frac{dA}{-2 \pi \chi(M)} \right| \]
\[ \leq 
\int_M (-\Delta  \phi +k) \log|-\Delta  \phi +k| \frac{dA}{-2 \pi \chi(M) }. \]
This along wit the obvious fact that 
$- \int_M || \nabla \phi||^2 dA \leq 0$ gives us the needed bound. 

\qed

Denote  the finite number $sup_{\bar{\frak{N}}}E$ as $m$.
From the lemma \ref{cont} on the continuity of $E$ in the weak topology
$K =  E^{-1} ([m,m+a])$ is  
weak closed and convex. 
In fact $K$ is weak compact.  
To see this it is enough by Alaglu's theorem
to see that the norm is bounded on $K$.  By Poincare 
inequality $|| \phi||_2 < C || \nabla \phi ||_2$
so 
\[ || \phi ||_B^2 < C_1(||\nabla \phi ||^2_2 +   
||\Delta \phi ||_{L \log^+(L)}) \]
\[< C_2( ||\nabla \phi ||^2_2 + \rho_{L \log^+(L)} ( \Delta \phi -k) )
< C_2 (E(\phi) + C_3) < C_2(m+a +C_3) . \]

Now just as in the discrete world we 
have a continuous function, $E$, on a compact set $K$
and hence we have at least one  point achieving the maximum value.

{\bf Step 2.}
Now just as in the discrete case we need to control the boundary.
Suppose a maximum occurs on the boundary at $\phi$. 
Just as in the discrete case we will construct a direction $\psi$ and a line 
$l(s) = \phi + s \psi$ 
is contained in $\bar{\frak{N}}$ such that  
$\frac{d E(l(s))}{ds} > 0$ for all 
$s \in (0,\epsilon)$ for some $\epsilon >0$; 
hence contradicting the maximality of $\phi$. 
%To construct $\psi$ recall by the Fredholm alternative that 
%$\Delta(H^{2,2}) = \{f \in L^2 \mid \int_M \phi dA =0\}$.  

$\phi$ being on the boundary implies that 
\[ M_{\delta} = \{ x \in M \mid \Delta \phi - k <  \delta \},\]
has measure  $m_{\delta} > 0$ for all $\delta >0$.    
Now since $\int_M \Delta \phi dA = - 2 \pi \chi(M) > 0$ there is certainly 
an interval  $[a,b]$ such that $a \geq  e$ and the set  
$S = \{ x \in M \mid \Delta -k \in [a,b] \}$ 
has measure $s >0$.
For each $\delta$ let 
\[ f_{\delta} = 
\frac{\chi_{M_{\delta}}}{m_\delta} -  \frac{\chi_{S}}{s} .\]
Note by the Fredholm alternative that there is a
$\psi_{\delta} \in H^{2,2,} \subset B$ such that $\Delta(\psi_{\delta}) 
= f_{\delta}$, and further more by altering this function with a constant that 
one can assume $\int_M \psi_{\delta} dA =0$.
Using this $\psi_{\delta} \in H^{2,2}$ direction we see that for small enough 
$s$ that indeed $l(s)$ is in $\bar{\frak{N}}$ and  that 
\[ \frac{d E(l(s))}{ds} = \int_M   f_{\delta} \cdot l(s)  
- f_{\delta} \log(-\Delta (l(s)) + k)dA .\]
Since  $\log|x|$ tends to $- \infty$ as $x$ tends to zero  
we see that for small enough that  $\delta$  that 
 \[ \int_M  
- f_{\delta} \log(-\Delta (l(s)) + k)dA \] 
can be made a large as we'd like.
Note that in $M_{\delta}$ we have $\Delta (l(s))$ is bounded and hence 
by the Green's function representation of $l(s)$ we have that 
$f_{\delta} \cdot l(s)$ is bounded simultaneously for all small enough 
$\delta$.  
So  indeed  $\delta$ can be chosen so 
$\frac{d E(l(s))}{ds} > 0$ for small $s$.  So boundary maxima 
are impossible, and we have an internal maxima.

{\bf Step 3.} So now we have that our point of maximal energy is internal.  
Note at such a maxima 
$\Delta \phi > k$ in essential supremum so 
$\ln(k_\phi) = -2 \phi + \log(\Delta \phi  -k) \in L^1(M)$. Using $l(s)$ as 
above we see for each 
$f \in C^{\infty}(M)$ such that $\int_M f dA =0$ we have 
\[  \frac{d E(l(s))}{ds} =  \int_M f \log( k_{\phi}) dA  =0. \]
So $\log( k_{\phi})$ is a constant as an $L^1$ function.  
In particular exponentiating we see that  $k_{\phi}$ is a constant.
In fact note 
\[ \Delta \phi = C e^{2 \phi} + k,\]
and by Trudinger's inequality $e^{2 \phi} \in L^2$ so $\phi \in H^{2,2}$ 
by elliptic regularity.  So  $e^{2 \phi} \in H^{1,2}$ and 
by elliptic regularity again $\phi \in H^{3,2}$.  
Continuing this $\phi$ and hence $k_{\phi}$ are in fact in all $H^{k,2}$ 
and hence by the Sobolev embedding theorem in $C^{\infty}$.

%\begin{observation}
%So far nothing about this argument changes for the sphere, and there exist a
%constant curvature metric uniform to any given one.  
%(Of course for the sphere this metric will not be unique since 
%there are conformal mapping which are not isometries, i.e.given such a 
%mapping $\eta$ and a uniform metric $u$ we have $\eta^*(u)$ are 
%both uniform with different conformal factors.)
%\end{observation}

{\bf Step 4.}
Now we can easily get uniqueness in the $\chi(M) < 0$ case, 
exactly as in the discrete case. Take two now $C^{\infty}(M)$ 
solution $\phi_1$ and $\phi_2$
and note form the second Frech'et derivative that the line $l(s)$ 
in $C^{\infty}$ connecting them satisfies
\[ \frac{d^2 E}{ds^2}(l(s)) = 
-\int_M \| \nabla \phi_1 - \nabla \phi_2||^2 
+ \frac{(\Delta \phi_1 - \Delta \phi_2)^2}{\Delta (l(s)) - k} dA < 0 \]  
 for all $s$ with $l(s)$ in $\frak{N}$.  
So as in the discrete case the point of maximum energy is unique. 
\end{subsubsection} 
 \end{subsection}

\end{section}

\end{chapter}

\begin{chapter}{Spheres and Tori}\label{tori}

This chapter is dedicated to exploring the ideas of the 
previous chapters in the cases of $\chi(M) \geq 0$. 
The cases of primary interest are the torus and sphere cases.
As in the $\chi(M) <0$ case the fundamental object 
needing exploration is the class of polyhedra related to a triangulations,
and the volume energy associated to this class.  In  both cases there
is an intrinsic class of such polyhedra but the story presented here in
the $\chi(M) <0$ case experiences difficulties. 
The nature of the difficulties is very different for the 
torus and the sphere.  

In the toroidal case the class of 
polyhedra needed has been studied in the literature in
\cite{Be} and \cite{Ri2},
and the issues in chapter two have already been essentially  dealt with. 
I will remind the reader of the issues involved in section \ref{toroid}, 
and for now
only highlight the differences with the $\chi(M) <0$ case.
In the torus case the  energy  is  significantly 
easier to deal with  and analogs to 
theorem \ref{tri1} and theorem \ref{cir1} exist and can be 
proved with the same methods.
 However 
the story is different in the important respect 
that the linear part of the problem sends one immediately into the land of 
zero curvature world, a phenomena which occurs 
in both the discrete and continuous cases.
So in the continuous metric case there is no non-linear 
{\bf metric} story at all.
This case in fact demonstrates the important fact that
perhaps a better  continuous analog  of the discrete world would be 
a connection or an affine like structure. This because, as we shall see, 
the discrete uniformization does not 
 produce  flat structures,  but
 rather affine structures.  In any case the ideas of  
chapter three breaks down at a rather fundamental level.

The spherical case is quite the opposite, and on some level the 
ideas don't break down at all.  As in the $\chi(M) \leq 0$ cases there is an 
appropriate class of ``intrinsic'' polyhedra, and the one gets a natural 
volume energy which has critical points precisely  at  uniform structures.  
Unlike in the toroidal case 
the randomization goes through perfectly  to produce an energy on the 
positive curvature metrics.  Every thing looks good except 
now the energy is 
no longer nice at all.  Analogs the theorem \ref{tri1} and \ref{cir1} exist 
(see \cite{Ri1})
but appear much harder to prove with the intrinsic methods used  in the 
$\chi(M) <0$ case. In particular the energy fails to be 
convex and the boundary behavior becomes very difficult to control.  
In fact both the 
discrete and   continuous proofs of the previous chapter 
fail in fundamental ways.  Of course intrinsic proofs may still exist, 
but as is often the case with spherical  uniformization (see \cite{os} 
and \cite{ha}) will involved significantly more drama. 
 
\begin{section}{The Toroidal Case}\label{toroid}
In the discrete toroidal  world one must first replace $\frak{N}$ with 
angle systems which  have zero rather than negative curvature.  This of course 
gives us a significantly smaller space of angles and in particular the 
conformal deformations preserve not only the $\theta^e$ but 
also this curvature condition, and are the span of the 
$w_v$ vectors in figure \ref{tor}
over all the vertices.

\begin{figure*}
\vspace{.01in}
\hspace*{\fill}
\epsfysize = 1.5in  
\epsfbox{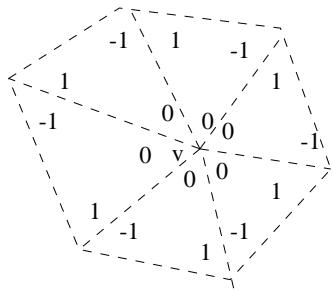}
\hspace*{\fill}
\vspace{.01in}
\caption{\label{tor} The $w_v$ Vector}
\end{figure*}

To   construct the polyhedra first view the Euclidean 
plane as the boundary at infinity of hyperbolic space 
in the upper-half space model. The polyhedra are now constructed by 
taking the union of the ideal simplexes over 
each of the Euclidean triangles.  
Note the energy formula in the $\chi(M) < 0$ agrees with this construction 
and becomes twice the sum of the volumes of all the 
ideal tetrahedra over a fundamental domain.  As before the volume is simply 
the sum of the volumes of the  ideal 
tetrahedra corresponding to the individual  Euclidean triangle angles
in the angle system. In other words 
\[ E(x) = 2 \sum_{t \in \frak P} V_t(x) \] 
where  if $d^t(x) = \{ A,B,C\}$ we have 
\[ V_t(x) = \Lambda(A) + \Lambda(B) + \Lambda(C). \]

This energy remains convex and boundary controllable 
and all of chapter two carries over with the most interesting 
point being why at the critical points of $E$  fit together.

\begin{observation}\label{tordif}
At a critical point of $E$ the triangles fit together 
to form an affine structure on the torus.
\end{observation}

{\bf Proof:}
As far as I'm aware the idea in this proof has its origin in 
Bragger's \cite{Be}.
The above formula tells us that 
from the formula for the Lobacevskii function we have    
$dE^x = \sum_{\alpha_i \in \frak{P}} E_i(x) \alpha^i$
with  
   \[E_i(x) =  - \ln\left(\sin(A^i)\right) \]

So at a critical point with $n$ faces $t_i$ in its flower 
at  $v$ and angles labeled $A^i_{\pm}$ we have  
\[0  = dE_x(w_v) =    \ln\left(\frac{\sin(A^1_+)}{\sin(A^2_-)}
\frac{\sin(A^2_+)}{\sin(A^3_-)} \cdots 
\frac{\sin(A^n_+)}{\sin(A^1_-)} \right).\]
Now lets attempt to fix our edge lengths.  
We will denote   the edge length opposite to $A^i_{\pm}$ as $a^i_{\pm}$.
Each triangle can be scaled with its angles
preserved since we are in the Euclidean plane.  So fix the size of $t_1$.  
Now scale $t_2$ so that $a^2_- = a^1_+$.  Continue this until the size of 
$a^n_-$ has been fixed.

Now  from the law of   
$\sin$s $0 =dE(w_v)$ gives us 
\[ 1  =   \frac{a^1_+}{a^2_-}
\cdots \frac{a^n_+}{a^1_-}   = \frac{a^n_+}{a^1_-}.\]
So the entire flower fits together.

Now we have open sets which are affine related  
in overlaps, so an affine structure.

\qed

With this observation the other ideas essentially work out in the same way.
\end{section}
\begin{section}{The Spherical Case}\label{sphere}
The initial discrete set up for the sphere is identical to the $\chi(M) <0$ 
case except the use of 
positive rather than negative curvature.  The class polyhedra is very 
simple to construct.  Given a triangulation of the sphere take the 
convex hull of the vertices.  The volume formula turns out to once again 
 agree with the 
negative curvature case.   
This is quite a nice fact.  
To see it  view the  sphere at infinity 
in the ball model from the origin.  Note form this view point 
the angles you see in the ideal polyhedra are precisely the 
angles in the triangulation as the sphere understands them.
Now for each triangle on the sphere at infinity with $d^t(x) = \{A,B,C\}$ 
 form 
the three ideal vertexes tetrahedra by taking the convex hull of the 
these vertices at infinity with the origin, see figure \ref{tet}.

\begin{figure*}
\vspace{.01in}
\hspace*{\fill}
\epsfysize = 1.5in  
\epsfbox{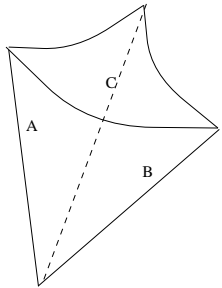}
\hspace*{\fill}
\vspace{.01in}
\caption{\label{tet} A Hyperbolic Tetrahedron}
\end{figure*}

\begin{fact} \label{spherevol}
The volume of the above tetrahedra $V_t(x)$ is given by 
\[2 V_t(x) = \Lambda(A) +\Lambda(B)+\Lambda(C)+ 
\Lambda\left(\frac{\pi-A-B-C}{2}\right) \]
\[
+\Lambda \left(\frac{\pi+A-B-C}{2}\right)
+\Lambda\left(\frac{\pi+B-A-C}{2}\right)+
\Lambda\left(\frac{\pi+C-A-B}{2}\right) .\]
\end{fact}
 {\bf Proof:}
To see this extend the geodesics in the tetrahedra 
and take the convex hull of this arrangement.  We get an ideal octahedron.
Using the three new points at infinity and the 
origin note we have in this octahedron
a symmetric copy of our original three ideal vertexed tetrahedra.  

\begin{figure*}
\vspace{.01in}
\hspace*{\fill}
\epsfysize = 1.5in  
\epsfbox{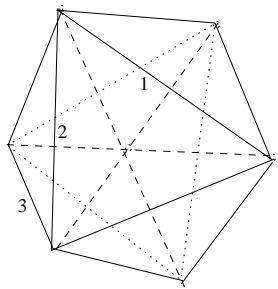}
\hspace*{\fill}
\vspace{.01in}
\caption{\label{oct} The Ideal Octagon}
\end{figure*}

Each edge $e_i$ of the octahedron corresponds to an ideal tetrahedra 
$T_i$, see figure \ref{oct} where three particularly 
relevant tetrahedra are labeled. 
Now simply note  that by using both copies of the 
three ideal vertexed tetrahedra that the needed volume can be expressed as 
half of $Vol(T_{1}) - Vol(T_{2}) + Vol(T_{3})$. Now the fact that an
ideal tetrahedra has its volume given by summing the Lobacevskii 
function over the angles meeting at a vertex gives
the needed formula.   
\qed

Note by spherical geometry that
 \[ \cos(a) = \frac{\cos(C) \cos(B) + \cos(A)}{\sin(B) \sin(C)},\] 
 so the same exact  
computation as in the $\chi(M) < 0$ case tells us a critical point is 
uniform.  However convexity and boundary control are both lost, 
and since we may 
use the same formula as in section \ref{prot} this is easy to see.  
The random computation goes through as in the $\chi(M) <0$ case, 
and we can see 
immediately from 
formula  \ref{hess} that we lose convexity.

It is worth noting that we 
 should have expected problems, at least with uniqueness.  
Namely there is a sort of Gauge group sitting around, 
and in the metric world
it corresponds to the fact   
there are 
conformal transformation of the standard sphere which 
fail to be isometries.
In fact given a uniform metric there is a 
three dimension space of distinct $\phi$ which remain constant curvature.
This is fun to witness in terms of triangulations where it indicates that 
we may expect distinct sets of angles to be be conformally equivalent and still
fit together.   
To see it  fix the vertexes of a triangulation and 
move what we view as the origin of  hyperbolic space
in the above construction away form the balls origin in the model.  
Then moving it back to the model's 
origin produces a topologically equivalent triangulation 
with distinct angles which clearly fits together.  
Although this observation is unfairly
 mixing our two notion of conformal
change it is still  indicates that uniqueness of a uniform structures 
in the  discrete world 
should not be common.

\end{section}
\end{chapter}

\begin{chapter}{Appendix: A Less Pleasant Proof of 
Lemma \ref{lit}}\label{eek1}

Here I will present an alternate proof of lemma \ref{lit}, and 
arrive at the slightly stronger condition of needing circles 
to only be on disks of
radius less than $min\{\frac{i}{6},\tau\}$ rather than on 
circles of radius less than $\frac{i}{8}$ as in section 3.1.2.

The trick to this proof of lemma  $\ref{lit}$  is to understand the curves satisfying $d(p,z) - d(q,z) = 0$; with $d(p,z)$  with less than $min \{ \frac{i}{6}, \tau \}$.  This because any point on such a curve corresponds to the center of a circle going through both $p$ and $q$, and if a triple $\{ p,q,o \}$ lives on a circle then the corresponding curves for each pair in the triple must intersect at the the point corresponding to the center of this circle.

 Before getting started there are a few basics pieces of notation convenient to introduce here: if $d(p,q) < i$ call  $\hat{\gamma}_{p,q}$ the shortest length geodesic segment between $p$ and $q$ (it is well defined by lemma $\ref{normal}$), and let $\hat{\gamma}_{p,q} \subset \gamma_{p,q}$ be the connected component of the geodesic contained in any set we happen to be exploring with $\hat{\gamma}_{p,q}$ in it.
For example relative to $B_{i}(p)$,   $\gamma_{p,q}$ is the geodesic splitting $B_{i}(p)$ into its two distinct  ``sides'' (simply look in normal coordinates).
Denote  as $\{ \hat{\gamma}_{p,q} \}^C$ the two components of 
$\{ \gamma_{p,q} \} - \overline{\hat{\gamma}_{p,q}} $.  We will also find it useful to name the midpoint of  $\hat{\gamma}_{p,q}$ - called it $m$.  (see figure $\ref{note}$ for periodic notation reminders).
  For the remainder of this section denote $d(p,z)$ as $ D_p(z)$, since the differential is usually represented with a $d$. 

\begin{figure*}
\vspace{.01in}
\hspace*{\fill}
\epsfysize = 2in 
\epsfbox{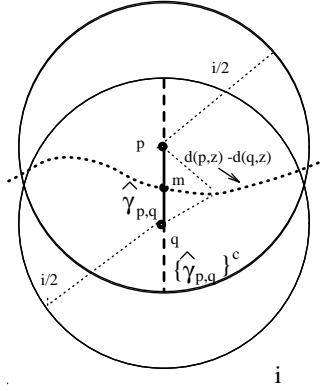}
\hspace*{\fill}
\vspace{.01in}
\caption{\label{note} The Notation}
\end{figure*}

 To get started it  is in fact  useful to consider the more general curve of the type  $D_p - D_q = c$, were $c \in \Bbb R$.  By the implicit function theorem when $d(D_p - D_p) \neq 0 $ and is defined, the solution to this equation is locally a curve $z(t)$ with $\dot{z}(t) \neq 0$; where $z(t)$ satisfies $ d(D_p - D_p)(\dot{z}(t)) =0$.  In fact any such curve will parameterize a solution.  It is useful to rephrase this one form business in terms of its dual object the gradient.
Recall the gradient of a function $f$ is the unique vector field  $\nabla f$ satisfying   $< \nabla f, v> = df(v)$ for all $v$ at every point.  So we may rewrite our differential relation $d(D_p - D_p)(\dot{z}(t)) = 0$ as  $<\nabla D_p  - \nabla D_q, \dot{z}> = 0$.

To get a grip on this differential relation it is first useful to spend a moment contemplating $\nabla D_p$.

\begin{sub-lemma}[The Distance Gradient Sub-lemma]
\label{silly}
Let $d(p,q) < \frac{i}{2}$ then:
\begin{enumerate}
\item
In $B_{\frac{i}{2}}(p) - \{ p \} $ we have $\nabla D_p$ is unit length with its integral curves the unit speed geodesics.
\item
In   $B_{\frac{i}{2}}(p) \bigcap B_{\frac{i}{2}}(q)$ we have  $\nabla D_p  = - \nabla D_q$ on $\hat{\gamma}_{p,q}$, $\nabla D_p  =  \nabla D_q$ on $\{ \hat{\gamma}_{p,q} \}^C$; and  outside $\gamma_{p,q}$ we have $\nabla D_p  \neq  c \nabla D_q$ for any $c \in \Bbb R$. 
\end{enumerate} 
\end{sub-lemma}
$\bf{Proof:}$
For the first part note the distance function's level sets are the spheres, so by Gauss's lemma its  integral curves are some re-parameterizations of geodesics from $p$.  
Now observe in geodesic polar coordinates  that the unit speed geodesics $\gamma$ satisfies $\dot{\gamma} = G_{\star} \left( \frac{\partial}{\partial r} \right)$, and 
\[ < \nabla D_p, \dot{\gamma}> = dD_p\left( G_{\star} \left( \frac{\partial}{\partial r} \right) \right) = G^{\star}(dD_p)  \left( \frac{\partial}{\partial r} \right) = dr \left( \frac{\partial}{\partial r} \right) = 1. \]
So indeed the integral curve of $\nabla D_p$ are precisely the unit speed geodesics.  

The equalities in the second part follows immediately from the first part and the fact $\gamma_{p,q}$ is a geodesic.
   
 To prove the last piece of the second part assume at some $z \in  B_{\frac{i}{2}}(p) \bigcap B_{\frac{i}{2}}(q)$ we have  $\nabla D_p  =  c\nabla D_q$.  First note from the above we have $c = \pm 1$. There are two cases, first we'll deal with $c =1$.   Since the geodesics satisfy a second order O.D.E they are uniquely determined by their position and tangent vector, so when $c = 1$ we have  both the geodesic from $p$ and the geodesic form $q$ are the same curves.  Without loss of generality $p$ is further away than $q$ and this point lies along the same minimal length geodesic (of length less than $\frac{i}{2}$) which connects $p$ and $q$, i.e. $\gamma_{p,q}$.  In the case  $c = -1$ we can follow the geodesic form $p$ to the point and then from the point back to $q$ forming a geodesic of length less than $i$ - which then must by lemma $\ref{normal}$ be the unique such one, i.e. $\gamma_{p,q}$.

\qed

Back to the relation  $<\nabla D_p  - \nabla D_q, \dot{z}> = 0$.  The first observation is that we can express a solution of this relation via a vector field. Using $\frac{\pi}{2}$ rotation field $\Theta$ (from section 1.1) we see that the solution to the differential relation $<\nabla D_p  - \nabla D_q, \dot{z}> = 0$ are re-parameterizations of  integral curves of the vector field $ \Theta (\nabla D_p  - \nabla D_q) $.  
Fortunately, as with $\nabla D_p$, we can say quite a bit about this vector field.  We are most interested in its integral curve corresponding to  $D_p - D_q = 0$.  Let $c_{p,q}(t)$ be the connected component of the integral curve of $ \Theta (\nabla D_p  - \nabla D_q) $  passing through $m$ in any set of interest to us;  and assume its parameterization   satisfies $c_{p,q}(0) = m$.  

\begin{lemma}[The Distance Difference Gradient Lemma]
\label{grad}
Assume $d(p,q) < \frac{i}{4}$ then:
\begin{enumerate}
\item
 In $B_{\frac{i}{4}}(p)$ we have that $c_{p,q}(t)$ is the unique component of $D_p - D_q = 0$, and  $m$ is the unique point of $\gamma_{p,q}$ on $c_{p,q}(t)$.
\item
 In $B_{\frac{i}{4}}(p) \bigcup  B_{\frac{i}{4}}(q)$ we have $d(c_{p,q}(t),p)$ and $d(c_{p,q}(t),q) $ strictly increase as the parameter  $|t|$ increases. 
\item
No geodesic from $p$ in  $B_{\frac{i}{4}}(p)$ is tangent to  $c_{p,q}(t)$ or cuts trough $c_{p,q}(t)$ twice on the same side of $\gamma_{p,q}(t)$. 
\item
In $B_{i}(o) - o$ it can never be the case that  $\nabla D_o  - \nabla D_p = c (\nabla D_o  - \nabla D_q) \neq 0$ at a point on  $c_{p,q}(t)$. 
\end{enumerate}
\end{lemma}
$\bf{Proof:}$
First for the uniqueness of $m$: suppose a point $l \neq m$ is $\gamma_{pq}$.  Then $l$ is within $\frac{i}{4}$ of $p$- hence 
$d(p,l)$ is determined by the length of the segment of  $\gamma_{p,q}$ form $p$ to $l$, similarly for $q$ (using $\frac{i}{2}$).  Now note that as we move from $m$ toward, say, $p$ that $D_p$ decreases while $D_q$ increase - so $D_p - D_q \neq 0$ at another point of $\hat{\gamma}_{p,q}$.  When on $\{ \hat{\gamma}_{p,q} \}^c$ say above $p$ the segment of $\gamma_{p,q}$ from $q$ to $l$ in fact covers the shorter segment from $p$ to $l$ - forcing  $D_p - D_q \neq 0$ once again.  So $l$ cannot satisfy $D_p - D_q =0$, forcing $m$ to indeed be the unique point of $\gamma_{p,q}$ on $D_p - D_q = 0$.

The remainder of the first part and the second part are intimately related.  To see why we first look at the component of $c_{p,q}(t)$ in $B_{\frac{i}{4}}(p)$ and note any component of  $D_p - D_q = 0$  would have to have a point closest to $p$.   
 This closest point  is tangent to a sphere emanating form $p$.  The same sort of phenomena must take place for the distance function to have a critical point; namely 
  if a point $z$ along any integral curve  of $\Theta (\nabla D_p  - \nabla D_q)$ is  a critical  point of the distance function $D(p, \cdot)$ then either $\nabla D_p  =  \nabla D_q$ or a circle is tangent to the solution curve. In the tangent case  $\Theta(\nabla D_p  - \nabla D_q) = c \Theta \nabla D_p$, or rather $\nabla D_p  - \nabla D_q = c  \nabla D_p$;  so both these situation we have forced the case  $\nabla D_p  =c\nabla D_q$.  The triangle inequality  tells $B_{\frac{i}{4}}(p) \subset B_{\frac{i}{2}}(p) \bigcap B_{\frac{i}{2}}(q)$ so we may  use sub-lemma $\ref{silly}$ to note that the point where this occurs is on  $\gamma_{p,q}$; but from above to be on $D_p - D_q = 0$ and  $\gamma_{p,q}$ means you must be exactly $m$.  So we have both that every  component of   $D_p - D_q = 0$ in $B_{\frac{i}{4}}(p)$ contains $m$ , and that the distance to $p$ parameterized by $t$ can have no critical points except at $m$ (similarly for $q$).

For the third part note that if a geodesic cuts twice on the same side of $\gamma_{p,q}(t)$ then some other geodesic must be tangent (via the mean value theorem in geodesic polar coordinates - see the see figure $\ref{geo}$).
 \begin{figure*}
\vspace{.01in}
\hspace*{\fill}
\epsfysize = 1.5in 
\epsfbox{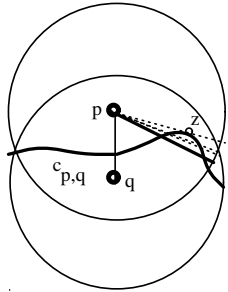}
\hspace*{\fill}
\vspace{.01in}
\caption{\label{geo}The Hunt for the Tangent}
\end{figure*}
%{ \bf (I'm pretty happy with the ``then'' in this argument but if one was in the mood to further justify it  one should  look at the picture in geodesic  polar coordinates.
%Now note that that from the above that $c_{pq}(t)$ is monotonically increasing with $r$, stays away form $\theta =0$ and $\theta =\pi$ (since  since $|D_p -D_q| > d(p,q)$ there, and hits $m$ with $\theta = \frac{\pi}{2}$ (since $\Theta(\nabla D_p  - \nabla D_q) = 2 \Theta(  \nabla D_p)$ at that point).
%So $c_{pq}$ in polar coordinates is a differentiable curve $f(r)$ as in figure $\ref{pole}$.
% \begin{figure*}
%\vspace{.01in}
%\hspace*{\fill}
%\epsfysize = 2in 
%\epsfbox{pole.eps}
%\hspace*{\fill}
%\vspace{.01in}
%\caption{\label{pole}The Hunt for the Tangent}
%\end{figure*}
 % Cutting through $c_{pq}$ twice with a geodesic exactly corresponds to cutting $f(r)$ twice at say  $r_1 > r_0$.  So by the mean value theorem there is a $r_2 \in [r_0,r_1]$ with the property $f(r_0) - f(r_1) = 0 = \frac{d f}{dr}(r_2) (r_1 -r_0)$.  Hence  $\frac{d f}{dr}(r_2) = 0$ or the corresponding $\theta = f(r_2)$ is tangent at this point, as needed.  }
So we are reduced to the tangent case. Let the point where this tangency occurs be called $z$. This tangency implies  $\Theta (\nabla D_p  - \nabla D_q) = c \nabla D_p $, or rather     $\Theta \nabla D_q = \Theta \nabla D_p -c \nabla D_p$.  But  $\nabla D_p$ and  $\Theta \nabla D_p$ are orthogonal, with  $\Theta \nabla D_p$ also  of unit length (by sub-lemma $\ref{silly}$) - so taking norms we have   $1 = \sqrt{1 + c^2}$; forcing $c=0$.  So at such a point $\Theta(\nabla D_p  - \nabla D_q) = 0$ or $\nabla D_p  = \nabla D_q$.  Once again using   $B_{\frac{i}{4}}(p) \subset B_{\frac{i}{2}}(p) \bigcap B_{\frac{i}{2}}(q)$ the above sub-lemma applies and forces  $z$ onto  $\{ \hat{\gamma}_{p,q} \}^C$; and 
now part one kicks in to eliminate this possibility of $c_{p,q}$ hitting  $\{ \hat{\gamma}_{p,q} \}^C$.

% \begin{figure*}
%\vspace{.01in}
%\hspace*{\fill}
%\epsfysize = 2in 
%\epsfbox{proj.eps}
%\hspace*{\fill}
%\vspace{.01in}
%\caption{\label{proj} An Cool Fact}
%\end{figure*}

 To prove the last part we will see first if  $\nabla D_o  - \nabla D_p = c (\nabla D_o  - \nabla D_q) \neq 0$ at $z$, then  $c=1$ and $ \nabla D_p =  \nabla D_q$ at $z$. 
 This because the $\nabla D_i$ are all unit vectors (by sub-lemma $\ref{silly}$ part one); and among unit vectors $v$ and $w \neq v$ the direction of $v - w$ is uniquely determined by $w$  (since we are re-parameterizing the radius one circle at located a $v$). Recall from the above sub-lemma that $ \nabla D_p =  \nabla D_q$ at $z$ implies  $z \in \{ \hat{\gamma}_{p,q} \}^C$. As above, no point of $c_{p,q}(t)$ can be on  $ \{ \hat{\gamma}_{p,q} \}^C$, as needed.

\qed

The last set of issues concerns not the $c_{p,q}$ curves themselves but the home of the curves.  Suppose we have   three points $p$, $q$ and $o$ on the boundary of circles of radius $r$  less than $\frac{i}{6}$; then using $\gamma_{p,o}$ and $\gamma_{q,o}$ we can split   $B_{\frac{i}{2}}(o)$ into four cones (as in figure $\ref{cones}$).  Let the forward length cone (FLC) be the cone containing $\hat{\gamma}_{p,q}$ , the backward length cone (BLC) be the forward length cone's mirror image (see figure $\ref{cones}$ and the following lemma to see this is well defined).  One nice thing about this cone notion is that it eliminates certain regions where the center of a small circle might have wanted to live - and it is at this point where we find the only place in the argument where using the strong convexity radius is necessary.  We shall explore these facts (and clear up an issue in the introduction) with the following  corollary to lemma $\ref{normal}$ in section 1.2:

\begin{figure*}
\vspace{.01in}
\hspace*{\fill}
\epsfysize = 2in 
\epsfbox{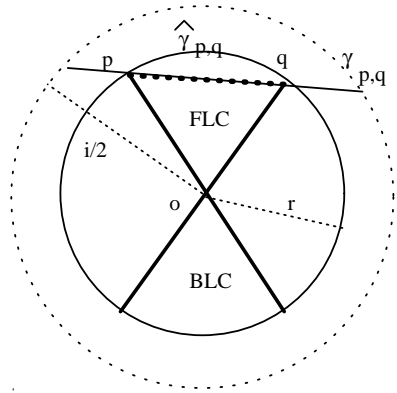}
\hspace*{\fill}
\vspace{.01in}
\caption{\label{cones}The Length Cones}
\end{figure*}

\begin{lemma}[The Length Cone Lemma] \label{cones1}
 Given three points $p$, $q$ and $o$ on the boundary of a circle of radius $r$ less than $\frac{i}{6}$ centered at $c$: 

\begin{enumerate}

\item
The forward length cone is well defined.
\item
The intersection of the three forward length cones is a triangle contained in $B_{\frac{i}{2}}(c)$; which, incidentally, can be parameterized in a way compatible with the abstract gluing to form $|K_{\{p_1, \dots , p_n\}}|$ (see figure  $\ref{tri}$).
\item
The $B_{\frac{i}{2}}(c)$ is dissected by the curves $\{ \gamma_{p,q} ,   \gamma_{p,o},   \gamma_{q,o}\}$  into seven regions - the unique FLC, three distinct BLCs, and the regions separating the BLCs (see in figure  $\ref{tri}$). 
\item
If in addition the balls radius is less than $\tau$, then circle's center cannot be in any of the backward length cones, and the simplex is contained in $B_{r}(c)$.

\end{enumerate}
\end{lemma}
 {\bf Proof:}  In the course of the this proof I use the triangle inequality in the sense of  $p \in B_{a}(q)$ then  $ B_{b}(p) \subset B_{a + b}(q)$.  I will simply call it the triangle inequality.

 First we attack well definedness of the FLC.  Since the radius of the circle is less than $\frac{i}{6}$ we have  by the triangle inequality
$p$ and $q$ are in $B_{\frac{2i}{6}}(o)$, and that $d(p,q) < \frac{2i}{6}$.  Now note $\hat{\gamma}_{p.q}$ is contained in    $B_{\frac{3i}{6}}(o)$ - since if it left $B_{\frac{3i}{6}}(o)$  then by the triangle inequality it would be a curve from $p$ to $q$ leaving  $B_{\frac{i}{6}}(p)$ and then afterwards entering  $B_{\frac{i}{6}}(q)$, so  would have length greater than $\frac{2i}{6}$; contradicting $d(p,q) < \frac{2 i}{6}$.  We would now like this geodesic segment to stays in one cone section.  If it did  not then it would cut a bounding geodesic, $\gamma_{o,p}$ or $\gamma_{o,q}$, in at least two points on the cone's boundary. Now the diameter geodesics are of length less than $i$, so between $c_1$ and $c_2$ we now have contradictory distinct geodesics of length less than $i$ connecting them.  So the cone is well defined.

 To get started on the second part first simply parameterize the needed simplex. Note we may choose  $\theta_1$ and  $\theta_2$ such that $\gamma_{o,p}(t) = exp_{o}(t v(\theta_1))$ and $\gamma_{o,q}(t) = exp_{o}(t v(\theta_2))$ bound the FLC, and by  the first argument  $\hat{\gamma}_{p,q} \subset B_{\frac{3 i}{6}}(o)$ so any point on $\hat{\gamma}_{p,q}$ has its shortest length geodesics from $o$ described by $exp_{o}(t v(\theta))$ with   $\theta \in [\theta_1, \theta_2]$.   The fact there are no double intersections of short geodesics now demonstrates that each  $exp_{o}(t v(\theta))$ with  $\theta \in [\theta_1, \theta_2]$ hits exactly one point on $\hat{\gamma}_{p,q}$, so we can parameterize  $\hat{\gamma}_{p,q}$ as  $exp_{o}(r(\theta) v(\theta))$ - with $r(\theta)$ continuous.
Now we have a homeomorphism  (once again since there are not double intersection between small geodesics) of the region trapped between the geodesics with the wedge  $\theta \in [\theta_1, \theta_2]$ and $ r \leq r(\theta)$.  This wedge is clearly homeomorphic to the disk.  However to parameterize the triangle we would like a parameterization which is compatible with the abstract gluing, allowing  $R: |K_{\{p_1, \dots , p_n  \}}| \rightarrow M$ to be continuous. To accomplish this first parameterize the edges with the unit interval using a  $d(p,q)$ speed parameterization of the geodesic.
Now for each triangle in the surface we will map on an  equilateral triangle by  first mapping the boundary of our triangle to the boundary of our simplex $M$ such that the linear unit length sides are mapped onto the geodesics at a speed of $d(p,q)$.
We have just found  that the image is in fact itself the boundary of a  disk, so  this mapping can be extended from a homeomorphism  of the equilateral triangle's boundary to a homeomorphism of the equilateral triangle onto the triangle in $M$.  So by its very construction the mapping on the boundary of the equilateral triangle is affine related to the mapping of the interval, and we indeed have the needed compatible parameterizations.  

We also need for this second part that  this is the same simplex for all three of the points, but by the very definition of the FLC with respect to a point being the cone section containing the opposite $\hat{\gamma}$ we have that the FLC must contain this simplex.  So the simplex is in the intersection - and precisely the intersection since the respective $\gamma_{p,q}$ curves all cut their respective $B_{i}(p)$ balls in to two distinct sides - and by the triangle inequality this entire discussion is taking place in any one of the radius $i$ balls. 

Lastly for the second part we need the  simplex is contained in $B_{\frac{i}{2}}(c)$ ball.  To see this note from above that  every point, $k$, in the simplex  is on a geodesic, from say $o$, of length less than $\frac{3i}{6}$ which hits $\hat{\gamma}_{p,q}$.  If $k$ is within $\frac{2 i}{6}$ of $o$ then $k$  is with in $\frac{i}{2}$ of $c$.  If not it is within $\frac{i}{6}$ of $\hat{\gamma}_{q,o}$.  The geodesic  segment  $\hat{\gamma}_{q,o}$ has  length less than $\frac{2i}{6}$, so the point where the geodesic from $k$ hits $\hat{\gamma}_{q,o}$ is within  $\frac{i}{6}$ of either $p$ or $q$; so $k$ is once again with in $\frac{3i}{6}$ of $c$ - as needed.

\begin{figure*}
\vspace{.01in}
\hspace*{\fill}
\epsfysize = 2in 
\epsfbox{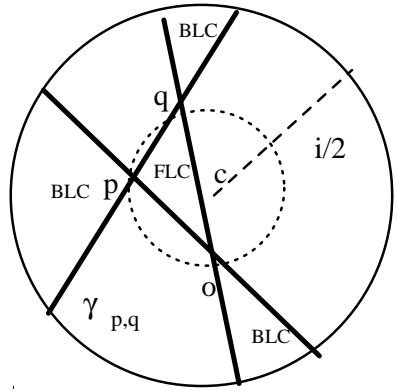}
\hspace*{\fill}
\vspace{.01in}
\caption{\label{tri}The Seven Regions}
\end{figure*}

For the third part note that by the triangle inequality a geodesic between two points in $B_{\frac{i}{2}}(c)$ is contained in the ball must be of length less than $i$ (use the sense of the triangle inequality introducing this proof at a point were the geodesic first hits $S_{\frac{i}{2}}(c)$).  Now extend  all the geodesic segments bounding the simplex out to the boundary in $B_{\frac{i}{2}}(c)$; and note, having length less than $i$  these segments   can intersect each other only once - at the vertices of the simplex (by part 2 above the simplex is in there).  So noting the definitions of all of our objects we indeed end up decomposing $U$ into the seven components pictured in the figure.

%\begin{figure*}
%\vspace{.01in}
%\hspace*{\fill}
%\epsfysize = 2in 
%\epsfbox{cut.eps}
%\hspace*{\fill}
%\vspace{.01in}
%\caption{\label{cut} A Typical Contradiction}
%\end{figure*}

 For the final part, the fact that the  simplex is contained in $B_{r}(c)$ is a trivial consequences of the $\hat{\gamma}$ curves having interior in this ball and never intersecting each other with in the ball.  
To handle  the center's placement, assume the contrary that the  center,$c$, is in the backward  length cone with respect to say $o$.  Then take the shortest length geodesic from $c$ to $o$ and continue it through the forward length cone until it hits $\hat{\gamma}_{p,q}$ at $h$.
Note this is a geodesic of length less than $\frac{4 i}{6}$, so in particular it is the unique shortest length geodesic between these points.  By assumption $p$ and $q$ are on the boundary of a disk at $c$ with radius less than $\tau$,
so by strong convexity $\hat{\gamma}_{p,q} \subset B_{r}(c)$.  In particular the shortest length geodesic to $h$ never leaves  $B_{r}(c)$, contradicting the fact it hit $o$ before $h$.  So the center must not be in the backward length cone.

\qed

We are finally in a position to prove the small disk uniqueness  theorem form the introduction (theorem $\ref{lit}$).

{\bf Proof of Theorem \ref{lit}:}
 As above call the points in the triple $p$, $q$ and $o$, and recall $\delta < min\{ \frac{i}{6}, \tau  \}$.   The whole discussion takes place with $B_{\frac{i}{2}} (c)$, since $d(c_{p,q}(t),c) \leq d(c_{p,q}(t), p) + d(c, p) \leq \frac{i}{3}$.  In particular the whole discussion takes place at points where we can use lemmas $\ref{grad}$ and $\ref{cones}$.  Note  a circle of radius less than $\frac{i}{6}$ passing through all three of these points corresponds to an intersection of   $c_{p,q}$ and $c_{p,o}$ (by part one of lemma  $\ref{grad}$). Further note  if  $c_{p,q}$ and $c_{p,o}$ intersect at $z$ then so does $c_{q,o}$, since then $ d(z,p)=  d(z,q) =  d(z,o)$.

    At these intersections we would like that the  curves must cut through each other; i.e. not be tangent.   If a pair were tangent then there are two possibilities: one where neither of the tangent curves is stopped when thought of as a parameterized  integral curve, and the other when one has stopped ($\dot{c}_{p,q} = 0$). If neither has stopped then  (up to index permutation) $\Theta(\nabla D_p  - \nabla D_q) = c \Theta(\nabla D_p  - \nabla D_o)\neq 0$, or rather  $(\nabla D_p  - \nabla D_q) = c (\nabla D_p  - \nabla D_o) \neq 0$.  This case implies by lemma $\ref{silly}$ that we are on   $\{ \hat{\gamma}_{q,o} \}^C $.   If one is stopped then    (once again up to index permutation) $\nabla D_p  - \nabla D_q =0$.   So by sub-lemma $\ref{silly}$  the point of intersection is on  $ \{ \hat{\gamma}_{p,q} \}^C $.  Be warned - I will refer to this same case analysis later in the argument. To finish off this impossibility argument simply note that all the curves intersect at this point so (up to index permutation) $ \{ \hat{\gamma}_{p,q} \}^C $ and $c_{p,q}$ would intersect. But by part one of lemma $\ref{grad}$ this point must then be $m_{p,q}$ which is not on   $ \{ \hat{\gamma}_{p,q} \}^C $.   So the intersection must include no tangencies.

Now suppose we had two distinct circles of radius less then $\frac{i}{6}$ passing through these triples.  This non-uniqueness of circles would correspond to at least a pair  of intersections of  $c_{p,q}$, $c_{p,o}$, and $c_{q,o}$ - at $z_1$ and  $z_2$. Note without loss of generality there are no further intersections of  any pair of $c$ curves between $z_1$ and $z_2$; 
since intersections occur only as triples at isolated possible distances by the above no tangency result.
Hence we are forced to have a picture something very much like figure $\ref{right}$, at least in the sense that we have one of the curves trapped between the other two; without loss of generality assume this curve is $c_{q,o}$.

%\begin{figure*}
%\vspace{.01in}
%\hspace*{\fill}
%\epsfysize = 2in 
%\epsfbox{con1.eps}
%\hspace*{\fill}
%\vspace{.01in}
%\caption{\label{mult} Multiple centers}
%\end{figure*}
\begin{figure*}
\vspace{.01in}
\hspace*{\fill}
\epsfysize = 1.75in 
\epsfbox{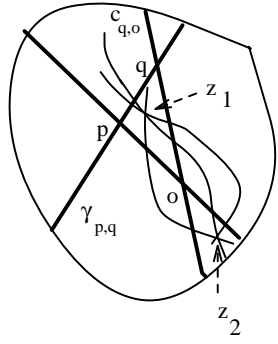}
\hspace*{\fill}
\vspace{.01in}
\caption{\label{right} The First Possibility}
\end{figure*}

Note    $\Theta (\nabla D_q  - \nabla D_o) = \Theta (\nabla D_p  - \nabla D_o) + \Theta (\nabla D_q  - \nabla D_p)$ , so $ \Theta (\nabla D_p  - \nabla D_o)$ and $ \Theta (\nabla D_q - \nabla D_p)$ are always on opposite sides of  $\Theta (\nabla D_p  - \nabla D_o)$; in the sense that the $\nabla D_q  - \nabla D_o$ component of these vectors must have opposite signs.  Now since the curves must cut through each other at $z_1$ and $z_2$ and no where in between,  the side of $c_{q,o}$ which  $ \Theta (\nabla D_p  - \nabla D_o)$ and $ \Theta (\nabla D_q  - \nabla D_p)$ live on must  switch roles while one moves along $c_{q,o}$ from $z_1$ to $z_2$. In order for these roles to switch and the vectors to always be on opposite sides,  we are forced  (as with our no tangency result above) to have at some point along $c_{q,o}$ either  $\Theta (\nabla D_p  - \nabla D_o) = c \Theta (\nabla D_p  - \nabla D_q) \neq 0 $, $\nabla D_p  - \nabla D_o = 0$, or $\nabla D_p  - \nabla D_q = 0$. Immediately,  the last part of lemma $\ref{grad}$ kicks in eliminating the first case. 

   The remaining case (up to index permutation) is were  $\nabla D_p  = \nabla D_o $  - which means this point lies on  $\{\hat{\gamma}_{p,o} \}^c$ (by the sub lemma).   We will be justifying each of the two pictured cases  - depending on which component of $\{ \hat{\gamma}_{p,o} \}^c$ the point lies (figure $\ref{right}$ and figure $\ref{left}$).  The idea in each case is this same - we will use this condition to force one of the centers $z_i$ into a backward length cone with respect to either  $o$ or $p$, and then use the previous lemma to observe this is impossible.   

Observe  by the first and third parts of lemma $\ref{grad}$ that $c_{q,o}(t)$  cannot be tangent to or cut twice any geodesic in the picture; and $c_{q,o}(t)$ cannot sneak around a geodesic in the picture (since the  decomposition in the third part of lemma $\ref{cones}$ has each of the geodesics heading all the way  out to $S_{\frac{i}{2}}(c)$).  So in the first case the $c_{q,o}$ is stuck to the right of $\gamma_{q,o}$ and must cross $\gamma_{p,o}$ to the right and can never cross it again.  So one of the $z_i$ is  stuck in the backwards length cone of $o$ as needed.  The other case is similar to this except we must justify why we must come from above and cross into the backward length cone with respect to $p$.  This is because  to come from below would mean to cut across $\gamma_{p,o}$ first to get there producing a double cut of $c_{q,o}$ by $\gamma_{p,o}$ to left of  $\gamma_{q,o}$,  contradicting part 3 of  lemma $\ref{grad}$.  So indeed once again we are forced into the contradictory setting of having a center in a backward length, this time with respect to $p$.

So the points $z_1$ and $z_2$ cannot simultaneously exist - and there is no non-uniqueness among small circles.

% and monotonically increasing distance (by lemma $\ref{silly}$.  
%(say good!  prev and beyond) 
%(sense of between from monotinicty) this because since we know the curves cannot be tangent at these points  (and are differentiable)  
 %since any intersection correspond to one of the monotonically increasing radii
 % (by lemma $\ref{grad}$) there must be no other intersection between any of these curves occurring between $z_1$ and $z_2$. 

\begin{figure*}
\vspace{.01in}
\hspace*{\fill}
\epsfysize = 1.75in 
\epsfbox{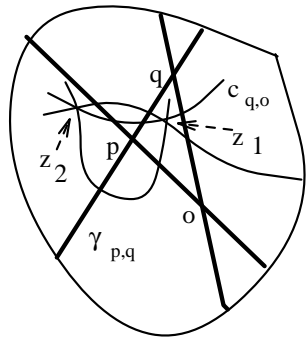}
\hspace*{\fill}
\vspace{.01in}
\caption{\label{left} The Second Possibility}
\end{figure*}

\qed

\end{chapter}
\bibliographystyle{plain}
\bibliography{foo}
\end{document}